\documentclass[12pt]{amsart}
\usepackage{a4wide,amsfonts,amssymb,stmaryrd,amscd,amsmath,latexsym,amsbsy}
\numberwithin{equation}{section}
\theoremstyle{plain}
\newtheorem{thm}{Theorem}[section]

\newtheorem{prop}[thm]{Proposition}

\newtheorem{defi}[thm]{Definition}
\newtheorem{lem}[thm]{Lemma}
\newtheorem{cor}[thm]{Corollary}

\theoremstyle{remark}
\newtheorem{rema}[thm]{Remark}

\newcommand{\solu}[1]{\begin{sol}{\bf (\ref{#1})}}

\newcommand{\AAA}{\mathcal A}

\newcommand{\CC}{\mathbb C}
\newcommand{\DD}{\Delta}

\newcommand{\gt}{\mathfrak g}
\newcommand{\ug}{U_q(\mathfrak g)}

\title[Spherical functions and Macdonald-Koornwinder polynomials]
{Vector valued spherical functions and Macdonald-Koornwinder
polynomials}

\author{Alexei A. Oblomkov}
\address{A. A. Oblomkov,
Department of Mathematics, MIT, 77 Mass. Ave, Cambridge, MA 02139}
\email{oblomkov@math.mit.edu}
\author{Jasper V. Stokman}
\address{J.V. Stokman,
KdV Institute for Mathematics, Universiteit van Amsterdam,
Plantage Muidergracht 24, 1018 TV Amsterdam, The Netherlands.}
\email{jstokman@science.uva.nl}

\begin{document}
\begin{abstract}
We interpret the five parameter family of Macdonald-Koornwinder polynomials
as vector valued spherical functions on quantum Grassmannians.
\end{abstract}

\maketitle

\section*{Introduction}
The representation theoretic construction of (quantum) conformal blocks in
certain conformal field theories is closely related to
harmonic analysis on (quantum) symmetric spaces of group type,
see, e.g., \cite{FR} and \cite{EV}. A striking consequence
is the interpretation of $A$-type Macdonald
polynomials as vector valued spherical functions for
the quantum analogue of the
symmetric pair ($\hbox{U}(n)\times \hbox{U}(n),
\hbox{diag}(\hbox{U}(n))$) of group type on the one hand
(see \cite{EK}), and as quantum conformal blocks on
the other hand (see \cite{EV}). With these interpretations
many properties of $A$-type Macdonald
polynomials, such as the
Macdonald-Ruijsenaars difference equations, quantum Khnizhnik-Zamolodchikov
equations, dualities and orthogonality relations, obtain their natural
representation theoretic and conformal field theoretic interpretations.
In this paper we consider the harmonic analytic part of these
constructions for the quantum
analogues of the symmetric pair $(U,K)=(\hbox{U}(2n),
\hbox{U}(n)\times \hbox{U}(n))$. This leads to the interpretation of
the five parameter family of Macdonald-Koornwinder polynomials as
vector valued spherical functions.

In \cite{O} the classical analogue of our main result was established.
It yields the interpretation of
$BC$-type Heckman-Opdam polynomials
as the restriction to the maximal torus of regular vector valued functions
\[
f: U\rightarrow \mathbb{C}\hbox{det}^{-\kappa_1}\otimes
S^{n\kappa}(\mathbb{C}^n)\hbox{det}^{\kappa_1-\kappa},\qquad
\kappa_1\in\mathbb{Z},\, \kappa\in\mathbb{Z}_{\geq 0}
\]
(with $S^{n\kappa}(\mathbb{C}^n)$ the homogeneous polynomials
of degree $n\kappa$) which transform under the
left (respectively right) regular $K$-action on $U$ according to the
natural $K$-action on the image space (respectively
the $K$-character $\hbox{det}^{-\kappa_2}\otimes
\hbox{det}^{\kappa_2}$ for some $\kappa_2\in\mathbb{Z}$).

In this paper we define a continuous one-parameter
family of quantum analogues of the symmetric pair $(U,K)$, following closely
Letzter's \cite{L}--\cite{L4} approach of constructing
quantum symmetric pairs as coideal subalgebras.
One may as well view this family of quantum symmetric pairs as
a continuous one-parameter family of quantum analogues of the
complex Grassmannian $U/K$. We proceed by defining the analogue
of the $K$-representation $\mathbb{C}\hbox{det}^{-\kappa_1}\otimes
S^{n\kappa}(\mathbb{C}^n)\hbox{det}^{\kappa_1-\kappa}$ for the
associated coideal subalgebras. This allows us to define
vector valued spherical functions for the one-parameter family
of quantum symmetric pairs in essentially the same manner as in the classical
case \cite{O}, with the exception that we have now the additional freedom to
choose different coideal subalgebras for the transformation behaviour
under the left and right regular action, respectively. We relate
the resulting vector valued spherical functions, which now depend on two
continuous and three discrete parameters, to the five
parameter family of Macdonald-Koornwinder polynomials.
The Macdonald-Koornwinder polynomials are Koornwinder's \cite{K} extension
of the $BC$-type Macdonald polynomials that contain all Macdonald polynomials
of classical type as special cases.

The interpretation of the five parameter family of
Macdonald-Koornwinder polynomials as vector valued
spherical functions contains several known results as special cases.
It entails the interpretation of
a two parameter subfamily of Macdonald-Koornwinder polynomials
as zonal spherical functions on quantum Grassmannians,
which was established by Noumi, Sugitani and Dijkhuizen \cite{NDS}
(see also \cite{DS}). This special case plays in fact an essential role
in establishing our general result. In rank one we reobtain
the interpretation of the four parameter
family of Askey-Wilson polynomials as matrix coefficients of quantum
$\hbox{sl}(2)$ representations, established before by Koornwinder \cite{Koo0}
(zonal case), and by Noumi and Mimachi \cite{NM}
and Koelink \cite{Koe} (general case).

The content of the paper is as follows. In \S 1 we give the main
definitions and formulate the main result. In \S 2 we define the
notion of expectation value for the quantum symmetric pairs under
investigation and establish branching rules using deformation theory.
In \S 3 we prove
the zonal case by translating the main results of \cite{NDS} to
our setup. In \S 4 we generalize the Chevalley restriction theorem
to the setup of vector valued spherical functions. Its
description involves the vector valued spherical function of the
smallest degree, which we call the ground state (in the zonal case
the ground state is the unit). The restriction of the ground state
to the quantum torus is computed explicitly in \S 5. In \S 6 we
establish the quantum Schur orthogonality relations for the vector
valued spherical functions and, combined with the results of
previous sections, we establish the explicit interpretation of
the Macdonald-Koornwinder polynomials as vector valued spherical functions.

The construction of expectation values in \S 2 and the dynamical
quantum group interpretation of the rank one results in \cite{St2}
hint at a natural interpretation of the
Macdonald-Koornwinder polynomials as quantum conformal blocks.
A more detailed study in this direction is subject of future research.

\vspace{1cm}

{\it Acknowledgments:} A substantial part of the research was done when
the second author visited MIT for a period of three months in the
beginning of 2002 and when the first author visited the KdV
institute in the summer of 2002. Both
authors are grateful for the hospitality of the institutes. Both
authors thank Pavel Etingof for stimulating discussions. The
second author is supported by the Royal Netherlands Academy of
Arts and Sciences (KNAW). The first author is supported by the
NSF grant DMS-9988796.


\section{Formulation of the main result}
In this section we give
definitions of the Macdonald-Koornwinder polynomials, the
quantum symmetric pairs and the associated class of vector valued
spherical functions, and we
formulate the main result of the paper. We add insightful proofs of some of the
intermediate results, but we postpone the more technical parts to later
sections. We fix a positive integer
$n\geq 1$ and a deformation parameter $0<q<1$ throughout the paper.


\subsection{Macdonald-Koornwinder polynomials.}
Koornwinder \cite{K} extended the definition of
the Macdonald polynomials \cite{M}
associated to the non-reduced, irreducible
root system $BC_n$ to a family of orthogonal polynomials depending,
besides on the deformation parameter $q$, on five additional coupling
parameters. This family of polynomials, known nowadays as
Macdonald-Koornwinder polynomials, reduces for $n=1$ to the
celebrated, four parameter family of Askey-Wilson polynomials.
In this subsection we recall their definition.

Denote $\Lambda_n$ for the lattice $\mathbb{Z}^n$ and let
$\Lambda_n^+\subset \Lambda_n$ be the set of partitions of length
$\leq n$. The dominance partial order on $\Lambda_n$ is defined by
\[\lambda\leq \mu \qquad \Leftrightarrow \qquad
\sum_{i=1}^j\lambda_i\leq \sum_{i=1}^j\mu_i,\quad j=1,\ldots,n.
\]
The Weyl group $W=S_n\ltimes \{\pm 1\}^n$ of $BC_n$, where $S_n$
is the symmetric group in $n$ letters, acts on $\Lambda_n$ by
permutations and sign changes. Each $W$-orbit in $\Lambda_n$
intersect $\Lambda_n^+$ exactly once.

Let $\CC[u^{\pm 1}]=\CC[u_1^{\pm 1},\dots,u_n^{\pm 1}]$ be the algebra
of Laurent polynomials in $n$ independent variables $u_k$ ($1\le k\le
n$), or, equivalently, the algebra of
regular functions on the complex $n$-torus $T^{\CC}=(\CC^*)^n$.
A basis of $\mathbb{C}[u^{\pm 1}]$ is given by the monomials
$u^\lambda=u_1^{\lambda_1}u_2^{\lambda_2}\cdots u_n^{\lambda_n}$
($\lambda=(\lambda_1,\lambda_2,\ldots,\lambda_n)\in \Lambda_n$).

The Weyl group $W$ acts on $\CC[u^{\pm 1}]$
by permutations and inversions of the $u_k$.
Let $\CC[u^{\pm 1}]^W\subset \CC[u^{\pm 1}]$ be the subalgebra of
$W$-invariant Laurent polynomials.
The orbit sums $m_\lambda=
\sum_{\mu\in W\lambda}u^{\mu}$ ($\lambda\in \Lambda_n^+$) form a linear
basis of $\mathbb{C}[u^{\pm 1}]^W$.

The Macdonald-Koornwinder polynomials form an orthogonal basis of
$\mathbb{C}[u^{\pm 1}]^W$ with respect
to a particular scalar product on $\mathbb{C}[u^{\pm 1}]^W$. The
scalar product, which we define now first, depends on five
additional coupling parameters $a,b,c,d$ and $t$.
It is defined in terms of an absolutely continuous
measure with respect to the normalized Haar measure
on the natural compact real form $T=\mathbb{T}^n$
of $T^{\CC}=(\mathbb{C}^*)^n$, where $\mathbb{T}$ is the unit
circle in the complex plane.
The corresponding weight function $\Delta$ is most conveniently expressed
in terms of the $q$-shifted factorial,
\[
(a;q)_k=\prod_{i=0}^{k-1}(1-a q^i),\qquad \forall k\in
\mathbb{Z}_+\cup \{\infty\}
\]
by $\Delta(u)=\Delta^+(u)\Delta^+(u^{-1})$ with
$\DD^+(u)=\DD^+(u;a,b,c,d;q,t)$ defined by
\begin{equation*}
\DD^+(u)=\prod_{i=1}^n\frac{(u_i^2;q)_\infty}{(au_i,bu_i,cu_i,du_i;q)_\infty}
\prod_{1\le i< j\le n}\frac{(u_i/u_j,u_iu_j;q)_\infty}
{(tu_i/u_j,tu_iu_j;q)_\infty}.
\end{equation*}
Here $(a_1,\dots,a_s;q)_k=\prod_{j=1}^s(a_j;q)_k$ is a short-hand
notation for products
of $q$-shifted factorials. If $a,b,c,d,t$ are real and
\begin{equation}\label{reg}
|a|,|b|,|c|,|d|<1,\quad 0<t<1,
\end{equation}
then $\DD(u)$ is a
positive continuous weight function on $T$. Under these assumptions, we
can define the Macdonald-Koornwinder polynomials as follows, see
\cite{K}.

\begin{thm}\label{MKpolynomial}
There exist unique $W$-invariant Laurent polynomials
\[
P_\lambda(u)=P_\lambda(u;a,b,c,d;q,t)\in\mathbb{C}[u^{\pm 1}]^W,\qquad
\lambda\in\Lambda_n^+
\]
satisfying the two conditions
\begin{equation*}
\begin{split}
P_\lambda(u)&=m_\lambda+\sum_{\mu\in \Lambda_n^+: \mu<\lambda}
c_{\lambda\mu}m_{\mu},\quad \hbox{ some }\,\,
c_{\lambda\mu}\in\mathbb{C},\\
\int_T&P_\lambda(u){\overline{P_\mu(u)}}
\Delta(u)\frac{du}{u}=0,\quad \hbox{ if }\,\,
\lambda\not=\mu,
\end{split}
\end{equation*}
where
$\frac{du}{u}=\frac{du_1}{u_1}\frac{du_2}{u_2}\cdots \frac{du_n}{u_n}$
is the Haar measure on $T$.
We call $P_\lambda(u)$ the monic Macdonald-Koornwinder polynomial
of degree $\lambda\in\Lambda_n^+$.
\end{thm}
The theorem does not follow by a straightforward Gram-Schmidt type procedure
since the dominance ordering is not a total ordening.
The key tool in proving the theorem is an explicit
self-adjoint difference operator which maps a symmetric monomial
$m_\lambda$ to a linear combination of symmetric monomials $m_\mu$
involving only degrees $\mu\leq\lambda$ (see \cite{K}).
In our set-up, the self-adjoint operator
arises as the radial part of the quadratic Casimir element
for the corresponding quantum symmetric pair (see \cite[Thm. 3.3]{NDS}).
\begin{rema}
Theorem \ref{MKpolynomial} is also valid for generic
values of the parameters $a,b,c,d$ outside the region
\eqref{reg} after deforming the compact torus $T$ in the definition of
the orthogonality relations in a suitable way.
In general one looses positivity of the weight function,
but for suitable values of the parameters (still violating
the condition that the modulus of all the
four parameters $a,b,c$ and $d$ is less than one), one can reobtain
a positive measure by shifting the deformed compact torus to $T$
while picking up residues, see \cite{St}.
\end{rema}

\subsection{The quantized universal enveloping algebra.}
We introduce here the notations for the quantized universal
enveloping algebra of $\mathfrak{g}\mathfrak{l}(m)$. For further
details and standard facts, we refer to \cite{NYM} and \cite{N}.

In notations below we surpress the dependence on $m\in\mathbb{Z}_{>0}$
as much as possible. For most of our applications, $m$ will be either
$2n$ (in which case we write $\gt$ for the
Lie algebra $\mathfrak g\mathfrak l(2n)$), or $m$ will be $n$ (in
which case we stick to the notation $\mathfrak{g}\mathfrak{l}(n)$). Let
$\delta_{i,j}$ be the usual Kronecker delta function ($=1$ if $i=j$,
and $=0$ otherwise). The quantized universal enveloping
algebra $U_q(\mathfrak{g}\mathfrak{l}(m))$ is the unital algebra over
$\mathbb{C}$ generated by $K_i^{\pm 1}$ $(i=1,\dots,m)$, $x_j,y_j$
$(j=1,\dots, m-1)$, subject to the relations
\begin{equation*}\label{firstqg}
\begin{split}
&K_iK_j=K_jK_j, \qquad\qquad\qquad\qquad\,\,\,\, K_iK_i^{-1}=1=K_i^{-1}K_i,\\
&K_ix_jK_i^{-1}=q^{\delta_{i,j}- \delta_{i,j+1}}x_j,\qquad\qquad
K_iy_jK_i^{-1}=q^{-\delta_{i,j}+\delta_{i,j+1}} y_j,\\
&x_iy_j-y_jx_i=\frac{K_iK_{i+1}^{-1}-K_i^{-1}K_{i+1}}{q-q^{-1}}\,
\delta_{i,j},
\\
&x_ix_j=x_jx_i\quad\mbox{ and }\quad y_iy_j=y_jy_i \,\,\,
\qquad\qquad\mbox{ when } |i-j|\ge 2,\\
&x_i^2x_j-(q+q^{-1})x_ix_jx_i+x_jx_i^2=0,\qquad\quad\mbox{ when } |i-j|=1,\\
&y_i^2y_j-(q+q^{-1})y_iy_jy_i+y_jy_i^2=0,
\,\,\,\qquad\quad\mbox{ when } |i-j|=1.
\end{split}
\end{equation*}
The quantized universal enveloping algebra $U_q(\mathfrak{g}\mathfrak{l}(m))$
is a Hopf $*$-algebra, with comultiplication
\begin{equation*}
\begin{split}
\DD(x_j)&=x_j\otimes 1+K_jK_{j+1}^{-1}\otimes x_j,\\
\DD(y_j)&=y_j\otimes K_j^{-1}K_{j+1}+1\otimes y_j,\\
\DD(K_i^{\pm 1})&=K_i^{\pm 1}\otimes K_i^{\pm 1},
\end{split}
\end{equation*}
counit
\[
\varepsilon(x_j)=\varepsilon(y_j)=0,\qquad\varepsilon(K_i^{\pm 1})=1,
\]
antipode
\begin{equation}\label{lastqg}
S(x_j)=-K_j^{-1}K_{j+1}x_j,\qquad
S(y_j)=-y_jK_jK_{j+1}^{-1}, \qquad S(K_i^{\pm 1})=K_i^{\mp 1}
\end{equation}
and $*$-structure
\begin{gather*}
x_j^*=q^{-1}y_jK_jK_{j+1}^{-1},\qquad
y_j^*=qK_j^{-1}K_{j+1}x_j,\qquad
(K_i^{\pm 1})^*=K_i^{\pm 1}.
\end{gather*}
This $*$-structure corresponds classically
to choosing the skew-hermitean matrices $\mathfrak{u}(m)$
as (compact) real form of $\mathfrak{g}\mathfrak{l}(m)$.
We use the (modified) Sweedler notation
$\Delta(X)=\sum X_{1}\otimes X_{2}$ ($X\in \ug$) for
the comultiplication.


\subsection{The quantum analogues of the
symmetric pair $(\mathfrak{g}\mathfrak{l}(2n),
\mathfrak{g}\mathfrak{l}(n)\times \mathfrak{g}\mathfrak{l}(n))$.}

The realization of a two parameter subfamily of the
Macdonald-Koornwinder polynomials as zonal spherical functions
on quantizations of the complex Grassmannian by Noumi, Dijkhuizen
and Sugitani \cite{NDS}, see also \cite{DS},
play a crucial role in the present paper.
A one-parameter family of quantum Grassmannians is defined in
\cite{NDS} in terms of invariance properties with respect to an
one-parameter family of two-sided coideals, which are analogues of the
Lie subalgebra $\mathfrak{k}=\mathfrak{g}\mathfrak{l}(n)\times
\mathfrak{g}\mathfrak{l}(n)\subset\mathfrak{g}=\mathfrak{g}\mathfrak{l}(2n)$.
Although there is a natural way to
define the analogue of the trivial representation of
$\mathfrak{k}$ to the coideal setup, this is
no longer the case for other representations. To get a grip on the
non-trivial representations, we replace the two-sided coideals by
coideal subalgebras of $\ug$, an idea which was suggested
by Noumi \cite{N} and developed in full generality in a series of
papers \cite{L}, \cite{L1}, \cite{L2}, \cite{L3}, \cite{L4},
\cite{L5} by Letzter. For the quantum symmetric pair
$(\gt,\mathfrak{k})=(\mathfrak{g}\mathfrak{l}(2n),
\mathfrak{g}\mathfrak{l}(n)\times \mathfrak{g}\mathfrak{l}(n))$
the coideal algebras have a particularly nice
presentation in terms of generators and relations, closely
resembling the relations of the Drinfeld-Jimbo quantized universal
enveloping algebras (see \cite{L1} and \cite{L4}).

One of the more technical parts of this paper is the translation
of the results of Noumi, Dijkhuizen and Sugitani \cite{NDS} in
the language of right coideal algebras.
Unfortunately, in doing so it turns out
to be more convenient to work
with a slightly modified version of Letzter's right coideal
algebra. In this
subsection we define the modified right coideal algebra for the
symmetric pair $(\gt,\mathfrak{k})$, and we present some of its important
properties.

Let $\{\epsilon_i\}_{i=1}^{2n}$ be the standard orthonormal basis
of the Euclidean space $\bigl(\mathbb{R}^{2n},\bigl(\cdot,\cdot\bigr)\bigr)$
and denote
\[\alpha_j=\epsilon_j-\epsilon_{j+1},\qquad j=1,\ldots,2n-1.
\]
\begin{defi}
Let $\AAA$ be the unital, associative algebra over $\CC$ with generators
$\gamma_i^{\pm 1}$ \textup{(}$i=1,\dots,2n$\textup{)} and
$\beta_j$ \textup{(}$j=1,\dots,2n-1$\textup{)} satisfying the relations
\begin{equation}
\gamma_i\gamma_j=\gamma_j\gamma_i,\qquad \gamma_i\gamma_i^{-1}=1=
\gamma_i^{-1}\gamma_i,\qquad \gamma_i=\gamma_{2n+1-i}
\end{equation}
for $i,j=1,\ldots,2n$;
\begin{equation}
\gamma_i\beta_j=q^{-(\epsilon_i+\epsilon_{2n+1-i},\alpha_j)}\beta_j\gamma_i
\end{equation}
for $i=1,\dots,2n$ and $j=1,\dots,2n-1$;
\begin{equation}
\beta_i\beta_j-\beta_j\beta_i=
q\left(\frac{\gamma_i^{-2}-\gamma_{i+1}^{-2}}{q-q^{-1}}\right)
\delta_{i,2n-j}
\end{equation}
for $i,j=1,\dots, 2n-1$ with $|i-j|\ge 2$;
\begin{equation}\label{relAl0}
q^{-1}\beta_i^2\beta_{i+1}-(q+q^{-1})\beta_i\beta_{i+1}\beta_i+
q\beta_{i+1}\beta_i^2=
q\beta_{n+1}\gamma_n^{-2}\delta_{i,n}
\end{equation}
for $i=1,\ldots,2n-2$ and
\begin{equation}\label{relAl}
q\beta_i^2\beta_{i-1}-(q+q^{-1})\beta_i\beta_{i-1}\beta_i+
q^{-1}\beta_{i-1}\beta_i^2=
q^{-1}\beta_{n-1}\gamma_n^{-2}\delta_{i,n}
\end{equation}
for $i=2,\ldots,2n-1$.
\end{defi}
Observe that $\AAA$ becomes a $*$-algebra with
$*$-structure defined by
\begin{equation*}
\gamma_i^*=\gamma_i,\qquad\qquad
\beta_j^*=\beta_{2n-j}
\end{equation*}
for $i=1,\ldots,2n$ and $j=1,\ldots,2n-1$.
\begin{prop}\label{embeddingsigma}
Let $\sigma\in \mathbb{R}$. The assignment
\begin{equation}\label{assignment}
\begin{split}
\pi_\sigma(\gamma_i^{\pm 1})&=K_i^{\pm 1}K_{2n+1-i}^{\pm 1},\\
\pi_\sigma(\beta_j)&=y_jK_{j+1}^{-1}K_{2n-j}^{-1}+
K_j^{-1}x_{2n-j}K_{2n-j}^{-1},\\
\pi_\sigma(\beta_n)&=y_nK_{n+1}^{-1}K_n^{-1}+K_n^{-1}x_nK_n^{-1}+
\left(\frac{q^{-\sigma}-q^{\sigma}}{q-q^{-1}}\right)K_n^{-2}
\end{split}
\end{equation}
for $i\in\{1,\ldots,2n\}$ and $j\in\{1,\ldots,2n-1\}\setminus \{n\}$,
uniquely extends to an injective $*$-algebra homomorphism
$\pi_\sigma:\AAA\rightarrow \ug$.
\end{prop}
\begin{proof}
This is essentially \cite[Thm. 7.1]{L4} for the
symmetric pair of type AIII (Case 2), see \cite[\S 7]{L4}.
More specifically, following the proof of \cite[Lem. 2.2]{L} one
shows that the assignment
\eqref{assignment} uniquely extends to a unital $*$-algebra homomorphism
$\pi_\sigma: \AAA\rightarrow \ug$ for $\sigma\in\mathbb{R}$.
The injectivity of $\pi_\sigma$ follows by a straightforward modification
of the proof of \cite[Prop. 2.3]{L}.
\end{proof}
\begin{defi}\label{defiLetzter}
Let $\sigma\in\mathbb{R}$.
We call $(\ug,\AAA_\sigma)$ with
$\AAA_\sigma=\pi_\sigma(\AAA)$,
a quantum analogue of the symmetric
pair $(\gt,\mathfrak{k})=(\mathfrak{g}\mathfrak{l}(2n),
\mathfrak{g}\mathfrak{l}(n)\times
\mathfrak{g}\mathfrak{l}(n))$. We denote
\begin{equation}\label{embed}
C_i=\pi_\sigma(\gamma_i),\qquad B_j=\pi_\sigma(\beta_j),
\qquad B_n^\sigma=\pi_\sigma(\beta_n)
\end{equation}
\textup{(}$i\in\{1,\ldots,2n\}$ and $j\in\{1,\ldots,2n-1\}\setminus
\{n\}$\textup{)} for the image of the generators of $\AAA$ under the embedding
$\pi_\sigma$.
\end{defi}
We fix $\sigma\in\mathbb{R}$ once and for all, unless specified differently.

The terminology in Definition \ref{defiLetzter} can formally be
justified as follows, cf., e.g., \cite[\S 3]{L}.
Taking the classical limit $q\rightarrow 1$ of the (modified) generators
of $\AAA_\sigma$ gives
\begin{equation}\label{classgenerator}
\begin{split}
&\frac{C_i-C_i^{-1}}{q-q^{-1}}\rightarrow E_{i,i}+E_{2n+1-i,2n+1-i},\\
&B_j\rightarrow E_{j+1,j}+E_{2n-j,2n+1-j},\\
&B_n^\sigma\rightarrow E_{n+1,n}+E_{n,n+1}-\sigma 1
\end{split}
\end{equation}
for $i\in\{1,\ldots,n\}$ and $j\in\{1,\ldots,2n-1\}\setminus
\{n\}$, where $E_{i,j}\in\mathfrak{g}\subset U(\mathfrak{g})$
is the matrix unit with
a one in row $i$ and column $j$, and zeroes elsewhere.
The unital subalgebra $\widetilde{\AAA}$ of $U(\mathfrak{g})$ generated by the
classical limits \eqref{classgenerator} of the generators of $\AAA_\sigma$
is independent of $\sigma$ and isomorphic to the subalgebra
$U(\mathfrak{k})\subset U(\mathfrak{g})$.
The isomorphism is induced by the inner
automorphism $\hbox{Ad}(g)$ of $\mathfrak{g}$,
with $g\in \hbox{GL}(2n;\mathbb{C})$
given explicitly by
\begin{equation}\label{forg}
g=\frac{1}{\sqrt{2}}\sum_{i=1}^n\bigl(E_{i,i}+
E_{i,2n+1-i}+E_{2n+1-i,i}-E_{2n+1-i,2n+1-i}\bigr).
\end{equation}
In fact, the associated inner automorphism of $U(\mathfrak{g})$
preserves all classical generators
\eqref{classgenerator} of $\widetilde{\AAA}$
besides $E_{n+1,n}+E_{n,n+1}-\sigma 1$,
which is mapped to $E_{n,n}-E_{n+1,n+1}-\sigma 1$.

The classical isomorphism ${\widetilde{\AAA}}\simeq
U(\mathfrak{k})$ via the inner automorphism $\hbox{Ad}(g)$ has the
following weak analogue on the quantum group level.
Let $U_q(\mathfrak{k})$ be the Hopf $*$-subalgebra of $\ug$
generated by $K_i^{\pm 1}$ ($i\in\{1,\ldots,2n\}$) and $x_j,y_j$
($j\in\{1,\ldots,2n-1\}\setminus \{n\}$). Denote
\begin{equation}\label{Bnalt}
\widehat{B}_n^\sigma=\frac{q^{-\sigma}K_{n+1}^{-2}-q^{\sigma}K_{n}^{-2}}
{q-q^{-1}}\in U_q(\mathfrak{k}),
\end{equation}
which has $E_{n,n}-E_{n+1,n+1}-\sigma 1$ as classical limit.
\begin{prop}\label{repA}
Let $I\subseteq U_q(\mathfrak{k})$ be the two-sided $*$-ideal
generated by $x_{n-1}y_{n+1}$. The assignment
\[\phi_\sigma(\gamma_i)=C_i+I,\qquad
\phi_\sigma(\beta_j)=B_j+I,\qquad
\phi_\sigma(\beta_n)=\widehat{B}_n^\sigma+I
\]
for $i\in\{1,\ldots,n\}$ and $j\in\{1,\ldots,2n-1\}\setminus\{n\}$
uniquely extends to a $*$-algebra homomorphism $\phi_\sigma: \AAA\rightarrow
U_q(\mathfrak{k})/I$.
\end{prop}
\begin{proof}
In view of Proposition \ref{embeddingsigma}, $\phi_\sigma$ extends
to an algebra homomorphism $\phi_\sigma:\AAA\rightarrow
U_q(\mathfrak{k})/I$ if $\phi_\sigma$ respects the defining
relations of $\AAA$ involving $\beta_n$. This is a straighforward,
but tedious computation. The $*$-ideal $I$ plays only a role for
the relations \eqref{relAl0} for $i=n-1$ and relation
\eqref{relAl} for $i=n+1$ (which in turn are each others
$*$-images).
In $U_q(\mathfrak{k})$ we have the relation
\begin{equation*}\label{repl2}
q^{-1}B_{n-1}^2\widehat{B}_n^\sigma
-(q+q^{-1})B_{n-1}\widehat{B}_n^\sigma B_{n-1}
+q\widehat{B}_n^\sigma B_{n-1}^2
=q^{-\sigma-1}(q^{-2}-q^{2})x_{n+1}y_{n-1}K_{n-1}^{-1}K_n^{-1}K_{n+1}^{-4},
\end{equation*}
which, modulo $I$, reduces to the $\phi_\sigma$-image of relation
\eqref{relAl0} for $i=n-1$. The case $i=n+1$ is checked similarly.
It is clear that the algebra homomorphism
$\phi_\sigma$ respects the $*$-structure.
\end{proof}

\begin{rema}\label{Rosremark}
Denote $U_q(\mathfrak{g}\mathfrak{l}(2))^{(n)}\subset \ug$
for the copy of the unital Hopf-$*$-algebra $U_q(\mathfrak{g}\mathfrak{l}(2))$
in $\ug$, generated by $x_n,y_n,K_n^{\pm 1},
K_{n+1}^{\pm 1}$. Note that both
$B_n^\sigma\in U_q(\mathfrak{g}\mathfrak{l}(2))^{(n)}$ and
$\widehat{B}_n^\sigma \in U_q(\mathfrak{g}\mathfrak{l}(2))^{(n)}$,
with $\widehat{B}_n^\sigma$ given by
\eqref{Bnalt}. By \cite{Ro}, there exists an invertible
element $x_\sigma$ in a suitable completion of
$U_q(\mathfrak{g}\mathfrak{l}(2))^{(n)}$
such that $x_{\sigma} B_n^\sigma
x_{\sigma}^{-1}=\widehat{B}_n^\sigma$. Explicitly, $x_\sigma$ can be given as
formal infinite series by
\begin{equation}\label{Rosengrenelement}
x_\sigma=\sum_{l,m=0}^{\infty}\frac{q^{-(l+m)\sigma}(-1)^m}
{\bigl(-q^{2-2\sigma};q^2\bigr)_l}
\frac{q^{l^2+2lm-l-m}}{\bigl(q^2;q^2\bigr)_l\bigl(q^2;q^2\bigr)_m}
(1-q^2)^{l+m}x_n^l\bigl(y_nK_nK_{n+1}^{-1}\bigr)^m.
\end{equation}
The element $x_\sigma$ is essentially
Babelon's \cite{B} vertex-IRF transformation,
which leads to an
interpretation of the parameter $\sigma$ as a dynamical parameter in
the sense of dynamical quantum groups, see \cite{St2} for a detailed discussion
and further references. Note
furthermore that conjugating by $x_\sigma$ fixes all other generators
$B_j$ and $C_i^{\pm 1}$ of $\AAA_\sigma$ besides
$B_{n-1}$ and $B_{n+1}$.
\end{rema}
Let $\mathcal{E}\subset \mathcal{A}_\sigma\subset \ug$ be the
unital subalgebra generated by the
$\sigma$-independent generators $C_i^{\pm 1}$ ($i\in\{1,\ldots,2n\}$)
and $B_j$ ($j\in\{1,\ldots,2n-1\}\setminus \{n\}$) of
$\mathcal{A}_\sigma$. We end this subsection by proving that
$\mathcal {E}$ is essentially the subalgebra
$U_q(\mathfrak{g}\mathfrak{l}(n))$, diagonally embedded in $\ug$.

Throughout this paper we identify
\[U_q(\mathfrak{g}\mathfrak{l}(n))\otimes
U_q(\mathfrak{g}\mathfrak{l}(n))\simeq
U_q(\mathfrak{k})\subset \ug
\]
as $*$-Hopf algebras, by identifying $x_j\otimes 1, y_j\otimes 1$
and $K_i^{\pm 1}\otimes 1$ (respectively $1\otimes x_j, 1\otimes
y_j$ and $1\otimes K_i^{\pm 1}$) with $x_j, y_j$ and $K_i^{\pm 1}$
(respectively $x_{n+j},y_{n+j}$ and $K_{n+i}^{\pm 1}$) for
$j=1,\ldots,n-1$ and $i=1,\ldots,n$. Let $\psi:
U_q(\mathfrak{g}\mathfrak{l}(n))\rightarrow
U_q(\mathfrak{g}\mathfrak{l}(n))$ be the $*$-algebra isomorphism
defined by
\begin{equation}\label{psifunction}
 \psi(x_j)=y_{n-j},\qquad \psi(y_j)=x_{n-j},\qquad
\psi(K_i^{\pm 1})=K_{n+1-i}^{\pm 1}
\end{equation}
for $j=1,\ldots,n-1$ and $i=1,\ldots,n$. Let $\Delta^{op}$ be the
opposite comultiplication of $U_q(\mathfrak{g}\mathfrak{l}(n))$.
\begin{lem}\label{diagonal}
For $j=1,\ldots,n-1$ and $i=1,\ldots,n$ we have
\begin{equation*}
\begin{split}
(\hbox{id}\otimes \psi)&\Delta^{op}(y_jK_{j+1}^{-1})=B_j,\\
(\hbox{id}\otimes \psi)&\Delta^{op}(x_jK_j^{-1})=B_{2n-j},\\
(\hbox{id}\otimes \psi)&\Delta^{op}(K_i^{\pm 1})=C_i^{\pm 1}.
\end{split}
\end{equation*}
Furthermore, $(\hbox{id}\otimes\psi)\Delta^{op}$ defines a $*$-algebra
isomorphism
\[(\hbox{id}\otimes\psi)\Delta^{op}:
U_q(\mathfrak{g}\mathfrak{l}(n))\overset{\sim}{\longrightarrow}
\mathcal{E}\subseteq \mathcal{A}_\sigma\cap U_q(\mathfrak{k})\subset \ug.
\]
\end{lem}
\begin{proof}
The proof is straightforward.
\end{proof}
Other crucial properties of the quantum symmetric pairs $(\ug,
\AAA_\sigma)$, such as the coideal property of $\AAA_\sigma\subset
\ug$, are discussed in later sections.


\subsection{Representations.}

We discuss those aspects of the
representation theory associated to the quantum symmetric pairs $(\ug,
\AAA_\sigma)$ which are relevant for the purposes of this paper.
We start by recalling some standard facts on the finite dimensional
representation theory of the quantized universal enveloping algebra
$U_q(\mathfrak{g}\mathfrak{l}(m))$. If no confusion can arise,
we surpress in the following definitions the dependence on $m$
as much as possible.

Let $P_m=\bigoplus_{i=1}^m\mathbb{Z}\varepsilon_i\simeq
\mathbb{Z}^{\times m}$ be the rational character lattice of
$\hbox{GL}(m;\mathbb{C})$ and
\[P_{m}^+=\{(\lambda_1,\dots,\lambda_{m})\in P_{m}\,|\,\lambda_1\ge\dots
\ge\lambda_{m}\}\]
the associated cone of dominant weights.

We only consider finite dimensional left
$U_q(\mathfrak{g}\mathfrak{l}(m))$-modules $M$ with weights in $P_m$,
i.e. we assume that $M$ has weight decomposition
\begin{equation*}
\begin{split}
M&=\bigoplus_{\mu\in P_m}M[\mu],\\
M[\mu]&=\{ v\in M \, | \, K_iv=q^{\mu_i}v\,\,\,i=1,\ldots,m \}.
\end{split}
\end{equation*}
Any such finite dimensional $\ug$-module $M$
is $\ug$-semisimple. The irreducible
ones are the irreducible finite dimensional highest
weight representations $L_\lambda$ of $U_q(\mathfrak{g}\mathfrak{l}(m))$ with
highest weight $\lambda\in P_m^+$.

The irreducible module $L_\lambda$ ($\lambda\in P_m^+$)
has a one-dimensional weight space $L_\lambda[\lambda]$.
A vector $v_\lambda$ spanning $L_\lambda[\lambda]$ is called a
highest weight vector of $L_\lambda$. A highest weight vector $v_\lambda$
is a cyclic vector of $L_\lambda$ satisfying $x_jv_\lambda=0$ for
$j=1,\ldots,m-1$. In particular,
\[ L_\lambda=\bigoplus_{\mu\in P_{m}: \mu\preceq \lambda} L_\lambda[\mu],
\]
with $\preceq$ the
$A$-type dominance order on $P_{m}$:
$\lambda\preceq \mu$ for $\lambda,\mu\in P_{m}$ if
\[\sum_{i=1}^j\lambda_i\leq \sum_{i=1}^j\mu_i\quad
(j\in\{1,\ldots,m-1\}),\quad
\quad \sum_{i=1}^{m}\lambda_i= \sum_{i=1}^{m}\mu_i.
\]

There exists a scalar product $\langle
\cdot,\cdot\rangle_\lambda$ on $L_\lambda$,
unique up to constant multiples, such that
\[
\langle Xv,w\rangle_\lambda=\langle v,X^*w\rangle_\lambda,\qquad
\forall\,v,w\in L_\lambda,\,\,\,\forall\,
X\in U_q(\mathfrak{g}\mathfrak{l}(m)).
\]
Consequently, any finite dimensional
$U_q(\mathfrak{g}\mathfrak{l}(m))$-module $M$ is $*$-unitarizable.

\begin{defi}
We write $M^\sigma$ for a left $\ug$-module $M$, viewed as
$\AAA_\sigma$-module by restriction of the action to the
subalgebra $\AAA_\sigma\subset \ug$.
\end{defi}
Since $\AAA_\sigma\subset \ug$ is $*$-invariant,
we have the following result, cf., e.g., \cite[Thm. 3.3]{L2}.
\begin{lem}\label{smodule}
If $M$ is a finite dimensional $\ug$-module, then $M^\sigma$ is
$\AAA_\sigma$-semisimple.
\end{lem}
An abstract representation theory for $\AAA_\sigma$ was
developed by Letzter \cite[\S 5-7]{L2} in the general context of
quantum symmetric pairs.
For our purposes it is more convenient to have concrete
realizations of certain special $\AAA_\sigma$-representations. These special
representations are constructed using the following corollary
of Proposition \ref{repA}.
\begin{cor}\label{twist}
Let $V$ be a finite dimensional Hilbert space and let
$\rho: U_q(\mathfrak{k})\rightarrow \hbox{End}_{\mathbb{C}}(V)$
be a $*$-representation such that $I\subseteq \hbox{Ker}(\rho)$.
Denote by $\overline{\rho}: U_q(\mathfrak{k})/I\rightarrow
\hbox{End}_{\mathbb{C}}(V)$ the associated $*$-representation
of $U_q(\mathfrak{k})/I$. Then
\[\rho_\sigma=\overline{\rho}\circ\phi_\sigma\circ\pi_\sigma^{-1}:
\AAA_\sigma\rightarrow \hbox{End}(V_\sigma),
\]
with $V_\sigma=V$ as finite dimensional Hilbert space,
defines a $*$-representation of $\AAA_\sigma$ such that
$\rho_\sigma|_{\mathcal{E}}=\rho|_{\mathcal{E}}$.
\end{cor}

\begin{rema}
With the conventions of Corollary \ref{twist} we
can define a $*$-representation
$\rho_{\sigma,\tau}=
\overline{\rho}\circ\phi_\tau\circ\pi_\sigma^{-1}$ of $\AAA_\sigma$
for arbitrary  $\sigma,\tau\in\mathbb{R}$. The $*$-representation
$\rho_{\sigma,\tau}$ thus satisfies
\[\rho_{\sigma,\tau}(B_n^\sigma)=\rho(\widehat{B}_n^\tau).
\]
The special role of $\rho_{\sigma}=\rho_{\sigma,\sigma}$ follows from
the fact that $B_n^\sigma$ and $\widehat{B}_n^\sigma$ are conjugate
in $\ug$ (cf. Remark \ref{Rosremark}) and from the fact that the formal
classical limit of $\rho_\sigma=\rho_{\sigma,\sigma}$
is equivalent to the formal classical limit of $\rho$ via the isomorphism
$\widetilde{\AAA}\simeq U_q(\mathfrak{k})$ constructed in \S 1.3.
\end{rema}

Let $\kappa_1\in \mathbb{Z}$ and $\kappa\in\mathbb{Z}_{\geq 0}$. Denote
$V(\kappa,\kappa_1)$ for the irreducible, finite dimensional
$U_q(\mathfrak{k})$-module
\[ V(\kappa,\kappa_1)=L_{(-\kappa_1)^n}\otimes
L_{(\kappa_1+(n-1)\kappa, (\kappa_1-\kappa)^{n-1})}.
\]
Here we have used the shorthand notation $(\kappa^m)\in P_m^+$
for the $m$-tuple with all entries equal to $\kappa$.
The representation map of $V(\kappa,\kappa_1)$
will be denoted by $\rho(\kappa,\kappa_1)$.
Clearly there exists a scalar product on
$V(\kappa,\kappa_1)$, unique up to constant multiples, such that
$\rho(\kappa,\kappa_1): U_q(\mathfrak{k})\rightarrow
\hbox{End}_{\mathbb{C}}(V(\kappa,\kappa_1))$ is a $*$-representation. Since
$L_{(-\kappa_1)^n}$ is one-dimensional,
$x_{n-1}\in U_q(\mathfrak{g}\mathfrak{l}(n))$ acts as zero,
hence $I\subseteq \hbox{ker}(\rho(\kappa,\kappa_1))$. Thus we obtain
a $*$-representation
\[\rho(\kappa,\kappa_1)_\sigma: \AAA_\sigma\rightarrow
\hbox{End}_{\mathbb{C}}(V(\kappa,\kappa_1)_\sigma)
\]
by the previous corollary. To simplify notations, we denote
\[(\rho(\kappa_1)_\sigma,V(\kappa_1)_\sigma)=
(\rho(0,\kappa_1)_\sigma,V(0,\kappa_1)_\sigma).
\]

The $U_q(\mathfrak{g}\mathfrak{l}(n))$-module
$L_{((n-1)\kappa,(-\kappa)^{n-1})}$ is the quantum analogue of the
$\hbox{GL}(n;\mathbb{C})$-module
$S^{n\kappa}(\mathbb{C}^n)\det^{-\kappa}$, where
$S^m(\mathbb{C}^n)$ is the space of homogeneous polynomials on
$\mathbb{C}^n$ of degree $m$. This is precisely the
$U_q(\mathfrak{g}\mathfrak{l}(n))$-module which plays the crucial
role in Etingof's and Kirillov's \cite{EK} theory on A-type
Macdonald polynomials and generalized quantum group characters.
The module $L_{((n-1)\kappa,(-\kappa)^{n-1})}$ admits a concrete
description (see, e.g., \cite[\S 5]{EK}) which
can be directly lifted to $V(\kappa,\kappa_1)_\sigma$,
leading to the following result.
\begin{lem}\label{explicit}
Let $\mathcal{J}_\kappa$ be the set of $n$-tuples $\underline{m}=
(m_1,\ldots,m_n)\in\mathbb{Z}_{\geq 0}^n$ which sum up to $n\kappa$.
There exists a $\sigma$-independent basis $\{ r_{\underline{m}}\, | \,
\underline{m}\in \mathcal{J}_\kappa\}$
of $V(\kappa,\kappa_1)$ such that
\begin{equation*}
\begin{split}
\rho(\kappa,\kappa_1)_\sigma(C_i)r_{\underline{m}}&=
q^{m_i-\kappa}r_{\underline{m}},\\
\rho(\kappa,\kappa_1)_\sigma(B_j)r_{\underline{m}}&=
q^{-m_{j+1}+\kappa}\left(\frac{q^{m_{j}}-q^{-m_{j}}}{q-q^{-1}}\right)
r_{\underline{m}-\varepsilon_{j}+\varepsilon_{j+1}},\\
\rho(\kappa,\kappa_1)_\sigma(B_{2n-j})r_{\underline{m}}&=
q^{-m_{j}+\kappa}\left(\frac{q^{m_{j+1}}-q^{-m_{j+1}}}{q-q^{-1}}\right)
r_{\underline{m}+\varepsilon_{j}-\varepsilon_{j+1}},\\
\rho(\kappa,\kappa_1)_\sigma(B_n^\sigma)r_{\underline{m}}&=
\left(\frac{q^{-\sigma+2(\kappa-\kappa_1-m_n)}-q^{\sigma+2\kappa_1}}
{q-q^{-1}}\right)r_{\underline{m}}
\end{split}
\end{equation*}
for $i=1,\ldots,n$ and $j=1,\ldots,n-1$, where we have used the
convention that $r_{\underline{m}}=0$ if
$\underline{m}\not\in \mathcal{J}_{\kappa}$.
\end{lem}
Lemma \ref{explicit} implies that the common eigenspace
\begin{equation}\label{zeroweight}
\widetilde{V}(\kappa,\kappa_1)_\sigma=\{v\in
V(\kappa,\kappa_1)_\sigma \,\,\, | \,\,\,
\rho(\kappa,\kappa_1)_\sigma(C_i)v=v\,\,\,\,\forall i\}
\end{equation}
is one-dimensional and spanned by $r_{(\kappa^n)}$.

The following lemma follows directly from Lemma \ref{explicit}.
\begin{lem}\label{restrictionmodule}
The $\AAA_\sigma$-module
$V(\kappa,\kappa_1)_\sigma$, viewed as
$U_q(\mathfrak{g}\mathfrak{l}(n))$-module by restricting the module
structure to $\mathcal{E}\subset \AAA_\sigma$ and using the isomorphism
$\mathcal{E}\simeq U_q(\mathfrak{g}\mathfrak{l}(n))$ of
Lemma \ref{diagonal}, is isomorphic to the simple
$U_q(\mathfrak{g}\mathfrak{l}(n))$-module
$L_{((n-1)\kappa,(-\kappa)^{n-1})}$.
In particular, $V(\kappa,\kappa_1)_\sigma$ is a simple $\AAA_\sigma$-module.
\end{lem}

We end this subsection by formulating
branching rules for the quantum symmetric pair
$(\ug,\AAA_\sigma)$. For the formulation it is convenient to
define an injective map ${}^\natural: \Lambda_n\rightarrow P_{2n}$
by
\begin{equation}
\mu^\natural=(\mu_1,\mu_2,\ldots,\mu_n,-\mu_n,\ldots,
-\mu_2,-\mu_1).
\end{equation}
Observe that ${}^\natural$ restricts to a map
${}^\natural: \Lambda_n^+\rightarrow P_{2n}^+$, and that ${}^\natural$
respects the dominance order, so
$\mu^\natural\preceq\nu^\natural$ for $\mu,\nu\in \Lambda_n$
if $\mu\leq\nu$.

\begin{prop}\label{branching}
Let $\kappa_2\in\mathbb{Z}$ and $\kappa,\kappa_1\in\mathbb{Z}_{\geq 0}$
with $-\kappa_1\leq\kappa_2\leq\kappa_1$, and set
\begin{equation}\label{rho}
\delta(\kappa,\kappa_1)=(\kappa_1+(n-1)\kappa,\kappa_1+(n-2)\kappa,
\ldots,\kappa_1)\in\Lambda_n^+.
\end{equation}
Let $\lambda\in P_{2n}^+$. The simple $\AAA_\sigma$-modules
$V(\kappa,\kappa_1)_\sigma$ and $V(\kappa_2)_\sigma$ are both
constituents of the $\AAA_\sigma$-module $L_\lambda^\sigma$ if and only if
$\lambda\in (\Lambda_n^++\delta(\kappa,\kappa_1))^\natural$.
For such $\lambda$, both
$V(\kappa,\kappa_1)_\sigma$ and $V(\kappa_2)_\sigma$ occur with
multiplicity one in $L_\lambda^\sigma$.
\end{prop}
For $\kappa=\kappa_1=\kappa_2=0$, Proposition \ref{branching}
states which of the simple, finite dimensional
$\ug$-modules $L_\lambda$ are spherical with respect to $\AAA_\sigma$.
This was proven for general quantum symmetric pairs in \cite[Thm. 4.3]{L2}.
This special case can also be derived from the results in \cite{NDS},
see Lemma \ref{zonal}.
The classical case ($q=1$) of Proposition \ref{branching} was
proven in \cite{O} using the
Littlewood-Richardson rule. Finally, for $n=1$ (in which case the
parameter $\kappa$ is reduntant), Proposition \ref{branching} is
due to Koornwinder \cite{Koo0}, who used $q$-special functions for
the proof. The general proof of Proposition \ref{branching} is
discussed in \S 2.
\begin{rema}\label{notationparameters}
In the remainder of the paper we always assume the conditions
$\kappa,\kappa_1\in\mathbb{Z}_{\geq 0}$ and
$-\kappa_1\leq\kappa_2\leq\kappa_1$ on
the representation labels $\kappa,\kappa_1,\kappa_2\in\mathbb{Z}$,
and we use the short hand notation
$\vec{\kappa}=(\kappa_1,\kappa_2,\kappa)$. There are three other
parameter domains for $\vec{\kappa}$ for which similar results can
be derived with only minor alterations, namely
$\kappa_1,\kappa_2\in \mathbb{Z}$, $\kappa\in\mathbb{Z}_{\geq 0}$
satisfying $\kappa_1\in\mathbb{Z}_{\leq 0}$ and $\kappa_1\leq \kappa_2\leq
-\kappa_1$, or $\kappa_2\in\mathbb{Z}_{\geq 0}$ and
$-\kappa_2\leq \kappa_1\leq
\kappa_2$, or $\kappa_2\in\mathbb{Z}_{\leq 0}$ and $\kappa_2\leq \kappa_1\leq
-\kappa_2$, compare with \cite{Koe} for the special case $n=1$ and with
\cite{O} for the classical case ($q=1$).
\end{rema}


\subsection{Vector valued spherical functions.}
We write $G=\hbox{GL}(2n;\mathbb{C})$ for the general linear group.
The quantized algebra of regular functions
$\mathbb{C}_q[G]\subseteq \ug^*$ is
the span of the matrix coefficients of the irreducible
$\ug$-representations
$L_\lambda$ ($\lambda\in P_{2n}^+$). The quantized
function algebra $\mathbb{C}_q[G]$ inherets from $\ug$
the structure of a Hopf $*$-algebra. In particular, the
$*$-structure on $\mathbb{C}_q[G]$ is defined by
\[f^*(X)=\overline{f(S(X)^*)},\qquad \forall\, X\in \ug.
\]
The left and right regular $\ug$-action on $\mathbb{C}_q[G]$ is
defined by
\[ (Y\cdot f\cdot Z)(X)=f(ZXY),\qquad X,Y,Z\in \ug.
\]
The Peter-Weyl decomposition
\begin{equation}\label{PW}
\mathbb{C}_q[G]=\bigoplus_{\lambda\in P_{2n}^+}W(\lambda)
\end{equation}
with $W(\lambda)$ the span of the matrix coefficients of $L_\lambda$,
is the irreducible decomposition of $\mathbb{C}_q[G]$ as $\ug$-bimodule.
For any vector space $V$, we may and will view elements
$f\in \mathbb{C}_q[G]\otimes V$ as linear maps
$f: \ug\rightarrow V$.

Recall the conventions and notations for the representation labels
$\vec{\kappa}$ from Remark \ref{notationparameters}. We fix
$\sigma,\tau\in\mathbb{R}$ once and for all, unless specified differently.
In the following definition we identify
$V(\kappa_2)_\sigma\simeq\mathbb{C}$ as vector spaces and view
$\rho(\kappa_2)_\sigma$ as character of $\AAA_\sigma$.
\begin{defi}
We call $f\in \mathbb{C}_q[G]\otimes
V(\kappa,\kappa_1)_\tau$ a vector
valued spherical function if for all $X\in \ug$,
\begin{equation*}
\begin{split}
f(Xa)&=\rho(\kappa_2)_\sigma(a)\,f(X),\qquad \forall\, a\in\AAA_\sigma,\\
f(bX)&=\rho(\kappa,\kappa_1)_\tau(b)\,f(X),\qquad \forall\,
b\in\AAA_\tau.
\end{split}
\end{equation*}
We denote by $F_{\vec{\kappa}}^{\sigma,\tau}$  the space of vector
valued spherical functions.
\end{defi}
By standard arguments (compare, e.g., with \cite{O} for the classical setup,
and with Lemma \ref{groundstateform}),
we obtain the following corollary of Lemma \ref{smodule} and
Proposition \ref{branching}.
\begin{cor}\label{Cormis}
The space of vector valued spherical
functions $F_{\vec{\kappa}}^{\sigma,\tau}$
decomposes as
\[ F_{\vec{\kappa}}^{\sigma,\tau}=\bigoplus_{\mu\in \Lambda_n^+}
F_{\vec{\kappa}}^{\sigma,\tau}(\mu),
\]
with $F_{\vec{\kappa}}^{\sigma,\tau}(\mu)$ the one-dimensional vector space
\[F_{\vec{\kappa}}^{\sigma,\tau}(\mu)=F_{\vec{\kappa}}^{\sigma,\tau}\cap
\Bigl(W\bigl((\mu+\delta(\kappa,\kappa_1))^\natural\bigr)\otimes
V(\kappa,\kappa_1)_\tau\Bigr),\qquad
\mu\in\Lambda_n^+.
\]
\end{cor}
We call a function $0\not=f\in
F_{\vec{\kappa}}^{\sigma,\tau}(\mu)$ an {\it elementary}
vector valued spherical function of degree $\mu\in\Lambda_n^+$.

We consider now the restriction of vector valued spherical functions
to the commutative
subalgebra $U^0\subset \ug$ generated by the
group-like elements $K_j^{\pm 1}$ ($j=1,\ldots,2n$). Denote
\[K^{\lambda}=K_1^{\lambda_1}K_2^{\lambda_2}\cdots
K_{2n}^{\lambda_{2n}},\qquad  \lambda\in P_{2n},
\]
which form a linear basis of $U^0$,
and define $\delta=\delta_{2n}\in P_{2n}^+$ by
\[\delta=(2n-1,2n-2,\ldots,1,0).
\]
Let $\mathbb{C}[z^{\pm 1}]$ be the algebra of Laurent polynomials
in $2n$ variables $z_1,\ldots,z_{2n}$.
\begin{defi}
Let $V$ be a complex vector space and $f\in \mathbb{C}_q[G]\otimes V$,
viewed as linear map $f: \ug\rightarrow V$. We define a regular function
$f|_T: (\mathbb{C}^*)^{2n}\rightarrow V$ \textup{(}or, equivalently,
$f|_T\in\mathbb{C}[z^{\pm 1}]\otimes V$\textup{)}
by the requirement that
\[f|_T(q^{\lambda})=f(K^{\lambda-\delta})\qquad \forall\,\lambda\in
P_{2n},
\]
where $q^{\lambda}=(q^{\lambda_1},\ldots,q^{\lambda_{2n}})$.
The resulting linear map
\[|_T: \mathbb{C}_q[G]\otimes V\rightarrow \mathbb{C}[z^{\pm 1}]\otimes V,
\qquad f\mapsto f|_T
\]
is called the restriction map.
\end{defi}
We define elements
$u_i\in \mathbb{C}[z^{\pm 1}]$ by
\begin{equation}
u_i=z_iz_{2n+1-i}^{-1},\qquad i=1,\ldots,n,
\end{equation}
and we write $\mathbb{C}[u^{\pm 1}]\subset \mathbb{C}[z^{\pm 1}]$
for the subalgebra generated by the $u_i^{\pm 1}$. Recall
the one-dimensional subspace
$\widetilde{V}(\kappa,\kappa_1)_\tau \subseteq
V(\kappa,\kappa_1)_\tau$ defined by \eqref{zeroweight}.
\begin{lem}\label{restriction}
If $f\in F_{\vec{\kappa}}^{\sigma,\tau}\subset
\mathbb{C}_q[G]\otimes V(\kappa,\kappa_1)_\tau$, then $f|_T\in
\mathbb{C}[u^{\pm 1}]\otimes \widetilde{V}(\kappa,\kappa_1)_\tau$.
Identifying the one-dimensional vector space
$\widetilde{V}(\kappa,\kappa_1)_\tau$ with $\mathbb{C}$, the
restriction map thus gives rise to a linear map
\[ |_T: F_{\vec{\kappa}}^{\sigma,\tau}\rightarrow
\mathbb{C}[u^{\pm 1}].
\]
\end{lem}
\begin{proof}
Fix $f\in F_{\vec{\kappa}}^{\sigma,\tau}$. For $i=1,\ldots,n$
and $X\in U^0$ we compute
\begin{equation*}
\begin{split}
\rho(\kappa,\kappa_1)_\tau(C_i)f(X)&=f(C_iX)\\
&=f(XC_i)\\
&=\rho(\kappa_2)_\sigma(C_i)f(X)=f(X),
\end{split}
\end{equation*}
hence $f(X)\in \widetilde{V}(\kappa,\kappa_1)_\tau$.
Identifying
$\widetilde{V}(\kappa,\kappa_1)_\tau\simeq\mathbb{C}$,
we thus conclude that $f|_T\in\mathbb{C}[z^{\pm 1}]$.
The formulas $f(XC_i)=f(X)$ for $X\in U^0$ and
$i=1,\ldots,n$ imply that $f|_T\in\mathbb{C}[u^{\pm 1}]$.
\end{proof}
The following main result of this paper gives a representation
theoretic interpretation of the full five parameter family of
Macdonald-Koornwinder polynomials (with two parameters continuous,
and three parameters discrete).

\begin{thm}\label{maintheorem}
Fix parameters $\sigma,\tau\in \mathbb{R}$ and
$\vec{\kappa}=(\kappa_1,\kappa_2,\kappa)\in \mathbb{Z}^{\times 3}$
with $\kappa,\kappa_1\in\mathbb{Z}_{\geq 0}$ and $-\kappa_1\leq \kappa_2\leq
\kappa_1$. Let $f_0\in F_{\vec{\kappa}}^{\sigma,\tau}(0)$ be an
elementary vector valued spherical function of degree $0\in\Lambda_n^+$.

{\bf (i)}
The restriction map $|_T$ defines a linear bijection
\[ |_T: F_{\vec{\kappa}}^{\sigma,\tau}\rightarrow
f_0|_T\,\mathbb{C}[u^{\pm 1}]^W.
\]

{\bf (ii)} We have
\begin{equation*}
\begin{split}
f_0|_T=C\,u^{\delta(\kappa,\kappa_1)}&\prod_{i=1}^n
\bigl(q^{1-\sigma+\tau}u_i^{-1};q^2\bigr)_{\kappa_1-\kappa_2}
\bigl(-q^{1+\sigma+\tau}u_i^{-1};q^2\bigr)_{\kappa_1+\kappa_2}\\
&\times \prod_{1\leq i<j\leq n}\bigl(q^2u_i^{-1}u_j,
q^2u_i^{-1}u_j^{-1};q^2\bigr)_{\kappa}
\end{split}
\end{equation*}
for some nonzero constant $C$.

{\bf (iii)}
If $f_\mu\in F_{\vec{\kappa}}^{\sigma,\tau}(\mu)$
is an elementary vector valued
spherical function of degree $\mu\in \Lambda_n^+$, then
\[\frac{f_\mu|_T}{f_0|_T}=
D\,P_\mu\bigl(u;-q^{\sigma+\tau+1+\kappa_1+\kappa_2},
-q^{-\sigma-\tau+1}, q^{\sigma-\tau+1}, q^{-\sigma+\tau+1+\kappa_1-\kappa_2};
q^2,q^{2\kappa+2}\bigr)
\]
for some nonzero constant $D$.
\end{thm}
\begin{rema}
In the special case $n=1$ Theorem \ref{maintheorem} gives a
representation theoretic
interpretation of the Askey-Wilson polynomials,
which is in accordance with the results from \cite{Koo0}, \cite{NM}
and \cite{Koe}. The classical case ($q=1$) of Theorem \ref{maintheorem}
was proven in \cite{O}. In \S 3 we show that
Theorem \ref{maintheorem} for $\kappa_1=\kappa_2=\kappa=0$ follows from
the main results in \cite{NDS}.
\end{rema}
The proof of Theorem \ref{maintheorem}{\bf (i)}, {\bf (ii)}
and {\bf (iii)} is given in \S 4, \S 5 and \S 6, respectively.

\section{Branching rules}
The main goal of this section is to establish the specific branching rules
for the quantum symmetric pairs $(\ug,\AAA_\sigma)$ as formulated
in Proposition \ref{branching}.
We start with defining and studying the notion of the expectation value
for the quantum symmetric pairs $(\ug,\AAA_\sigma)$.


\subsection{The expectation value}
The expectation value for intertwiners plays a crucial role in
the representation theoretic approach to quantum field theory, see, e.g.,
\cite{EV}. The notion of expectation value
is naturally associated to symmetric pairs $(H\times H,
\hbox{diag}(H))$ with $H$ a semisimple Lie group and
$\hbox{diag}(H)$ the diagonal embedding of $H$ in $H\times H$, as well
as to their quantum analogues $(U_q(\mathfrak{h})\otimes U_q(\mathfrak{h}),
\Delta(U_q(\mathfrak{h})))$ with $\mathfrak{h}$ a semisimple
Lie algebra.
In this subsection we define the expectation value for
the quantum symmetric pair $(\ug,\AAA_\sigma)$. As we shall see, this
directly leads to information on branching rules for the quantum symmetric
pairs $(\ug,\AAA_\sigma)$.

For $\lambda\in P_{2n}^+$ we choose a highest weight vector
$0\not=v_\lambda\in L_\lambda[\lambda]$ and a lowest weight vector
$0\not=v_{w_0\lambda}\in L_\lambda[w_0\lambda]$, where $w_0\in S_{2n}$
is the longest Weyl group element, $w_0(i)=2n+1-i$ for all
$i$. Both $v_\lambda$ and $v_{w_0\lambda}$ are unique up to a
multiplicative constant.
\begin{defi}
Let $M$ be a left $\AAA_\sigma$-module and denote
$\hbox{Hom}_{\AAA_\sigma}(L_\lambda^\sigma,M)$ for the space of
$\AAA_\sigma$-intertwiners $L_\lambda^\sigma\rightarrow M$.
The maps $h_M^\lambda, l_M^\lambda:
\hbox{Hom}_{\AAA_\sigma}(L_\lambda^\sigma,M)
\rightarrow M$, defined by $h_M^\lambda(\phi)=\phi(v_\lambda)$ and
$l_M^\lambda(\phi)=\phi(v_{w_0\lambda})$, are called
the higher and lower expectation maps, respectively.
The values $h_M^\lambda(\phi)$ and $l_M^\lambda(\phi)$ are called
the higher and lower expectation value of
$\phi$.
\end{defi}

Define a map $\flat: P_{2n}\rightarrow \Lambda_n$ by
\[\lambda^\flat=(\lambda_1+\lambda_{2n},
\lambda_2+\lambda_{2n-1},\cdots,\lambda_n+\lambda_{n+1}),\qquad
\lambda\in P_{2n}.
\]
Observe that the kernel of $\flat$ equals $\Lambda_n^\natural$ and that
$(w_0\lambda)^\flat=\lambda^\flat$ for all $\lambda\in P_{2n}$.
If $M$ is a left $\AAA_\sigma$-module and
$\mu\in\Lambda_n$, then we define
\begin{equation}\label{weightA}
 M_{\mu}=\{m\in M \, | \, C_im=q^{\mu_i}m,\quad \forall\,
i=1,\ldots,n \}.
\end{equation}
\begin{prop}\label{expectationvalue}
Let $M$ be a left $\AAA_\sigma$-module and let $\lambda\in P_{2n}^+$.
The expectation maps $h_M^\lambda$ and $l_M^\lambda$
are injective, with image contained in $M_{\lambda^\flat}$.
\end{prop}
\begin{proof}
We prove the proposition for the expectation map $h_M^\lambda$,
the proof for $l_M^\lambda$ is analogous.
For $i=1,\ldots,n$ and $\phi\in
\hbox{Hom}_{\AAA_\sigma}(L_\lambda^\sigma,M)$ we have
\[C_ih_M^\lambda(\phi)=C_i\phi(v_\lambda)=\phi(C_iv_\lambda)=
q^{\lambda_i+\lambda_{2n+1-i}}h_M^\lambda(\phi),
\]
so $h_M^\lambda(\phi)\in M_{\lambda^\flat}$. Thus the image of
$h_M^\lambda$ is contained in $M_{\lambda^\flat}$.

Suppose $\phi\in \hbox{Hom}(L_\lambda^\sigma,M)$ is in the kernel of
$h_M^\lambda$, so $\phi(v_\lambda)=0$. We show that $\phi(v)=0$ for $v\in
L_\lambda[\nu]$ and $\nu\preceq \lambda$
by induction on the height
$\hbox{ht}(\lambda-\nu)\in\mathbb{Z}_{\geq 0}$, where
$\hbox{ht}(\lambda-\nu)=\sum_{i=1}^{2n-1}m_i$ if
$m_i$ are the unique positive integers such that
$\lambda-\nu=\sum_im_i\alpha_i$. So suppose that $\phi(v)=0$ for
all vectors $v\in L_\lambda[\nu]$ with $\nu\preceq\lambda$ and
$\hbox{ht}(\lambda-\nu)\leq N$.

For arbitrary $\nu\in P_{2n}$ we have $L_\lambda[\nu]=\{0\}$
unless $\nu\preceq\lambda$ and
\[L_\lambda[\nu]=\sum_{j=1}^{2n-1}\widetilde{y}_jL_\lambda[\nu+\alpha_j],
\]
where $\widetilde{y}_j=y_jK_{j+1}^{-1}K_{2n-j}^{-1}$.
To prove the induction step, it
thus suffices to show that $\phi(\widetilde{y}_jv)=0$ for all $j$ when
$v\in L_\lambda[\nu]$, $\nu\preceq\lambda$
and $\hbox{ht}(\lambda-\nu)=N$.
If $v$ is such a vector, then by the explicit form of
$B_j$ and $B_n^\sigma$ and by the induction hypothesis,
\begin{equation*}
\begin{split}
\phi(\widetilde{y}_jv)&=\phi(B_jv)=B_j\phi(v)=0,\\
\phi(\widetilde{y}_nv)&=\phi(B_n^\sigma v)=B_n^\sigma\phi(v)=0
\end{split}
\end{equation*}
for $j\in\{1,\ldots,2n-1\}\setminus \{n\}$, as desired.
\end{proof}
Since $L_\lambda^\sigma$ is $\AAA_\sigma$-semisimple
by Lemma \ref{smodule}, the previous proposition has the
following immediate consequence.
\begin{cor}\label{lh}
If $M$ is a simple $\AAA_\sigma$-module and
$\lambda\in P_{2n}^+$, then the number of summands
isomorphic to $M$ in the irreducible
decomposition of $L_\lambda^\sigma$
is bounded by $\hbox{Dim}(M_{\lambda^\flat})$.
\end{cor}
Recall the conditions on the representation labels
$\vec{\kappa}=(\kappa_1,\kappa_2,\kappa)$ from Remark
\ref{notationparameters}.
\begin{lem}\label{weakbranching}
Let $\lambda\in P_{2n}^+$. If both $V(\kappa,\kappa_1)_\sigma$ and
$V(\kappa_2)_\sigma$ occur as constituent in the irreducible
decomposition of the $\AAA_\sigma$-semisimple module
$L_\lambda^\sigma$, then $\lambda\in
(\Lambda_n^++\delta(\kappa,\kappa_1))^\natural$. For such $\lambda$,
the modules $V(\kappa,\kappa_1)_\sigma$ and $V(\kappa_2)_\sigma$
occur at most with multiplicity one in $L_\lambda^\sigma$.
\end{lem}
\begin{proof}
The weight spaces $M_\mu$ (see \eqref{weightA})
for the simple $\AAA_\sigma$-modules $M=V(\kappa,\kappa_1)_\sigma$
and $M=V(\kappa_2)_\sigma$ are at most one-dimensional by
Lemma \ref{explicit}.
Thus for any $\lambda\in P_{2n}^+$,
the number of summands isomorphic to $V(\kappa,\kappa_1)_\sigma$ or to
$V(\kappa_2)_\sigma$ in the irreducible
decomposition of $L_\lambda^\sigma$ is at most one
by Corollary \ref{lh}.

If $v\in V(\kappa_2)_\sigma$ then $C_iv=v$ for all $i$, hence
$V(\kappa_2)_\sigma$ can only occur as constituent in the irreducible
decomposition of $L_\lambda^\sigma$ when $\lambda\in P_{2n}^+\cap
\hbox{Ker}\bigl(\flat\bigr)=\bigl(\Lambda_n^+\bigr)^\natural$.

Fix $\mu\in\Lambda_n^+$ such that $V(\kappa,\kappa_1)_\sigma$ occurs
as constituent of the irreducible decomposition of
$L_{\mu^\natural}^\sigma$. It remains to show that this assumption imposes
the additional restrictions $\mu_n\geq \kappa_1$ and
$\mu_j-\mu_{j+1}\geq \kappa$ ($j=1,\ldots,n-1$) on $\mu$.
Throughout the proof,
we fix a nonzero intertwiner
$\phi\in \hbox{Hom}_{\AAA_\sigma}\bigl(L_{\mu^\natural}^\sigma,
V(\kappa,\kappa_1)_\sigma\bigr)$.

We start with proving $\mu_n\geq \kappa_1$. We use the notations
of Remark \ref{Rosremark}.
Standard $U_q(\mathfrak{g}\mathfrak{l}(2))$-representation theory
implies that
$V=U_q(\mathfrak{g}\mathfrak{l}(2))^{(n)}v_{\mu^\natural}\subseteq
L_{\mu^\natural}$ is the $(2\mu_n+1)$-dimensional irreducible
representation of $U_q(\mathfrak{g}\mathfrak{l}(2))^{(n)}\simeq
U_q(\mathfrak{g}\mathfrak{l}(2))$ with highest weight $(\mu_n,-\mu_n)$ and
highest weight vector $v_{\mu^\natural}$, and that
$\widehat{B}_n^\sigma|_V\in\hbox{End}_{\mathbb{C}}(V)$ is semisimple
with simple spectrum
\[s_l=\frac{q^{-\sigma+2l}-q^{\sigma-2l}}{q-q^{-1}},
\qquad l=-\mu_n,\ldots,\mu_n-1,\mu_n.
\]
Observe furthermore that $V$ is
$x_\sigma$-stable.
Since $x_\sigma B_n^\sigma x_\sigma^{-1}=\widehat{B}_n^\sigma$,
it follows that $B_n^\sigma|_V\in\hbox{End}_{\mathbb{C}}(V)$ is
semisimple with simple spectrum
$\{s_l\, | \, l=-\mu_n,\ldots,\mu_n-1,\mu_n\}$.
Proposition \ref{expectationvalue} and the fact that $C_i$ centralizes
$U_q(\mathfrak{g}\mathfrak{l}(2))^{(n)}$ for $i=1,\ldots,n$ imply that
$\phi(V)=\widetilde{V}(\kappa,\kappa_1)_\sigma$ with
$\widetilde{V}(\kappa,\kappa_1)_\sigma$ the one-dimensional space
defined by \eqref{zeroweight}. On the other hand,
Lemma \ref{explicit} shows that
$\widetilde{V}(\kappa,\kappa_1)_\sigma$ is the eigenspace of
$\rho(\kappa,\kappa_1)_\sigma(B_n^\sigma)$ with eigenvalue
$s_{-\kappa_1}$,
hence $s_{-\kappa_1}$ must be in the spectrum of $B_n^\sigma|_V$. Consequently,
$\mu_n\geq\kappa_1$.

Next we show that $\mu_j-\mu_{j+1}\geq \kappa$ for $j=1,\ldots,n-1$.
Fix $j\in\{1,\ldots,n-1\}$.
Since
$0\not=\phi(v_{\mu^\natural})\in\widetilde{V}(\kappa,\kappa_1)_\sigma$
by Proposition \ref{expectationvalue}, we
conclude from the explicit description of $V(\kappa,\kappa_1)_\sigma$
(see Lemma \ref{explicit}) that
\[
\phi(B_j^{\kappa}v_{\mu^\natural})=
\rho(\kappa,\kappa_1)_\sigma(B_j^{\kappa})\phi(v_{\mu^\natural})\not=0.
\]
This implies that
$B_j^{\kappa}v_{\mu^\natural}\not=0$ in $L_{\mu^\natural}$.
We write $B_j=\widetilde{y}_j+\widetilde{x}_{2n-j}$ with
\[\widetilde{y}_j=y_jK_{j+1}^{-1}K_{2n-j}^{-1},\qquad
\widetilde{x}_{2n-j}=K_j^{-1}x_{2n-j}K_{2n-j}^{-1}.
\]
Since $v_{\mu^\natural}\in L_{\mu^\natural}$ is a highest weight vector,
the commutation relation $\widetilde{x}_{2n-j}\widetilde{y}_{j}=
q^2\widetilde{y}_{j}\widetilde{x}_{2n-j}$ in $\ug$ implies
\[\widetilde{y}_j^\kappa
v_{\mu^\natural}=B_j^\kappa v_{\mu^\natural}\not=0,
\]
and consequently $y_j^\kappa v_{\mu^\natural}\not=0$.
Standard arguments from
$U_q(\mathfrak{g}\mathfrak{l}(2))\simeq \mathbb{C}\langle
x_j,y_j,K_j^{\pm 1}, K_{j+1}^{\pm 1}\rangle$ representation theory
imply that $\mu_j-\mu_{j+1}$ is the largest positive integer $m$
for which $y_j^mv_{\mu^\natural}\not=0$ in $L_{\mu^\natural}$.
We conclude that $\mu_j-\mu_{j+1}\geq \kappa$, as desired.
\end{proof}

In, e.g., \cite{N} and \cite{NDS}
highest weight considerations are used to determine branching
rules for quantum symmetric pairs.
These highest weight arguments can be formalized and generalized
in the present setting as follows.
\begin{lem}\label{highcontri}
Let $\lambda\in P_{2n}^+$ and let
$\{0\}\not=M\subseteq L_\lambda^\sigma$ be an $\AAA_\sigma$-submodule.
Then
\[M\not\subseteq \bigoplus_{\nu\prec\lambda}L_\lambda[\nu],
\qquad M\not\subseteq \bigoplus_{\nu\succ w_0\lambda}L_\lambda[\nu].
\]
In other words, there exists a vector $m_+\in M$
\textup{(}resp. $m_-\in M$\textup{)}
whose decomposition in $U^0$-weights has a nonzero
highest \textup{(}resp. lowest\textup{)} weight contribution.
\end{lem}
\begin{proof}
Suppose $M\subseteq L_\lambda^\sigma$ is an $\AAA_\sigma$-submodule and
\begin{equation}\label{condition}
M\subseteq \bigoplus_{\nu\prec\lambda}L_\lambda[\nu].
\end{equation}
Let $M^\perp$ be the orthocomplement of $M$ in $L_\lambda$
with respect to the $*$-unitary scalar product $\langle
\cdot,\cdot\rangle_\lambda$ on $L_\lambda$.
Then $M^\perp$ is an $\AAA_\sigma$-module complement of
$M$ in $L_\lambda^\sigma$.
Let $\pi\in \hbox{Hom}_{\AAA_\sigma}(L_\lambda^\sigma, M)$ be the
projection onto $M$ along $M^\perp$.
By the $*$-selfadjointness of the $K_i$'s, we have
$L_\lambda[\nu]\perp L_\lambda[\nu^\prime]$ for $\nu\not=\nu^\prime$.
Condition \eqref{condition} thus implies $L_\lambda[\lambda]\subseteq
M^\perp$. Hence the higher expectation value $h_M^\lambda(\pi)$ is zero.
Proposition \ref{expectationvalue}
then implies $\pi\equiv 0$, hence $M=\{0\}$.
The second statement is proved similarly.
\end{proof}

\begin{rema}
For an arbitrary quantum symmetric pairs $(U, B_\theta)$ in the sense
of Letzter (see \cite[\S 7]{L4} for a complete list), one can make a similar
definition of higher and lower expectation maps and of expectation values.
Following the reasoning of this section, one can establish the analogues
of Proposition \ref{expectationvalue},
Corollary \ref{lh} and Lemma \ref{highcontri} in this general setup
(cf. the proof of \cite[Lem. 5.6]{L5}).
\end{rema}

\subsection{Deformation arguments}

In this subsection we complete the proof of Proposition~\ref{branching}
by using deformation arguments. Without loss of generality
we may take $q=e^h$ to be a formal deformation parameter.
We use part four of Kassel's
book \cite{Kas} as main reference for facts and notations related to
topological $\mathbb{C}[[h]]$-modules and topological
$\mathbb{C}[[h]]$-algebras.

Let $U_h(\mathfrak{g})$ be the
topological version of the quantized
universal enveloping algebra $\ug$, with topological generators
$E_{i,i}$, $x_j$ and $y_j$ ($i\in\{1,\ldots,2n\}$, $j\in \{1,\ldots,2n-1\}$),
where $K_i=e^{hE_{i,i}}$.
Let $U_h(\mathfrak{sl}(2n))$ be the subalgebra of $U_h(\mathfrak{gl}(2n))$
topologically generated by the elements  $E_{j,j}-E_{j+1,j+1}$,
$x_j$ and $y_j$ ($j\in\{1,\ldots,2n-1\}$). Note that
$U_h(\mathfrak{g})$ is the
central extension of $U_h(\mathfrak{sl}(2n))$
by the central element $E_{1,1}+E_{2,2}+\cdots +E_{2n,2n}$.

The topological version of $\mathcal{A}_\sigma$ is the subalgebra
of $U_h(\mathfrak{g})$ topologically generated by
$B_j,B_n^\sigma$ and $E_{i,i}+E_{2n+1-i,2n+1-i}$ for
$j\in\{1,\ldots,2n-1\}\setminus \{n\}$ and $i\in\{1,\ldots,n\}$.
We denote by $\overline{\mathcal A}_\sigma$  the subalgebra of
$U_h(\mathfrak{sl}(2n))$ topologically generated by $B_jC_{j+1}$,
$B_n^{\sigma}C_n$ and $E_{i,i}-E_{i+1,i+1}-
E_{2n-i,2n-i}+E_{2n+1-i,2n+1-i}$ for $j\in
\{1,\dots,2n-1\}\setminus \{n\}$ and $i\in \{1,\ldots,n-1\}$. Note
that $\mathcal{A}_\sigma$ is the central extension of
$\overline{\mathcal A}_\sigma$ by the central element
$E_{1,1}+E_{2,2}+\cdots +E_{2n,2n}\in U_h(\mathfrak{g})$.

We fix $\lambda=\bigl(\mu+\delta(\kappa,\kappa_1)\bigr)^\natural$
with $\mu\in\Lambda_n^+$.
To complete the proof of Proposition \ref{branching}, it suffices to
show that the irreducible $\AAA_\sigma$-module
$V(\kappa,\kappa_1)_\sigma$ occurs as consituent of $L_\lambda^\sigma$.
Since the central element $E_{1,1}+\cdots +E_{2n,2n}$ acts as zero on
the topological versions of
$V(\kappa,\kappa_1)_\sigma$ and $L_\lambda$, it suffices to show
that $V(\kappa,\kappa_1)_\sigma$ is a constituent of $L_\lambda^\sigma$
when viewed as $\overline{\mathcal{A}}_\sigma$-modules.

The formal computations after Definition 1.5 and \cite[Thm. 7.5]{L3}
imply that the algebra $\overline{\mathcal{A}}_\sigma$ is a deformation of
the universal enveloping algebra $U(\overline{\mathfrak{k}})$ with
\[\overline{\mathfrak{k}}=\hbox{Ad}(g)\bigl(\mathfrak{k}\bigr)\cap
\mathfrak{sl}(2n)=\overline{\mathfrak{k}}_{ss}\oplus\mathbb{C}Z,
\]
where $g$ is given by \eqref{forg},
$\overline{\mathfrak{k}}_{ss}=\hbox{Ad}(g)\bigl(\mathfrak{sl}(n)\oplus
\mathfrak{sl}(n)\bigr)$, and
$Z\in \overline{\mathfrak{k}}$ the central element
\[
Z=\hbox{Ad}(g)(Z^\prime)=-\frac{1}{2n}\sum_{j=1}^{2n}E_{j,2n+1-j}
\]
where
\[
Z^\prime=-\frac{1}{2n}
\bigl(E_{1,1}+E_{2,2}+\cdots +E_{n,n}-E_{n+1,n+1}-E_{n+2,n+2}-\cdots -
E_{2n,2n}\bigr).
\]
Since $\overline{\mathfrak{k}}$ is reductive with one-dimensional
center, we have $H^2(\overline{\mathfrak{k}},U(\overline{\mathfrak{k}}))=
\{0\}$. By \cite[Thm. XVIII.2.2]{Kas} there exists a topological
algebra isomorphism
\[
\beta_\sigma: \overline{\mathcal{A}}_\sigma\rightarrow
U(\overline{\mathfrak{k}})[[h]]
\]
which is the identity modulo $h$. We define
\[\overline{\mathcal{A}}_\sigma^{ss}=\beta_\sigma^{-1}\bigl(
U(\overline{\mathfrak{k}}_{ss})[[h]]\bigr),\qquad
Z_\sigma=\beta_\sigma^{-1}(Z),
\]
then $\overline{\mathcal{A}}_\sigma^{ss}$ is a deformation of
$U(\overline{\mathfrak{k}}_{ss})$, $Z_\sigma$ is central in
$\overline{\mathcal{A}}_\sigma$, and
$\overline{\mathcal{A}}_\sigma$ is topologically generated by
$\overline{\mathcal{A}}_\sigma^{ss}$ and $Z_\sigma$. We denote by
$(\beta_\sigma^{-1})^*(V(\kappa,\kappa_1)_\sigma)$ the module
$V(\kappa,\kappa_1)_\sigma$ viewed as module over
$U(\overline{\mathfrak{k}})[[h]]$ via the isomorphism
$\beta_\sigma: \overline{\mathcal{A}}_\sigma \rightarrow
U(\overline{\mathfrak{k}})[[h]]$ (we will use similar notations
to indicate the pull back of an action by some
algebra morpism). The way we constructed the
module $V(\kappa,\kappa_1)_\sigma$ in \S 1 immediately implies
that modulo $h$, we reobtain the irreducible, finite dimensional
$\bigl(\mathfrak{gl}(n)\times \mathfrak{gl}(n)\bigr)\cap
\mathfrak{sl}(2n)\simeq \overline{\mathfrak{k}}$ module
$V^{clas}(\kappa,\kappa_1)$ of highest weight
\begin{equation}\label{hw}
\bigl((-\kappa_1)^n,\bigl(\kappa_1+(n-1)\kappa,(\kappa_1-\kappa)^{n-1}\bigr)
\bigr)\in P_n^+\times P_n^+.
\end{equation}
Observe that the central element $Z$ acts as multiplication by $\kappa_1$
on $V^{clas}(\kappa,\kappa_1)$, while $V^{clas}(\kappa,\kappa_1)$
is independent of $\kappa_1$
when viewed as module over $\overline{\mathfrak{k}}_{ss}$.
Since $\overline{\mathfrak{k}}_{ss}$ is semisimple, we
conclude that $(\beta_\sigma^{-1})^*(V(\kappa,\kappa_1)_\sigma)$, viewed as
$U(\overline{\mathfrak{k}}_{ss})[[h]]$-module,
is isomorphic to the $\mathbb{C}[[h]]$-linear extension of the
$\overline{\mathfrak{k}}_{ss}$-module
$V^{clas}(\kappa,\kappa_1)$.
\begin{lem}\label{lemp}
There exists a $\overline{\mathcal{A}}_{\sigma}$-submodule
of $\overline{V}(\kappa,\kappa_1)_\sigma\subset L_\lambda^\sigma$
such that\\
{\bf (i)} $\overline{V}(\kappa,\kappa_1)_\sigma\simeq
V(\kappa,\kappa_1)_\sigma$ as $\overline{\mathcal{A}}_\sigma^{ss}$-modules.\\
{\bf (ii)} The central element $Z_\sigma$ acts as multiplication by a scalar
${\overline{\xi}}_\sigma\in\mathbb{C}[[h]]$ on
$\overline{V}(\kappa,\kappa_1)_\sigma$.\\
{\bf (iii)} Modulo $h$, $\overline{V}(\kappa,\kappa_1)_\sigma$ is
the $\overline{\mathfrak{k}}$-module $V^{clas}(\kappa,\kappa_1)$.
\end{lem}
\begin{proof}
Since $U_h(\mathfrak{sl}(2n))$ is a deformation of the universal enveloping
algebra of the simple Lie algebra $\mathfrak{sl}(2n)$,
there exists by \cite[Thm. XVIII.2.2]{Kas}
a topological algebra
isomorphism
\[\alpha: U_h(\mathfrak{sl}(2n))\rightarrow U(\mathfrak{sl}(2n))[[h]]
\]
which is the identity modulo $h$. The natural embedding
$\iota_\sigma: \overline{\mathcal{A}}_\sigma\rightarrow U_h(\mathfrak{sl}(2n))$
of topological algebras can thus be pulled back to obtain an
embedding
\[\gamma_\sigma=\alpha\circ\iota_\sigma\circ\beta_\sigma^{-1}:
U(\overline{\mathfrak{k}})[[h]]\rightarrow U(\mathfrak{sl}(2n))[[h]]
\]
of topological algebras, which is the identity modulo $h$.
Since $\overline{\mathfrak{k}}_{ss}$ is semisimple,
we have $H^1(\overline{\mathfrak{k}}_{ss},U(\mathfrak{sl}(2n)))=\{0\}$,
hence there exists an invertible element
$F_\sigma\in U(\mathfrak{sl}(2n))[[h]]$ which is $1$ modulo $h$
such that
\begin{equation}\label{trivpart}
\gamma_\sigma(X)=F_\sigma X F_{\sigma}^{-1}\qquad
\forall X\in U(\overline{\mathfrak{k}}_{ss})[[h]]\subset
U(\overline{\mathfrak{k}})[[h]]\subset
U(\mathfrak{sl}(2n))[[h]]
\end{equation}
by \cite[Thm. XVIII.2.1]{Kas}.

Let $L_\lambda^{clas}$ be the irreducible $U(\mathfrak{s}\mathfrak{l}(2n))$-
module corresponding to the highest weight $\lambda$.
We denote $V\subset L_\lambda^{clas}$ for its
$V^{clas}(\kappa,\kappa_1)$-isotypical component viewed as
$U(\overline{\mathfrak{k}}_{ss})$-module.
The central element $Z\in \overline{\mathfrak{k}}$
acts semisimply on $V$ and the corresponding $Z$-eigenspaces
$V_r$ of $V$ (with $r$ the $Z$-eigenvalue) yield the
$V^{clas}(\kappa,r)$ -isotypical components of $L_\lambda^{clas}$
viewed now as $\overline{\mathfrak{k}}$-module.
By the classical branching rules \cite[Lemma 3]{O} for the symmetric pair
$(\mathfrak{g},\mathfrak{gl}(n)\times \mathfrak{gl}(n))$ we have
$V_{\kappa_1}\simeq V^{clas}(\kappa,\kappa_1)$.

Let $L_\lambda^{clas}[[h]]$ be the
$U(\mathfrak{s}\mathfrak{l}(2n))[[h]]$-module obtained by
$\mathbb{C}[[h]]$-linear extension of the
$U(\mathfrak{s}\mathfrak{l}(2n))$-module structure on $L_\lambda^{clas}$.
Then  $L_\lambda^\sigma\simeq
(\alpha\circ \iota_\sigma)^*(L_\lambda^{clas}[[h]])$
as topological $\overline{\mathcal{A}}_\sigma$-module.
Now consider the $U(\overline{\mathfrak{k}})[[h]]$-module
\[(\beta_\sigma^{-1})^*(L_\lambda^\sigma)=
\gamma_\sigma^*(L_\lambda^{clas}[[h]]).
\]
By \eqref{trivpart}, $F_{\sigma}V[[h]]$ is the
$V^{clas}(\kappa,\kappa_1)[[h]]$-isotypical component
of the $U(\overline{\mathfrak{k}}_{ss})[[h]]$-module
$\gamma_\sigma^*(L_\lambda^{clas}[[h]])$.
The proof is now completed by observing that
$F_{\sigma}V[[h]]$ contains a nonzero
$Z$-eigenspace $\overline{V}(\kappa,\kappa_1)_\sigma$ for some eigenvalue
$\overline{\xi}_\sigma\in \mathbb{C}[[h]]$ satisfying
$\overline{\xi}_\sigma=\kappa_1$ modulo $h$.
\end{proof}

The central element $Z_\sigma\in \overline{\mathcal{A}}_\sigma$
acts by a scalar $\xi_\sigma\in \mathbb{C}[[h]]$ on
$V(\kappa,\kappa_1)_\sigma$. To complete the arguments it thus
remains to prove the following lemma.
\begin{lem}
$\xi_\sigma=\overline{\xi}_\sigma$ in $\mathbb{C}[[h]]$.
\end{lem}
\begin{proof}
We use the fact that
$B_n^{\sigma}C_n\in \overline{\mathcal{A}}_\sigma$ is conjugate
to the Cartan type element
\[
\overline{B}_n^\sigma=\frac{q^{-\sigma}K_nK_{n+1}^{-1}-
q^{\sigma}K_n^{-1}K_{n+1}}{q-q^{-1}}
\]
in $U_h(\mathfrak{sl}(2n))$, see \cite{Ro} and Remark
\ref{Rosremark}. We first show that $B_n^\sigma C_n$ acts semisimply on
$\overline{V}(\kappa,\kappa_1)_\sigma$, with spectrum contained in
\[\left\{t_l=\frac{q^{-\sigma+l}-q^{\sigma-l}}{q-q^{-1}}\, | \,
l\in\mathbb{Z} \right\}.
\]
Indeed, $B_n^\sigma C_n$ acts semisimply on
$L_\lambda^\sigma$, with spectrum
contained in $\{t_l\}_{l\in \mathbb{Z}}$.
Hence any vectors $v\in \overline{V}(\kappa,\kappa_1)_\sigma$ can
be written as $v=\sum_l v_l$ with $v_l\in \ker (B_n^\sigma
C_n-t_l)\subseteq L_\lambda^\sigma$. Since $B_n^\sigma C_n$ preserves
$\overline{V}(\kappa,\kappa_1)_\sigma$, we have
$(B_n^\sigma C_n)^k(v)=\sum_l t_l^k
v_l\in\overline{V}(\kappa,\kappa_1)_\sigma$ for
all $k\in\mathbb{Z}_{\geq 0}$.
Simple calculations involving the
Vandermonde determinant shows that there exists
$c_k^l\in\mathbb{C}[[h]]$ such that $v_l=\sum_k c^l_k (B_n^\sigma
C_n)^k v$ for all $l$, hence $v_l\in
\overline{V}(\kappa,\kappa_1)_\sigma$ for all $l$.
By Lemma \ref{explicit}, the same statement about the spectrum is valid
for the action of $B_n^\sigma C_n$ on $V(\kappa,\kappa_1)_\sigma$.

Lemma \ref{lemp} and the fact that the eigenvalues
$t_l|_{h=0}=\sigma-l$ ($l\in \mathbb{Z}$) remain seperated modulo
$h$ imply that the dimensions over $\mathbb{C}[[h]]$
of the $B_n^\sigma C_n$-eigenspaces of $V(\kappa,\kappa_1)_\sigma$
and $\overline{V}(\kappa,\kappa_1)_\sigma$ corresponding to a given
eigenvalue $t_l$, are the same. In particular, we have
\begin{equation}\label{trace=}
\hbox{Tr}_{V(\kappa,\kappa_1)_\sigma}(B_n^\sigma C_n)=
\hbox{Tr}_{\overline{V}(\kappa,\kappa_1)_\sigma}(B_n^\sigma C_n).
\end{equation}
Next we expand $B_n^{\sigma}C_n$ in its semisimple and central parts,
\[B_n^\sigma C_n=\sum_{j=0}^{\infty}a_jZ_\sigma^j
\]
in $\overline{\mathcal{A}}_\sigma$, with
$a_j\in\overline{\mathcal{A}}_\sigma^{ss}$. Then \eqref{trace=}
leads to the equality
\begin{equation}\label{equalityxi}
\sum_{j=0}^{\infty}b_j\xi_\sigma^j
=\sum_{j=0}^{\infty}b_j\overline{\xi}_\sigma^j
\end{equation}
as formal power series in $h$, with
\[b_j=\hbox{Tr}_{V(\kappa,\kappa_1)_\sigma}(a_j)\in\mathbb{C}[[h]].
\]
Now the computations after Definition 1.5 show that
\[
\hbox{Ad}(g)\bigl(B_n^{\sigma}C_n|_{h=0}\bigr)=
E_{n,n}-E_{n+1,n+1}-\sigma \,1=a_0|_{h=0}+a_1|_{h=0}Z^\prime
\]
with
\begin{equation*}
\begin{split}
a_0|_{h=0}&=-\frac{1}{n}E_{1,1}-\cdots -\frac{1}{n}E_{n-1,n-1}
+\frac{n-1}{n}E_{n,n}\\
&+\frac{1-n}{n}E_{n+1,n+1}+\frac{1}{n}E_{n+2,n+2}+
\cdots +\frac{1}{n}E_{2n,2n}-\sigma\,1
\end{split}
\end{equation*}
and $a_1|_{h=0}=-2$. Thus $b_1|_{h=0}\not=0$ and
$b_j|_{h=0}=0$ for $j\geq 2$. Now we can expand \eqref{equalityxi}
in powers of $h$, and compare coefficients of $h^k$ on each side.
An easy induction argument using $b_1|_{h=0}\not=0$ and
$b_j|_{h=0}=0$ for $j\geq 2$ then leads to $\xi_\sigma=\overline{\xi}_\sigma$
(this may be viewed as a formal version of the inverse function theorem).
\end{proof}


\section{Zonal spherical functions}

The purpose of this section is to establish Theorem
\ref{maintheorem} for $\vec{\kappa}=\vec{0}=(0,0,0)$. This is done
by translating the results for zonal spherical functions on
quantum Grassmannians proved in \cite{NDS} and \cite{DS} to the
present setting. In \cite{NDS} and \cite{DS}, extensive use is made
of the $L$-operators of $\ug$. Roughly speaking, the $L$-operators
$L^{\pm}$ for $\ug$ are the $2n\times 2n$ matrices with entries in
$\ug$ obtained from the universal $R$-matrix of $\ug$ by applying
the vector representation of $\ug$ to one of its components. We
refer to Letzter \cite[\S 6]{L1} and
Noumi \cite{N} for precise definitions, here
we only recall those properties of the $L$-operators which are
used in this paper. The matrix $L^+$ (resp. $L^-$) is upper (resp.
lower) triangular. If $l_{ij}^{\pm}\in \ug$ is the $(i,j)$th
matrix coefficient of $L^{\pm}$, then
\begin{equation}\label{lcoeff}
\begin{split}
l_{ii}^+&=K_i,\qquad\qquad\qquad\qquad l_{ii}^-=K_i^{-1},\\
l_{j,j+1}^+&=(q-q^{-1})K_jy_j,\qquad
l_{j+1,j}^-=(q^{-1}-q)x_jK_j^{-1}
\end{split}
\end{equation}
for $i=1,\ldots,2n$ and $j=1,\ldots,2n-1$. More generally,
$l_{ij}^+$ (resp. $l_{ij}^-$) is an analogue of the root vector of
$\mathfrak{g}$ corresponding to the root $\epsilon_i-\epsilon_j$
(resp. $\epsilon_j-\epsilon_i$). Furthermore,
\[
\Delta(l_{ij}^{\pm})=\sum_{k=1}^{2n}l_{ik}^{\pm}\otimes
l_{kj}^{\pm},\qquad
\epsilon(l_{ij}^{\pm})=\delta_{i,j},\qquad
(l_{ij}^{\pm})^*=S(l_{ji}^{\mp})
\]
and $S(L^\pm)=(L^\pm)^{-1}$, where the antipode is applied
componentwise.

Define a scalar valued $2n\times 2n$ matrix $J^\sigma$, depending on
an auxiliary parameter $\sigma\in\mathbb{R}$, by
\[J^\sigma=(1-q^{2\sigma})\sum_{k=1}^nE_{k,k}-
q^{\sigma}\sum_{k=1}^{2n}E_{k,2n+1-k}.
\]
The matrix $J^\sigma$ is a solution of a reflection equation, see
\cite[Prop. 2.2]{NDS} or \cite[\S 6]{DS}
for details. Observe furthermore that $J^\sigma$ is invertible, with inverse
\[(J^\sigma)^{-1}=(1-q^{-2\sigma})\sum_{k=1}^{n}E_{2n+1-k,2n+1-k}-
q^{-\sigma}\sum_{k=1}^{2n}E_{k,2n+1-k}.
\]
\begin{defi}[\cite{NDS}, \cite{DS}]
Let $\mathfrak{k}_\sigma\subset \ug$ be the linear subspace spanned
by the $4n^2$ matrix coefficients of
$L^+J^{\sigma}-J^{\sigma}L^-$.
\end{defi}
Observe that $\mathfrak{k}_\sigma\subset \ug$ is a two-sided
coideal, i.e. $\Delta(\mathfrak{k}_\sigma)\subseteq
\mathfrak{k}_\sigma\otimes \ug+ \ug\otimes \mathfrak{k}_\sigma$
and $\epsilon(\mathfrak{k}_\sigma)=0$. Note furthermore that
$\mathfrak{k}_\sigma^*=S(\mathfrak{k}_\sigma)$ since $J^\sigma$ is symmetric,
so $\mathfrak{k}_\sigma$ is {\it not} $*$-stable but
$\mathfrak{k}_\sigma$ is $\omega$-stable, where
$\omega:\ug\rightarrow \ug$ is the involutive, anti-linear algebra
isomorphism $\omega=*\circ S$. Note that the involution $\omega$ is related to
the $*$-structure on $\mathbb{C}_q[G]$ by
\begin{equation}\label{omegaster}
f^*(X)=\overline{f(\omega(X))},\qquad f\in
\mathbb{C}_q[G],\,\, X\in \ug.
\end{equation}
It is convenient to write
\begin{equation}\label{phi}
\omega_\delta(X)=K^{-\delta}\omega(X)K^{\delta},\qquad X\in \ug.
\end{equation}
Note that $\omega_\delta$ is an anti-linear algebra involution of
$\ug$, which acts on
on the algebraic generators of $\ug$ by
\[\omega_\delta(K_i^{\pm 1})=K_i^{\mp 1},\qquad
\omega_\delta(x_j)=-y_j,\qquad
\omega_\delta(y_j)=-x_j
\]
for $i=1,\ldots,2n$ and $j=1,\ldots,2n-1$.

For a subspace $U\subset \ug$, let $I(U)\subseteq \ug$ be the left
ideal of $\ug$ generated by $U$ and let $A(U)\subseteq \ug$
be the unital subalgebra generated by $I(U)$. Applied to
$\mathfrak{k}_\sigma$ and $\mathfrak{k}_\sigma^*$, we thus obtain
two unital subalgebras $A(\mathfrak{k}_\sigma)$ and
$A(\mathfrak{k}_\sigma^*)$ of $\ug$.
In the following lemma we link the algebras
$A(\mathfrak{k}_\sigma)$ and $A(\mathfrak{k}_\sigma^*)$ to $\AAA_\sigma$.

\begin{lem}\label{linksubalgebras}
{\bf (i)}
$\AAA_\sigma\cup \omega(\AAA_\sigma)\subseteq A(\mathfrak{k}_\sigma)$.

{\bf (ii)} $\omega_\delta(\AAA_\sigma)\subseteq A(\mathfrak{k}_\sigma^*)$.
\end{lem}
\begin{proof}
{\bf (i)} Note that $A(\mathfrak{k}_\sigma)$ is $\omega$-stable,
because $\mathfrak{k}_\sigma$ is $\omega$-stable. It thus suffices
to prove that $\AAA_\sigma\subseteq A(\mathfrak{k}_\sigma)$.

The left ideal $I(\mathfrak{k}_\sigma)$ contains the
matrix coefficients of the $\ug$-valued matrix
\[(t_{i,j})_{i,j}=S(L^+)\bigl(L^+J^{\sigma}-J^{\sigma}L^-\bigr)=
J^{\sigma}-S(L^+)J^\sigma L^-.
\]
By direct computations one verifies that
\begin{equation*}
\begin{split}
&t_{i,2n+1-i}=-q^{\sigma}+q^\sigma C_i^{-1},\qquad
t_{j,2n-j}=q^{\sigma}(q^{-1}-q)B_j,\\
&t_{n,n}=q^{\sigma}(q^{-1}-q)B_n^\sigma+1-q^{2\sigma}
\end{split}
\end{equation*}
for $i\in\{1,\ldots,n\}$ and $j\in\{1,\ldots,2n-1\}\setminus \{n\}$.
The left ideal $I(\mathfrak{k}_\sigma)$ also contains the matrix
coefficients of
\[S(L^-)(J^{\sigma})^{-1}(L^+J^{\sigma}-J^{\sigma}L^-)(J^{\sigma})^{-1}=
S(L^-)(J^{\sigma})^{-1}L^+-(J^{\sigma})^{-1}.
\]
Computing the coefficients of the antidiagonal shows that
$q^{-\sigma}-q^{-\sigma}C_i\in
I(\mathfrak{t}_\sigma)$ for $i\in\{1,\ldots,n\}$. Consequently, all algebraic
generators of $\AAA_\sigma$ are contained in $A(\mathfrak{k}_\sigma)$,
hence $\AAA_\sigma\subseteq A(\mathfrak{k}_\sigma)$.

{\bf (ii)}
Observe that
\begin{equation}\label{twistedgenerators}
\begin{split}
\omega_\delta(C_i^{\pm 1})&=C_i^{\mp 1},\\
\omega_\delta(B_j)&=-q\bigl(K_{2n-j}y_{2n-j}K_j+K_{2n-j}K_{j+1}x_j\bigr),\\
\omega_\delta(B_n^\sigma)&=-\left(K_ny_nK_n+K_nK_{n+1}x_n+
\left(\frac{q^{\sigma}-q^{-\sigma}}{q-q^{-1}}\right)K_n^2\right)
\end{split}
\end{equation}
for $i\in\{1,\ldots,n\}$ and $j\in\{1,\ldots,2n-1\}\setminus\{n\}$
are algebraic generators of $\omega_\delta(\AAA_\sigma)$.
Using similar arguments as in the proof of {\bf (i)}, it can be shown that
all these generators are contained in $A(\mathfrak{k}_\sigma^*)=
A(S(\mathfrak{k}_\sigma))$.
\end{proof}
For a left $\ug$-module $M$ and an unital subalgebra $A\subset \ug$
we define the subspace of $A$-invariant vectors in $M$ by
\[
M^{A}=\{m\in M \, | \, am=\epsilon(a)m\quad \forall\,
a\in A\}.
\]

\begin{lem}\label{onesimple}
For the $\ug$-subalgebras $A=\AAA_\sigma$, $A=\omega(\AAA_\sigma)$
and $A=\omega_\delta(\AAA_\sigma)$, we have

{\bf (i)}  $\hbox{Dim}\bigl(L_\lambda^A\bigr)=0$ if
$\lambda\in P_{2n}^+\setminus (\Lambda_n^+)^\natural$.

{\bf (ii)} $\hbox{Dim}\bigl(L_\lambda^A\bigr)\leq 1$ if $\lambda\in
(\Lambda_n^+)^\natural$.

\end{lem}
\begin{proof}
The statement for $A=\AAA_\sigma$ follows from the results in \S 2.

For $\omega(\AAA_\sigma)$ and $\omega_\delta(\AAA_\sigma)$, we first
remark that any simple, finite dimensional $\ug$-module
$L_\lambda$ ($\lambda\in P_{2n}^+$) is semisimple as module over the
subalgebras $\omega(\AAA_\sigma)$ and $\omega_\delta(\AAA_\sigma)$.
Clearly it suffices to prove this for $\omega_\delta(\AAA_\sigma)$.
In this case we first note that
\begin{equation}\label{starcommute}
\omega_\delta(X^*)=\bigl(\omega_\delta(X)\bigr)^*,\qquad
X\in\ug,
\end{equation}
which follows easily from the fact that
\begin{equation}\label{Sinner}
S^2(X)=K^{-2\delta}XK^{2\delta},\qquad X\in \ug.
\end{equation}
We conclude that $\omega_\delta(\AAA_\sigma)$ is a $*$-subalgebra
of $\ug$, hence $L_\lambda$ is $\omega_\delta(\AAA_\sigma)$-semisimple.

Straightforward adjustments of the arguments in \S 2 now lead
to the construction of expectation values and to the analogue of Proposition
\ref{expectationvalue} for the algebras $\omega(\AAA_\sigma)$ and
$\omega_\delta(\AAA_\sigma)$. Since the elements $C_i$ ($i=1,\ldots,n$)
are in $\omega(\AAA_\sigma)$ and $\omega_\delta(\AAA_\sigma)$, the
lemma now follows for these two subalgebras in the same way as for
$\AAA_\sigma$.
\end{proof}
If the algebra $A\subseteq \ug$ is of the form
$A=A(U)$ for some vector space $U\subset \ug$ and
$\epsilon|_U\equiv 0$, then $M^{A(U)}$ consists of the vectors
$m\in M$ which are annihilated by $u\in U$. In this situation,
$M^{A(U)}$ is called the subspace of $U$-fixed vectors, cf. \cite{NDS},
\cite{DS}.
In particular, this applies to
$U=\mathfrak{k}_\sigma$ and $U=\mathfrak{k}_\sigma^*$.
We now recall the following result, see \cite[Thm. 2.6]{NDS} and
\cite[Thm. 6.6]{DS}.
\begin{prop}\label{Noumi}
For $U=\mathfrak{k}_\sigma$ and $U=\mathfrak{k}_\sigma^*$ we have
\begin{equation*}
\begin{split}
\hbox{Dim}\bigl(L_\lambda^{A(U)}\bigr)&=0\,\,\, \hbox{ if }\,\,\,
\lambda\in P_{2n}^+
\setminus (\Lambda_n^+)^\natural,\\
\hbox{Dim}\bigl(L_\lambda^{A(U)}\bigr)&=1\,\,\, \hbox{ if }\,\,\,
\lambda\in (\Lambda_n^+)^\natural.
\end{split}
\end{equation*}
\end{prop}
Combined with Lemma \ref{linksubalgebras} and Lemma \ref{onesimple},
we obtain
\begin{lem}\label{zonal}
For $\lambda\in P_{2n}^+$ we have

{\bf (i)} $L_\lambda^{\AAA_\sigma}=L_\lambda^{\omega(\AAA_\sigma)}=
L_{\lambda}^{A(\mathfrak{k}_\sigma)}$.

{\bf (ii)} $L_\lambda^{\omega_\delta(\AAA_\sigma)}=
L_{\lambda}^{A(\mathfrak{k}_\sigma^*)}$.\\
\end{lem}
\begin{proof}
By Lemma \ref{onesimple} and Proposition \ref{Noumi}, all
spaces in {\bf (i)} and {\bf (ii)} are $\{0\}$
when $\lambda\in P_{2n}^+\setminus
(\Lambda_n^+)^\natural$.

Let $\lambda\in (\Lambda_n^+)^\natural$. Lemma
\ref{linksubalgebras} implies
$L_\lambda^{\AAA_\sigma}\supseteq
L_\lambda^{A(\mathfrak{k}_\sigma)}$,
$L_\lambda^{\omega(\AAA_\sigma)}\supseteq
L_\lambda^{A(\mathfrak{k}_\sigma)}$ and
$L_\lambda^{\omega_\delta(\AAA_\sigma)}\supseteq
L_\lambda^{A(\mathfrak{k}_\sigma^*)}$.
Since $\hbox{Dim}(L_\lambda^A)=1$ for
$A=A(\mathfrak{k}_\sigma)$ and $A=A(\mathfrak{k}_\sigma^*)$
by Proposition \ref{Noumi}, we conclude from Lemma \ref{onesimple}
that $L_\lambda^{\AAA_\sigma}=
L_\lambda^{A(\mathfrak{k}_\sigma)}$, $L_\lambda^{\omega(\AAA_\sigma)}=
L_\lambda^{A(\mathfrak{k}_\sigma)}$ and
$L_\lambda^{\omega_\delta(\AAA_\sigma)}=
L_\lambda^{A(\mathfrak{k}_\sigma^*)}$.
\end{proof}

The following result is a direct consequence of the definition
of $\omega_\delta$ and of Lemma \ref{zonal}.
\begin{cor}\label{spacesequiv}
$L_\lambda^{A(\mathfrak{k}_\sigma^*)}
=K^{-\delta}L_\lambda^{A(\mathfrak{k}_\sigma)}$ for all $\lambda\in P_{2n}^+$.
\end{cor}

The main object of study in \cite{NDS} and \cite{DS} is the
space of $(\mathfrak{k}_\sigma, \mathfrak{k}_\tau)$-fixed
regular functions on $\ug$,
which is defined as follows.
\begin{defi}
The space of $(\mathfrak{k}_\sigma, \mathfrak{k}_\tau)$-fixed
functions on $\ug$ is defined by
\[\mathcal{H}^{\sigma,\tau}=\{f\in \mathbb{C}_q[G] \,\, | \,\,
Y\cdot f=f\cdot Z=0\,\,\,\forall\,
Y\in\mathfrak{k}_\sigma,\,\forall\,Z\in\mathfrak{k}_\tau \}.
\]
\end{defi}

\begin{prop}\label{Hlink}
The subspaces $\mathcal{H}^{\sigma,\tau}$ and
$F_{\vec{0}}^{\sigma,\tau}$ of $\mathbb{C}_q[G]$ are related by
\[\mathcal{H}^{\sigma,\tau}=F_{\vec{0}}^{\sigma,\tau}\cdot
K^{-\delta}.
\]
\end{prop}
\begin{proof}
Any $f\in \mathcal{H}^{\sigma,\tau}$ can be written as
a linear combination of functions
\[\ug\ni X\mapsto
\langle Xv_\sigma,\widetilde{v}_\tau\rangle_\lambda
\]
with $v_\sigma\in L_\lambda^{A(\mathfrak{k}_\sigma)}$,
$\widetilde{v}_\tau\in L_\lambda^{A(\mathfrak{k}_\tau^*)}$
and $\lambda\in\bigl(\Lambda_n^+\bigr)^\natural$.
On the other hand, any $g\in F_{\vec{0}}^{\sigma,\tau}$ can be
written as a linear combination of functions
\[\ug\ni X\mapsto
\langle Xv_\sigma,\widehat{v}_\tau\rangle_\lambda
\]
with $v_\sigma\in L_\lambda^{\AAA_\sigma}$,
$\widehat{v}_\tau\in L_\lambda^{\AAA_\tau}$ and $\lambda\in
\bigl(\Lambda_n^+\bigr)^\natural$ since
$\rho(0,0)_\sigma=\rho(0)_\sigma=\epsilon|_{\AAA_\sigma}$.
The proposition follows
now from the fact that $L_{\lambda}^{A(\mathfrak{k}_\sigma)}=
L_{\lambda}^{\AAA_\sigma}$,
$L_\lambda^{A(\mathfrak{k}_\tau^*)}=
K^{-\delta}L_\lambda^{A(\mathfrak{k}_\tau)}=
K^{-\delta}L_\lambda^{\AAA_\tau}$ and from the fact
that $K^{-\delta}$ is
$*$-selfadjoint.
\end{proof}

\begin{rema}\label{alt}
{\bf (i)}
The fact that $\mathfrak{k}_\sigma$ and $\mathfrak{k}_\tau$ are
$\omega$-stable two-sided coideals implies that
$\mathcal{H}^{\sigma,\tau}\subseteq \mathbb{C}_q[G]$ is a unital
$*$-subalgebra, cf. \cite{NDS} and \cite{DS}. Furthermore,
$K^{-\delta}\in \ug$ is a group-like element, so the previous
proposition shows that $F_{\vec{0}}^{\sigma,\tau}\subseteq
\mathbb{C}_q[G]$ is a unital subalgebra (but it is not $*$-stable!).

{\bf (ii)} Observe that $\epsilon(\omega_\delta(X))=\overline{\epsilon(X)}$
for $X\in \ug$. Combined with Lemma \ref{zonal},
the proof of Proposition \ref{Hlink} and \eqref{starcommute},
this leads to the alternative description
\[\mathcal{H}^{\sigma,\tau}=\{f\in \mathbb{C}_q[G]
 \, | \, f(\omega_\delta(b)Xa)=
\epsilon(a)\overline{\epsilon(b)}f(X),\qquad \forall\, X\in \ug,\,\,
\forall\,a\in \AAA_\sigma,\,\,\forall\, b\in \AAA_\tau \}
\]
for the space of $(\mathfrak{k}_\sigma,\mathfrak{k}_\tau)$-fixed
regular functions on $\ug$.
\end{rema}

Recall that the restriction map $|_T:F_{\vec{0}}^{\sigma,\tau}\rightarrow
\mathbb{C}[u^{\pm 1}]$ involved a $\delta$-shift:
$f|_T(q^{\lambda})=f(K^{\lambda-\delta})$ for all $\lambda\in P_{2n}$.
When dealing with the $*$-subalgebra $\mathcal{H}^{\sigma,\tau}$,
it is convenient to use the restriction map without $\delta$-shift,
\[ \hbox{Res}_T: \mathcal{H}^{\sigma,\tau}\rightarrow \mathbb{C}[u^{\pm
  1}], \qquad \hbox{Res}_T(f)(q^{\lambda})=f(K^{\lambda})\quad\,\,
\forall\,\lambda\in P_{2n},
\]
since then we have
\[ \hbox{Res}_T(f\cdot K^{-\delta})=f|_T,\qquad \forall\,
f\in F_{\vec{0}}^{\sigma,\tau},
\]
cf. Proposition \ref{Hlink}.
Note that $\hbox{Res}_T$ is a $*$-algebra homomorphism, with
the $*$-structure on $\mathbb{C}[u^{\pm 1}]$ defined by
\begin{equation}\label{starpol}
p^*=\sum_{\mu\in\Lambda_n}\overline{c_{\mu}}u^{-\mu}\quad
\hbox{ if }\quad p=\sum_{\mu\in\Lambda_n}
c_{\mu}u^{\mu}\in\mathbb{C}[u^{\pm 1}].
\end{equation}
We are now in a position to translate the main result from \cite{NDS}
to our present setup.
\begin{prop}\label{validity}
Theorem \ref{maintheorem} is valid for $\vec{\kappa}=\vec{0}$.
\end{prop}
\begin{proof}
Observe that $\delta(0,0)=0$, hence we can take the unit $1\in\mathbb{C}_q[G]$
as nonzero element in $F_{\vec{0}}^{\sigma,\tau}(0)$.
Clearly $1|_T$ is the unit of $\mathbb{C}[u^{\pm 1}]$. This gives
Theorem \ref{maintheorem}{\bf (ii)} for $\vec{\kappa}=\vec{0}$.

The generalized Chevalley restriction theorem (Theorem
\ref{maintheorem}{\bf (i)}) for $\vec{\kappa}=\vec{0}$
reduces to the statement that $|_T: F_{\vec{0}}^{\sigma,\tau}\rightarrow
\mathbb{C}[u^{\pm 1}]^W$ is a linear isomorphism. In view of
Proposition \ref{Hlink} and the $\delta$-shift in the definition
of the restriction map $|_T$, this is equivalent to the statement that
$\hbox{Res}_T: \mathcal{H}^{\sigma,\tau}\rightarrow \mathbb{C}[u^{\pm
  1}]^W$ is a linear isomorphism, which follows from \cite[Thm. 3.2]{NDS}, see
also \cite[Thm. 7.5]{DS}. For a direct proof of this result for
an arbitrary quantum symmetric pair, see \cite{L4}.

Let $0\not=f_\mu\in
F_{\vec{0}}^{\sigma,\tau}(\mu)$ be an elementary
spherical function of degree $\mu\in\Lambda_n^+$, then
$f_\mu\cdot K^{-\delta}\in \mathcal{H}^{\sigma,\tau}\cap W(\mu^\natural)$.
Consequently, \cite[Thm. 3.4]{NDS} shows that the restriction
$f_\mu|_T=\hbox{Res}_T(f_\mu\cdot K^{-\delta})
\in \mathbb{C}[u^{\pm 1}]^W$ is a nonzero scalar multiple of
the Macdonald-Koornwinder polynomial
\begin{equation}\label{Pzonal}
P_{\mu}^{\sigma,\tau}(u)=
P_\mu(u;-q^{\sigma+\tau+1},-q^{-\sigma-\tau+1},q^{\sigma-\tau+1},
q^{-\sigma+\tau+1};q^2,q^{2}),
\end{equation}
(see also \cite[Thm. 7.5]{DS}).
This proves Theorem \ref{maintheorem}{\bf (iii)}
for $\vec{\kappa}=\vec{0}$.
\end{proof}


\section{The generalized Chevalley theorem}

In this section we prove the generalized Chevalley restriction
theorem (see Theorem \ref{maintheorem}{\bf (i)}).

To simplify notations, we write $\lbrack
\cdot,\cdot\rbrack_\mu$ for the
scalar product $\langle \cdot,
\cdot\rangle_{(\mu+\delta(\kappa,\kappa_1))^\natural}$
on the $\ug$-module
$L_{(\mu+\delta(\kappa,\kappa_1))^\natural}$ ($\mu\in \Lambda_n^+$).
\begin{lem}\label{groundstateform}
Let $\mu\in\Lambda_n^+$ and identify $V(\kappa,\kappa_1)_\tau$ with
its unique copy in the semisimple $\AAA_\tau$-module
$L_{(\mu+\delta(\kappa,\kappa_1))^\natural}^\tau$
and $V(\kappa_2)_\sigma$ with its unique copy in the semisimple
$\AAA_\sigma$-module
$L_{(\mu+\delta(\kappa,\kappa_1))^\natural}^\sigma$.
Let $f_\mu\in
F_{\vec{\kappa}}^{\sigma,\tau}(\mu)$ be an elementary vector valued
spherical function of degree $\mu$. Then
\[f_\mu|_T(q^\lambda)=\lbrack
K^{\lambda-\delta}v,w\rbrack_{\mu},
\qquad \forall\,\lambda\in P_{2n}
\]
for suitably normalized vectors $0\not=v\in V(\kappa_2)_\sigma$
and $0\not=w\in\widetilde{V}(\kappa,\kappa_1)_\tau$. Furthermore,
\[f_\mu|_T=\sum_{
-\mu-\delta(\kappa,\kappa_1)\leq\nu\leq\mu+\delta(\kappa,\kappa_1)}
c_\nu u^{\nu}\in \mathbb{C}[u^{\pm 1}]
\]
for certain constants $c_\nu\in\mathbb{C}$
\textup{(}$\nu\in\Lambda_n$\textup{)}, and
$c_{\pm\mu\pm\delta(\kappa,\kappa_1)}\not=0$.
\end{lem}
\begin{proof}
Choose an orthonormal basis $\{b_1,\ldots,b_m\}$ of
$V(\kappa,\kappa_1)_\tau$ with respect to the restriction of
$\lbrack \cdot,\cdot\rbrack_\mu$ to $V(\kappa,\kappa_1)_\tau$
such that $b_1\in \widetilde{V}(\kappa,\kappa_1)_\tau$.
In particular, $b_1$ is a nonzero
constant multiple of $w$. Then
$f_\mu\in F_{\vec{\kappa}}^{\sigma,\tau}(\mu)$
can be realized as
\[f_\mu(X)=\sum_{i=1}^m\lbrack Xv, b_i\rbrack_\mu\, b_i,\qquad X\in\ug
\]
for some $0\not=v\in V(\kappa_2)_\sigma$.
We have
$f_\mu(X)\in\widetilde{V}(\kappa,\kappa_1)_\tau$
for $X\in U^0$, see (the proof of) Lemma \ref{restriction}.
Using the identification
$\widetilde{V}(\kappa,\kappa_1)_\tau\simeq\mathbb{C}$ of vectorspaces
given by the assignment $b_1\mapsto 1$, we obtain
\[f_\mu|_T(q^\lambda)=f_\mu(K^{\lambda-\delta})=
\lbrack K^{\lambda-\delta}v,b_1\rbrack_\mu,\qquad \forall\,
\lambda\in P_{2n}.
\]
Now both vectors $v$ and $b_1$ satisfy $C_iv=v$ and $C_ib_1=b_1$ for
$i=1,\ldots,n$, hence $v[\lambda]=0$ and $b_1[\lambda]=0$ for
$\lambda\in P_{2n}$ unless $\lambda\in \bigl(\Lambda_n\bigr)^\natural$.
Since $\natural: \Lambda_n\rightarrow \bigl(\Lambda_n\bigr)^\natural\subset
P_{2n}$ is a bijection of partially ordered sets and
$(-\nu)^\natural=w_0(\nu^\natural)$
for $\nu\in\Lambda_n$, we conclude that
the weight decompositions of $v$ and $b_1$ are of the form
\[v=\sum_{-\mu-\delta(\kappa,\kappa_1)\leq
\nu\leq\mu+\delta(\kappa,\kappa_1)}v[\nu^\natural],
\qquad\qquad
b_1=\sum_{-\mu-\delta(\kappa,\kappa_1)\leq
\nu\leq\mu+\delta(\kappa,\kappa_1)}b_1[\nu^\natural],
\]
with $v[\nu^\natural]$ and $b_1[\nu^\natural]$ the components of the
vectors $v$ and $b_1$ in the weight space of the $\ug$-module
$L_{(\mu+\delta(\kappa,\kappa_1))^\natural}$ of weight
$\nu^\natural$.

A simple application of Lemma \ref{highcontri} shows that
$v[(\pm\mu\pm\delta(\kappa,\kappa_1))^\natural]\not=0$. Denote
$\widetilde{L}^\tau_{(\mu+\delta(\kappa,\kappa_1))^\natural}$ for
the vectors $v\in L^\tau_{(\mu+\delta(\kappa,\kappa_1))^\natural}$
satisfying $C_iv=v$ for $i=1,\ldots,n$. This subspace contains the
one dimensional weight spaces
$L_{(\mu+\delta(\kappa,\kappa_1))^\natural}
[(\pm\mu\pm\delta(\kappa,\kappa_1))^\natural]$ and
$\widetilde{L}^\tau_{(\mu+\delta(\kappa,\kappa_1))^\natural}\cap
\widetilde{V}(\kappa,\kappa_1)_\tau=\mathbb{C}\{b_1\}$. Lemma
\ref{highcontri} implies
$b_1[(\pm\mu\pm\delta(\kappa,\kappa_1))^\natural]\not=0$. Hence
\begin{equation*}
\begin{split}
f_\mu|_T&=
\sum_{-\mu-\delta(\kappa,\kappa_1)\leq\nu\leq\mu+\delta(\kappa,\kappa_1)}
\lbrack K^{-\delta}v[\nu^\natural],
b_1[\nu^\natural]\rbrack_\mu\,u^{\nu}\\
&=\sum_{-\mu-\delta(\kappa,\kappa_1)\leq
\nu\leq\mu+\delta(\kappa,\kappa_1)}c_\nu
u^\nu
\end{split}
\end{equation*}
with $c_\nu\in\mathbb{C}$ satisfying the conditions as stated in the lemma.
\end{proof}

A more detailed analysis of the elementary vector valued spherical functions
requires the right coideal algebra structure of $\AAA_\sigma$.
Recall that a right $\ug$-comodule $(M,\delta_M)$ is a
vectorspace $M$ together with a linear map
$\delta_M: M\rightarrow M\otimes \ug$ satisfying
\[\bigl(\delta_M\otimes {\hbox{Id}}_{\ug}\bigr)\circ\delta_M=
\bigl({\hbox{Id}}_M\otimes \Delta\bigr)\circ\delta_M,\qquad
\bigl({\hbox{Id}}_M\otimes\epsilon\bigr)\circ\delta_M={\hbox{Id}}_M.
\]
If $M$ is a unital algebra and $\delta_M$ is
a (unit preserving)
algebra homomorphism, then $(M,\delta_M)$ is called a right
$\ug$-comodule algebra. In particular, a unital subalgebra
$A\subseteq\ug$ satisfying $\Delta(A)\subseteq A\otimes \ug$
gives rise to the right comodule algebra $(A,\Delta|_A)$.
In this case, $A$ is called a right coideal subalgebra of $\ug$.

By computing the action of the comultiplication $\Delta$ on the
algebraic generators of $\AAA_\sigma$, one verifies that
$\AAA_\sigma\subseteq\ug$ is a right $\ug$-coideal algebra.
Using the isomorphism $\pi_\sigma: \AAA\rightarrow \AAA_\sigma$ of
algebras, see Proposition \ref{embeddingsigma}, $\AAA$ inherets a
unique right $\ug$-comodule structure such that $\pi_\sigma:
\AAA\rightarrow \AAA_\sigma$ is an isomorphism of right
$\ug$-comodule algebras. The right comodule map for $\AAA$
turns out to be independent of $\sigma$, and is given as follows,
cf., e.g., \cite{L3}.

\begin{prop}\label{coalgebra}
The algebra $\AAA$ has the structure of a right
$\ug$-comodule algebra, with comodule map $\delta_\AAA$ given
explicitly by
\begin{equation*}
\begin{split}
\delta_\AAA(\gamma_i^{\pm 1})&=\gamma_i^{\pm 1}\otimes
K_i^{\pm 1}K_{2n+1-i}^{\pm 1},\\
\delta_\AAA(\beta_j)&=\beta_j\otimes K_j^{-1}K_{2n-j}^{-1}+
\gamma_{j+1}^{-1}\otimes y_jK_{j+1}^{-1}K_{2n-j}^{-1}+
\gamma_j^{-1}\otimes K_j^{-1}x_{2n-j}K_{2n-j}^{-1}
\end{split}
\end{equation*}
for $i=1,\ldots,2n$ and $j=1,\ldots,2n-1$. The map
$\pi_\sigma: \AAA\rightarrow\AAA_\sigma$ is an isomorphism of right
$\ug$-comodule algebras.
\end{prop}
Finite dimensional $\AAA_\sigma$-modules do not form
a tensor category, since $\AAA_\sigma\subset \ug$ is not a coalgebra.
The right coideal algebra structure of $\AAA_\sigma$ implies
though that the tensor product $M\otimes N$ of an $\AAA_\sigma$-module
$M$ and an $\ug$-module $N$ is an $\AAA_\sigma$-module by the usual formula
\[a\cdot(m\otimes n)=\sum a_1\cdot m\otimes a_2\cdot n,\qquad
m\in M,\,\, n\in N,
\]
where (recall) $\Delta(a)=\sum a_1\otimes a_2$ for $a\in\AAA_\sigma$
with the $a_1$'s
from $\AAA_\sigma$. In the special case that $M$ is the one-dimensional
$\AAA_\sigma$-module $V(k)_\sigma$ ($k\in\mathbb{Z}$), we can naturally
relate the $\AAA_\sigma$-module $V(k)_\sigma\otimes N$ to the
$\AAA_{\sigma+2k}$-module $N^{\sigma+2k}$ in the following manner.

Let $\chi_{k}^\sigma: \AAA_\sigma\rightarrow \CC$ be the
character of $V(k)_\sigma$.
The $\AAA_\sigma$-action on $V(k)_\sigma\otimes N$
can then be rewritten as
\[a\cdot(m\otimes n)=m\otimes \zeta_k^\sigma(a)\cdot n,\qquad
a\in\AAA_\sigma,\,\,m\in V(k)_\sigma,\,\, n\in N,
\]
with $\zeta_k^\sigma: \AAA_\sigma\rightarrow \ug$ the
unital algebra homomorphism
\begin{equation}\label{pushforward}
\zeta_{k}^\sigma(a)=\sum \chi_{k}^\sigma(a_1)a_2,\qquad a\in \AAA_\sigma.
\end{equation}
\begin{lem}\label{pushforwardlem}
The map $\zeta_{k}^\sigma$ defines an algebra isomorphism
$\zeta_{k}^\sigma: \AAA_\sigma\rightarrow \AAA_{\sigma+2k}$,
with inverse $\zeta_{-k}^{\sigma+2k}:
\AAA_{\sigma+2k}\rightarrow \AAA_\sigma$.
Furthermore,
$\zeta_k^\sigma=\pi_{\sigma+2k}\circ \pi_\sigma^{-1}|_{\AAA_\sigma}$
\textup{(}see Proposition \ref{embeddingsigma} for the definition of
$\pi_\sigma$\textup{)}, which means that
\begin{equation}\label{zetaactie}
\zeta_{k}^\sigma(C_i)=C_i,\qquad
\zeta_{k}^\sigma(B_j)=B_j,\qquad
\zeta_{k}^\sigma(B_n^\sigma)=B_n^{\sigma+2k}
\end{equation}
for $i\in\{1,\ldots,2n\}$ and $j\in\{1,\ldots,2n-1\}\setminus\{n\}$.
\end{lem}
\begin{proof}
By Lemma \ref{explicit} the values of $\chi_{k}^\sigma$
on the algebraic generators of $\AAA_\sigma$ are given by
\begin{equation}\label{chiactie}
\chi_{k}^\sigma(C_i)=1,\qquad
\chi_{k}^\sigma(B_j)=0,\qquad
\chi_{k}^\sigma(B_n^\sigma)=\vartheta_0(q^{\sigma+2k})
\end{equation}
for $i\in\{1,\ldots,2n\}$ and $j\in\{1,\ldots,2n-1\}\setminus \{n\}$,
with $\vartheta_l(s)$ ($l\in\mathbb{Z}$) defined by
\[
\vartheta_l(s)=\frac{s^{-1}-sq^{-2l}}{q-q^{-1}}.
\]
Combined with Proposition \ref{coalgebra} one now
easily derives \eqref{zetaactie}.
The remaining statements follow now immediately.
\end{proof}

\begin{cor} \label{pushforwardcor}
Let $N$ be a finite dimensional $\ug$-module and suppose that
$N^\prime\subset N^{\sigma+2k}$ is an $\AAA_{\sigma+2k}$-submodule.

{\bf (i)} $V(k)_\sigma\otimes N^\prime$ is an $\AAA_\sigma$-submodule
of $V(k)_\sigma\otimes N$.

{\bf (ii)} If $N^\prime\simeq V(\kappa,\kappa_1)_{\sigma+2k}$, then
$V(k)_\sigma\otimes N^\prime\simeq V(\kappa,\kappa_1+k)_\sigma$ as
$\AAA_\sigma$-modules.

\end{cor}
\begin{proof}
{\bf (i)} follows from Lemma \ref{pushforwardlem}, and
{\bf (ii)} follows from Lemma \ref{pushforwardlem} and
Lemma \ref{explicit}.
\end{proof}
As a first application of these tensor product constructions
we derive (a refinement of) the generalized Chevalley restriction theorem
(see Theorem \ref{maintheorem}{\bf (i)}). Before stating the result,
we first recall some basic facts on tensor products of
finite dimensional $\ug$-modules which we need in the proof.
For any $\lambda,\mu\in P_{2n}^+$, the irreducible decomposition of
the finite dimensional
$\ug$-module $L_\lambda\otimes L_\mu$ is of the form
\begin{equation}\label{triang}
L_\lambda\otimes L_\mu= L_{\lambda+\mu}
\oplus \bigoplus_{\stackrel{\scriptstyle{\nu\in P_{2n}^+}}
{\nu\prec \lambda+\mu}}
L_\nu^{\oplus d_{\nu}^{\lambda,\mu}}
\end{equation}
for certain multiplicities $d_{\nu}^{\lambda,\mu}\in\mathbb{Z}_{\geq
  0}$. The copy $L_{\lambda+\mu}$ has $v_\lambda\otimes v_\mu$
as highest weight vector. We write $\hbox{pr}_{\lambda,\mu}:
L_\lambda\otimes L_\mu\rightarrow L_{\lambda+\mu}$ for the
projection along the direct sum decomposition \eqref{triang},
and $\hbox{pr}_\mu=\hbox{pr}_{\mu^\natural,\delta(\kappa,\kappa_1)^\natural}$
for $\mu\in\Lambda_n^+$.

\begin{prop}\label{triangular}
Let $\mu\in \Lambda_n^+$ and choose elementary vector valued spherical
functions
$f_\mu\in F_{\vec{\kappa}}^{\sigma,\tau}(\mu)$ and
$f_0\in F_{\vec{\kappa}}^{\sigma,\tau}(0)$. Then
$f_\mu|_T$ is divisible by $f_0|_T$ in $\mathbb{C}[u^{\pm 1}]$, and
\begin{equation}\label{triangularform}
\frac{f_\mu|_T}{f_0|_T}=
\sum_{\stackrel{\scriptstyle{\nu\in\Lambda_n^+}}{\nu\leq\mu}}
c_{\nu}m_{\nu}\in\mathbb{C}[u^{\pm 1}]^W
\end{equation}
for some constants $c_\nu\in \mathbb{C}$ with $c_\mu\not=0$.
\end{prop}
\begin{proof}
The proof is by induction on $\mu\in\Lambda_n^+$ along the
dominance order $\leq$. For $\mu=0$ the proposition is trivial.
We identify $V(0)_\sigma$ with its unique copy in
$L_{\mu^\natural}^\sigma$ and
$V(\kappa_2)_\sigma$ with its unique copy in
$L_{\delta(\kappa,\kappa_1)^\natural}^\sigma$, and
we choose nonzero vectors
\[
w_\sigma\in V(0)_\sigma\subset L_{\mu^\natural}^\sigma,\qquad
u_\sigma\in V(\kappa_2)_\sigma\subset
L_{\delta(\kappa,\kappa_1)^\natural}^\sigma.
\]
We furthermore fix nonzero intertwiners
\[\phi_\mu\in\hbox{Hom}_{\AAA_\tau}(L_{\mu^\natural}^\tau,
V(0)_\tau),\qquad
\phi_0\in \hbox{Hom}_{\AAA_\tau}(L_{\delta(\kappa,\kappa_1)^\natural}^\tau,
V(\kappa,\kappa_1)_\tau)
\]
and we identify $V(0)_\tau\simeq \mathbb{C}$ as vector spaces.
Then
\[
g_\mu(X)=\phi_\mu(Xw_\sigma),\qquad f_0(X)=\phi_0(Xu_\sigma)
\]
for $X\in \ug$ define elementary vector valued spherical functions
$g_\mu\in F_{\vec{0}}^{\sigma,\tau}(\mu)$ and
$f_0\in F_{\vec{\kappa}}^{\sigma,\tau}(0)$. Furthermore,
$g_\mu|_T$ is of the form
\begin{equation}\label{gexpansion}
g_\mu|_T=\sum_{\stackrel{\nu\in\Lambda_n^+}{\nu\leq\mu}}
d_\nu m_{\nu}
\end{equation}
for some constants $d_\nu\in\mathbb{C}$ with $d_\mu\not=0$ in view of
Proposition \ref{validity} (in fact, $g_\mu|_T$ is a nonzero constant
multiple of a Macdonald-Koornwinder polynomial of degree $\mu$).

We now consider the linear map $f: \ug\rightarrow V(\kappa,\kappa_1)_\tau$
defined by
\[f(X)=\phi\bigl(X(w_\sigma\otimes u_\sigma)\bigr),\qquad X\in \ug,
\]
with $\phi: L_{\mu^\natural}\otimes L_{\delta(\kappa,\kappa_1)^\natural}
\rightarrow V(\kappa,\kappa_1)_\tau$ the linear map defined by
\[\phi(u\otimes v)=\phi_\mu(u)\phi_0(v),\qquad u\in L_{\mu^\natural},\,\,
v\in L_{\delta(\kappa,\kappa_1)^\natural}.
\]
Observe that
$\mathbb{C}\{w_\sigma\otimes u_\sigma\}$ is an $\AAA_\sigma$-submodule
of $\bigl(L_{\mu^\natural}\otimes
L_{\delta(\kappa,\kappa_1)^\natural}\bigr)^\sigma$ which is isomorphic to
$V(\kappa_2)_\sigma$ by Corollary \ref{pushforwardcor}. Furthermore,
\[
\phi\in \hbox{Hom}_{\AAA_\tau}\bigl(
(L_{\mu^\natural}\otimes L_{\delta(\kappa,\kappa_1)^\natural})^\tau,
V(\kappa,\kappa_1)_\tau),
\]
since
\begin{equation*}
\begin{split}
\phi(a(v\otimes w))&=\sum\phi_\mu(a_1v)\phi_0(a_2w)\\
&=\phi_\mu(v)\phi_0(\zeta_0^{\tau}(a)w)\\
&=\rho(\kappa,\kappa_1)_\tau(a)\phi(v\otimes w)
\end{split}
\end{equation*}
for
$a\in\AAA_\tau$, $v\in L_{\mu^\natural}$ and
$w\in L_{\delta(\kappa,\kappa_1)^\natural}$,
where, besides the intertwining properties of $\phi_\mu$ and $\phi_0$,
we have used that $\zeta_0^\tau=\hbox{Id}_{\AAA_\tau}$. Hence we conclude
that $f\in F_{\vec{\kappa}}^{\sigma,\tau}$.

By \eqref{triang}
and Corollary \ref{Cormis} we have an expansion
\begin{equation}\label{fexpansion}
f=\sum_{\stackrel{\nu^\prime\in\Lambda_n^+}{\nu^\prime\leq\mu}}
f_{\nu^\prime}
\end{equation}
with $f_{\nu^\prime}\in F_{\vec{\kappa}}^{\sigma,\tau}(\nu^\prime)$
and $f_\mu(\cdot)=\phi(\cdot\,\hbox{pr}_\mu(w_\sigma\otimes u_\sigma))$.
Lemma \ref{highcontri} implies that
$\hbox{pr}_\mu(w_\sigma\otimes u_\sigma)$ spans the unique copy of
$V(\kappa_2)_\sigma$ in $L_{(\mu+\delta(\kappa,\kappa_1))^\natural}^\sigma$.
Furthermore, $\phi\circ\hbox{pr}_{\mu}\not=0$ since
\begin{equation}\label{extra}
\phi(v_{\mu^\natural}\otimes v_{\delta(\kappa,\kappa_1)^\natural})=
\phi_\mu(v_{\mu^\natural})\phi_0(v_{\delta(\kappa,\kappa_1)^\natural})\not=0
\end{equation}
by Proposition \ref{expectationvalue}. Consequently, $f_\mu$ is nonzero.

Up to a nonzero multiplicative constant, $f_\mu$ is the
unique elementary vector valued spherical function of degree $\mu$.
It thus suffices to prove the induction step for $f_\mu$.
By the construction
of the vector valued spherical function $f$, we have
\begin{equation}\label{prodtorus}
f|_T=(g_\mu|_T)(f_0|_T)
\end{equation}
in $\mathbb{C}[u^{\pm 1}]$.
Substituting the expansions \eqref{gexpansion} and \eqref{fexpansion}
of $g_\mu|_T$ and $f$ respectively in \eqref{prodtorus}, gives
\[\sum_{
\stackrel{\scriptstyle{\nu\in\Lambda_n^+}}{\nu\leq\mu}}
d_\nu\,f_0|_Tm_{\nu}=f|_T
=f_\mu|_T+
\sum_{\stackrel{\scriptstyle{\nu^\prime\in\Lambda_n^+}}
{\nu^\prime<\mu}}
f_{\nu^\prime}|_T.
\]
By the induction hypothesis applied to the $f_{\nu^\prime}|_T$'s
($\nu^\prime<\mu$), we conclude that \eqref{triangularform}
is valid for some $c_\nu\in\mathbb{C}$ and that $c_\mu=d_\mu\not=0$. This
completes the proof of the induction step.

\end{proof}

As an immediate corollary of Proposition \ref{triangular},
we obtain the generalized Chevalley restriction theorem
(see Theorem \ref{maintheorem}{\bf (i)}):
\begin{cor}
Let $f_0\in F_{\vec{\kappa}}^{\sigma,\tau}(0)$ be an elementary vector
valued spherical function of degree zero.
The restriction map $|_T$ induces a linear bijection
$|_T: F_{\vec{\kappa}}^{\sigma,\tau}\rightarrow
f_0|_T\,\mathbb{C}[u^{\pm 1}]^W$.
\end{cor}


\section{The ground state}
We define a ground state $f_0$ (with respect to $(\vec{\kappa},\sigma,\tau)$)
to be an elementary vector valued spherical function
$f_0\in F_{\vec{\kappa}}^{\sigma,\tau}(0)$ of degree zero.
Note that a ground state $f_0$ is unique up to nonzero scalar multiples.
In this section we compute the radial part $f_0|_T$ of
a ground state $f_0$ explicitly. In the first subsection we show
that it suffices to compute the
ground state for $\vec{\kappa}=(0,0,\kappa)$ and
$\vec{\kappa}=(\kappa_1,\kappa_2,0)$ seperately.
The second and third subsections
are devoted to the explicit computation of the ground state
for these two special cases.


\subsection{Splitting of the ground state.}

For the computation of the radial part of the ground state it is
sufficient to compute the radial part for the special cases
$\kappa_1=\kappa_2=0$ and $\kappa=0$ in view of the following lemma.

\begin{lem}\label{splitting}
Let $f\in F_{(\kappa_1,\kappa_2,0)}^{\sigma,\tau}(0)$ and
$f^\prime\in
F_{(0,0,\kappa)}^{\sigma+2\kappa_2,\tau+2\kappa_1}(0)$ be ground
states with respect to $(\kappa_1,\kappa_2,0,\sigma,\tau)$ and
$(0,0,\kappa,\sigma+2\kappa_2,\tau+2\kappa_1)$, respectively.
Then $f|_Tf^\prime|_T=f_0|_T$
in $\mathbb{C}[u^{\pm 1}]$ with $f_0\in F_{\vec{\kappa}}^{\sigma,\tau}(0)$
a ground state with respect to $(\vec{\kappa},\sigma,\tau)$.
\end{lem}
\begin{proof}
We proceed as in the first part of the proof of Proposition
\ref{triangular}.
We identify $V(\kappa_2)_\sigma$ with its unique copy in
$L_{\delta(0,\kappa_1)^\natural}^\sigma$ and $V(0)_{\sigma+2\kappa_2}$
with its unique copy in $L_{\delta(\kappa,0)^\natural}^{\sigma+2\kappa_2}$,
and we fix nonzero vectors $u\in V(\kappa_2)_\sigma$ and
$v\in V(0)_{\sigma+2\kappa_2}$. We furthermore fix nonzero intertwiners
\[\varphi\in \hbox{Hom}_{\AAA_\tau}
\bigl(L_{\delta(0,\kappa_1)^\natural}^\tau, V(\kappa_1)_\tau\bigr),
\qquad
\phi\in\hbox{Hom}_{\AAA_{\tau+2\kappa_1}}
\bigl(L_{\delta(\kappa,0)^\natural}^{\tau+2\kappa_1},
V(\kappa,0)_{\tau+2\kappa_1}\bigr),
\]
then the formulas
\[f(X)=\varphi(Xu),\qquad f^\prime(X)=\phi(Xv)
\]
for $X\in \ug$ define ground states
$f\in F_{(\kappa_1,\kappa_2,0)}^{\sigma,\tau}(0)$ and
$f^\prime\in F_{(0,0,\kappa)}^{\sigma+2\kappa_2,\tau+2\kappa_1}(0)$.
We now consider the function $f_0(X)=\psi(X(u\otimes v))$ for $X\in \ug$, with
\begin{equation}\label{psiintertwiner}
\psi:L_{\delta(0,\kappa_1)^\natural}\otimes L_{\delta(\kappa,0)^\natural}
\rightarrow V(\kappa_1)_\tau\otimes V(\kappa,0)_{\tau+2\kappa_1}\simeq
V(\kappa,\kappa_1)_\tau
\end{equation}
defined by $\psi(w\otimes w^\prime)=\varphi(w)\otimes\phi(w^\prime)$ for
$w\in L_{\delta(0,\kappa_1)^\natural}$ and $w^\prime\in
L_{\delta(\kappa,0)^\natural}$ (see Corollary \ref{pushforwardcor}
for the isomorphism $V(\kappa_1)_\tau\otimes V(\kappa,0)_{\tau+2\kappa_1}\simeq
V(\kappa,\kappa_1)_\tau$ of $\AAA_\tau$-modules in \eqref{psiintertwiner}).
Similarly as in the proof of
Proposition \ref{triangular} we have
\[f_0|_T=(f|_T)(f^\prime|_T)
\]
in $\mathbb{C}[u^{\pm 1}]$ and $0\not=f_0\in F^{\sigma,\tau}_{\vec{\kappa}}$.
In the expansion of $f_0$ in elementary vector valued spherical functions,
there is only a contribution of degree zero in view of
\eqref{triang} and Corollary \ref{Cormis}, hence we conclude that
$f_0\in F_{(\kappa_1,\kappa_2,\kappa)}^{\sigma,\tau}(0)$
is the ground state with respect to $(\vec{\kappa},\sigma,\tau)$.
\end{proof}


\subsection{The ground state for $\kappa_1=\kappa_2=0$.}

We compute the radial part $f_0|_T\in\mathbb{C}[u^{\pm 1}]$ of the
ground state
$f_0\in F_{(0,0,\kappa)}^{\sigma,\tau}(0)$ by relating $f_0|_T$ to the
ground state for vector-valued
$U_q(\mathfrak{g}\mathfrak{l}(n))$-characters, which in turn was
computed explicitly by Etingof and Kirillov \cite{EK}.

We use the notations of \S 1.
Denote $(\rho_{\kappa}, V_\kappa)$ for the finite
dimensional, irreducible $U_q(\mathfrak{g}\mathfrak{l}(n))$-representation
$L_{(n-1)\kappa,(-\kappa)^{n-1}}$.
We view an element
$f\in \mathbb{C}_q[\hbox{GL}(n;\mathbb{C})]^{\otimes 2}\otimes V$ as a
linear map  $f: U_q(\mathfrak{g}\mathfrak{l}(n))^{\otimes 2}\rightarrow V$.
\begin{defi}
Let $\widetilde{F}_\kappa$ be the space of functions $f\in
\mathbb{C}_q[\hbox{GL}(n;\mathbb{C})]^{\otimes 2}\otimes V_\kappa$ satisfying
\begin{equation*}
\begin{split}
f(X\Delta^{op}(a))&=\epsilon(a)f(X),\qquad \forall\, a\in
U_q(\mathfrak{g}\mathfrak{l}(n)),\\
f(\Delta^{op}(b)X)&=\rho_\kappa(b)f(X),\qquad \forall\, b\in
U_q(\mathfrak{g}\mathfrak{l}(n))
\end{split}
\end{equation*}
for all $X\in U_q(\mathfrak{g}\mathfrak{l}(n))^{\otimes 2}$.
\end{defi}
Recall from Lemma \ref{restrictionmodule} that the
isomorphism $U_q(\mathfrak{g}\mathfrak{l}(n))\simeq\mathcal{E}$ of
Lemma \ref{diagonal} implies that the restriction of $V(\kappa,\kappa_1)_\tau$
to $\mathcal{E}$ is isomorphic to $V_\kappa$ as
$U_q(\mathfrak{g}\mathfrak{l}(n))$-modules. In the following lemma
we furthermore use the natural identification
$U_q(\mathfrak{g}\mathfrak{l}(n))^{\otimes 2}\simeq U_q(\mathfrak{k})\subset
\ug$, cf. \S 1.3.
\begin{lem}
The assignment
\[f\mapsto f\circ(\hbox{Id}\otimes \psi),
\]
with $\psi$ defined by \eqref{psifunction},
defines a linear map $F_{\vec{\kappa}}^{\sigma,\tau}\rightarrow
\widetilde{F}_\kappa$.
\end{lem}
\begin{proof}
With the identifications discussed above, the lemma follows
immediately from  the explicit form of the isomorphism
$U_q(\mathfrak{g}\mathfrak{l}(n))\simeq\mathcal{E}$,
see Lemma \ref{diagonal}.
\end{proof}
For a left $U_q(\mathfrak{g}\mathfrak{l}(n))$-module $V$,
we consider its linear dual $V^*$ as $U_q(\mathfrak{g}\mathfrak{l}(n))$-module
by
\[
(X\phi)(v)=\phi(S^{-1}(X)v),\qquad
X\in U_q(\mathfrak{g}\mathfrak{l}(n)),\, \phi\in V^*,\,v\in V.
\]
Let $(V\otimes V^*)^{op}$ be the $U_q(\mathfrak{g}\mathfrak{l}(n))$-module
with representation space $V\otimes V^*$ and
$U_q(\mathfrak{g}\mathfrak{l}(n))$-action
\[X(v\otimes \phi)=\sum X_2v\otimes X_1\phi,
\qquad X\in U_q(\mathfrak{g}\mathfrak{l}(n)),\,\,v\in V,\,\,
\phi\in V^*.
\]
\begin{lem}
Let $V$ be a finite dimensional
$U_q(\mathfrak{g}\mathfrak{l}(n))$-module. The
$U_q(\mathfrak{g}\mathfrak{l}(n))$-in\-ter\-twi\-ners
$\xi: (V\otimes V^*)^{op}\rightarrow V_\kappa$ are in bijective
correspondence with $U_q(\mathfrak{g}\mathfrak{l}(n))$-intertwiners
$\xi^*: V\rightarrow V\otimes V_\kappa$, where $V\otimes V_\kappa$
is now considered as $U_q(\mathfrak{g}\mathfrak{l}(n))$-module by
\[X(v\otimes w)=\sum X_1v\otimes X_2w,\qquad X\in
U_q(\mathfrak{g}\mathfrak{l}(n)),\, v\in V,\, w\in V_\kappa.
\]
The bijective correspondence is explicitly given by the formula
\[(\phi\otimes \hbox{Id})\xi^*(v)=\xi(v\otimes
K^{2\delta_n}\phi),\qquad \phi\in V^*,\,v\in V
\]
where \textup{(}recall\textup{)}
$K^{2\delta_n}=K_1^{2(n-1)}K_2^{2(n-2)}\cdots K_{n-1}^2\in
U_q(\mathfrak{g}\mathfrak{l}(n))$.
\end{lem}
\begin{proof}
The standard proof is left to the reader.
\end{proof}
Let $V$ be a finite dimensional $U_q(\mathfrak{g}\mathfrak{l}(n))$-module
with basis $\{v_i\}$ and corresponding dual basis $\{v_i^*\}$.
Let $\xi: (V\otimes V^*)^{op}\rightarrow V_\kappa$ be an
$U_q(\mathfrak{g}\mathfrak{l}(n))$-intertwiner. Then
\[g_{V,\xi}(\cdot)=\sum_i\xi(\cdot\,(v_i\otimes v_i^*))
\]
defines an element in $\widetilde{F}_{\kappa}$.
Any element in $\widetilde{F}_{\kappa}$ is of this form for some $V$ and
$\xi$.

\begin{lem}\label{groundstate1}
The radial part $f|_T\in\mathbb{C}[u^{\pm 1}]$ of a vector valued spherical
function $f\in F_{\vec{\kappa}}^{\sigma,\tau}$ is divisible by
\begin{equation*}
I=u^{\delta(\kappa,0)}\prod_{1\leq i<j\leq n}
\bigl(q^2u_i^{-1}u_j;q^2\bigr)_{\kappa}
\end{equation*}
in $\mathbb{C}[u^{\pm 1}]$. Moreover, the quotient
$f|_T/I\in \mathbb{C}[u^{\pm 1}]$ is $S_n$-invariant.
\end{lem}
\begin{proof}
Let $f\in F_{\vec{\kappa}}^{\sigma,\tau}$ be a vector valued
spherical function and set $g=f\circ(\hbox{Id}\otimes
\psi)\in\widetilde{F}_\kappa$. Then $g=g_{V,\xi}$ for some finite
dimensional $U_q(\mathfrak{g}\mathfrak{l}(n))$-module $V$ and some
$U_q(\mathfrak{g}\mathfrak{l}(n))$-intertwiner $\xi: (V\otimes
V^*)^{op}\rightarrow V_{\kappa}$. Let $\lambda\in P_{2n}$, then
\[f|_T(q^\lambda)=f(K^{\lambda-\delta_{2n}})=
g(K^{\lambda^\ddagger-(1^n)}\otimes K^{2\delta_n})
\]
with
\begin{equation}\label{lambdaddagger}
\lambda^\ddagger=(\lambda_1-\lambda_{2n},\lambda_2-\lambda_{2n-1},\ldots,
\lambda_n-\lambda_{n+1})\in P_n.
\end{equation}
By the previous lemma we thus obtain
\begin{equation*}
\begin{split}
f|_T(q^\lambda)&=\sum_i\xi\bigl(K^{\lambda^\ddagger-(1^n)} v_i\otimes
K^{2\delta_n}v_i^*\bigr)\\
&=\sum_i(v_i^*\otimes \hbox{Id})\xi^*(K^{\lambda^\ddagger-(1^n)}v_i)\\
&=\hbox{Tr}|_{V}\bigl(\xi^*(K^{\lambda^\ddagger-(1^n)}(\,\cdot\,))\bigr).
\end{split}
\end{equation*}
The results of Etingof and Kirillov \cite{EK} imply
\[\hbox{Tr}|_V\bigl(\xi^*(K^{\mu}(\,\cdot\,))\bigr)=
I(q^{\mu})p(q^\mu),\qquad \forall\mu\in P_n
\]
with $q^\mu=(q^{\mu_1},q^{\mu_2},\ldots,q^{\mu_n})$ and
$p\in \mathbb{C}[u^{\pm 1}]^{S_n}$.
Since $I\in \mathbb{C}[u^{\pm 1}]$ is a homogeneous Laurent polynomial,
we arrive at the desired result.
\end{proof}

We continue our analysis of the ground state for $\kappa_1=\kappa_2=0$
following the method of Kirillov \cite{Kir}. We express the
elements $B_j$ ($j\not=n$) and $B_n^\sigma$ of $\AAA_\sigma$ by
\[B_j=\widetilde{y}_j+\widetilde{x}_{2n-j},\qquad
B_n^\sigma=\widetilde{y}_n+\widetilde{x}_n+\vartheta_0(q^\sigma)K_n^{-2}
\]
with
\[\widetilde{x}_j=K_{2n-j}^{-1}x_jK_j^{-1},\qquad
\widetilde{y}_j=y_jK_{j+1}^{-1}K_{2n-j}^{-1},
\]
and with $\vartheta_l$ given by
\begin{equation}\label{lambdak}
\vartheta_l(s)=\frac{s^{-1}-sq^{-2l}}{q-q^{-1}},\qquad l\in\mathbb{Z},
\end{equation}
cf. Lemma \ref{pushforwardlem}. We define elements
$E_{n-1},E_n,F_{n-1},F_n\in \ug$ by
\begin{equation*}
\begin{split}
E_{n-1}&=\widetilde{x}_{n+1},\qquad E_n=\widetilde{x}_n\widetilde{x}_{n+1}-
\widetilde{x}_{n+1}\widetilde{x}_n,\\
F_{n-1}&=\widetilde{y}_{n-1},\qquad
F_n=\widetilde{y}_n\widetilde{y}_{n-1}
-q^{-2}\widetilde{y}_{n-1}\widetilde{y}_n.
\end{split}
\end{equation*}
Direct computations lead to the following lemma.
\begin{lem}\label{technical}
{\bf (i)} For any $\sigma$,
\begin{equation*}
\begin{split}
E_n&=B_n^{\sigma}E_{n-1}-E_{n-1}B_n^{\sigma},\\
F_n&=B_n^{\sigma}F_{n-1}-q^{-2}F_{n-1}B_n^{\sigma}.
\end{split}
\end{equation*}
{\bf (ii)} The following commutation relations hold in $\ug$,
\begin{equation*}
\begin{split}
&E_{n-1}F_{n-1}=q^2F_{n-1}E_{n-1},\qquad
E_{n-1}F_n=q^2F_nE_{n-1},\\
&E_nF_{n-1}=F_{n-1}E_n,\qquad
E_nF_n=q^2F_nE_n+(q^{-2}-1)F_{n-1}E_{n-1}C_n^{-2}.
\end{split}
\end{equation*}
\end{lem}
For an ordered $r$-tuple $J=(j_1,j_2,\ldots,j_r)$ with $j_s\in \{n-1,n\}$
and $r\leq \kappa$ we write
\[F_J=F_{j_1}F_{j_2}\cdots F_{j_r},\qquad E_J=E_{j_r}\cdots E_{j_2}E_{j_1}
\]
and
\[r_{J,n-1}=\#\{s\, | \, j_s=n-1\},\qquad r_{J,n}=r-r_{J,n-1}.
\]
We furthermore define the Laurent polynomial
$\Delta_J\in\mathbb{C}[u^{\pm 1}]\subset \mathbb{C}[z^{\pm 1}]$ by
\[\Delta_J=\prod_{s,t}\left(1-
q^{2s^2+2st+2t^2}(u_{n-1}^{-1}u_n)^s(u_{n-1}^{-1}u_n^{-1})^{t}\right),
\]
where the product is taken over $(s,t)\in \mathbb{Z}_{\geq 0}^{\times 2}$
with $0<s+t\leq r$ and $t\leq r_{J,n}$.
For $J=\emptyset$, we write $F_\emptyset=E_\emptyset=1$,
$r_{\emptyset,n-1}=r_{\emptyset,n}=0$ and $\Delta_\emptyset=1$.
Using the notations of
Lemma \ref{explicit}, we can now formulate the following result.

\begin{lem}\label{polardiv}
Let $f\in F_{(0,0,\kappa)}^{\sigma,\tau}$ be a vector valued spherical
function.
Let $r\in\mathbb{Z}_{\geq 0}$ with $r\leq\kappa$ and fix an ordered
$r$-tuple $J$ with coefficients from $\{n-1,n\}$.
Then there exists a Laurent polynomial
$P_J\in\mathbb{C}[z^{\pm 1}]$ such that
\[
\Delta_J(q^\lambda)f(F_JK^{\lambda-\delta})=
P_J(q^\lambda)f|_T(q^\lambda)r_{(\kappa^{n-2},\kappa-r,\kappa+r)}
\]
for all $\lambda\in P_{2n}$.
\end{lem}
\begin{proof}
First note that the four elements
$Y=E_{n-1},E_n,F_{n-1}, F_n$ all have the same restricted weight,
\[C_jY=q^{\delta_{j,n}-\delta_{j,n-1}}YC_j,\qquad \forall \,j=1,\ldots,n.
\]
In particular, we conclude that
\[f(F_JK^{\lambda-\delta})\in\hbox{span}
\{r_{(\kappa^{n-2},\kappa-r,\kappa+r)}\},
\]
cf. Lemma \ref{explicit}, hence
\[f(F_JK^{\lambda-\delta})=f_J(q^\lambda)r_{(\kappa^{n-2},\kappa-r,\kappa+r)},
\qquad \forall\,\lambda\in P_{2n}
\]
for a unique Laurent polynomial $f_J\in \mathbb{C}[z^{\pm 1}]$. Furthermore,
we may and will
view all formulas below as identities in $\mathbb{C}[z^{\pm 1}]$
by identifying $\mathbb{C}r_{(\kappa^{n-2},\kappa-r,\kappa+r)}\simeq
\mathbb{C}$.

We thus need to prove the existence of a Laurent polynomial
$P_J\in\mathbb{C}[z^{\pm 1}]$ such that $\Delta_Jf_J=P_Jf|_T$.
To avoid unnecessary technicalities we avoid the use of
the underlying root data as much as possible. Let $\mathcal{F}_r$ be the
finite set of $r$-tuples $J=(j_1,j_2,\ldots,j_r)$ with $j_s\in\{n-1,n\}$,
and write $\mathcal{F}^\kappa=\cup_{0\leq r\leq \kappa}\mathcal{F}_r$.
For $I,J\in \mathcal{F}^\kappa$ we write $I<J$ if
$\#I\lneq \#J$ and $r_{I,n}\leq r_{J,n}$, or if
$\#I=\#J$ and $r_{I,n}\lneq r_{J,n}$.
We prove the lemma by induction
to $J\in \mathcal{F}^{\kappa}$ along $<$. For the smallest element
$J=\emptyset\in \mathcal{F}_0$, we have $\Delta_{\emptyset}=1$,
$f_{\emptyset}=f|_T$, hence we can take $P_{\emptyset}=1$.
Suppose there exists a Laurent polynomial $P_J$ such that
$\Delta_Jf_J=P_Jf|_T$ for all $J\in \mathcal{F}^\kappa$ with
$J<I=(i_1,i_2,\ldots,i_r)\in \mathcal{F}_r$. For simplicity we
write $r_1=r_{I,n-1}$ and $r_2=r_{I,n}$.

We write $g\equiv g^\prime$ for
$g,g^\prime\in\mathbb{C}[z^{\pm 1}]$ if $g-g^\prime$ lies in the
ideal of $\mathbb{C}[z^{\pm 1}]$
generated by the $f_J$'s for $J\in \mathcal{F}^\kappa$ with $J<I$.
Sometimes it is convenient to write
$g(q^\lambda)\equiv g^\prime(q^\lambda)$, where we have formally
evaluated the equivalence $g\equiv g^\prime$
of the Laurent polynomials $g$ and $g^\prime$
in $z=q^\lambda$ for arbitrary $\lambda\in P_{2n}$.

To prove the existence of the Laurent polynomial $P_I$ satisfying
$\Delta_If_I=P_If|_T$, it suffices to show that
\begin{equation}\label{indstep}
\left(1-q^{2r_1^2+2r_1r_2+2r_2^2}
(u_{n-1}^{-1}u_n)^{r_1}(u_{n-1}^{-1}u_n^{-1})^{r_2}
\right)f_I\equiv 0
\end{equation}
due to the induction hypothesis and the fact that $\Delta_J$
divides $\Delta_I$ in
$\mathbb{C}[u^{\pm 1}]\subset\mathbb{C}[z^{\pm 1}]$ when $J<I$.
We start by noting that
\begin{equation}\label{a}
\begin{split}
&f(F_{n-1}X)=-f(E_{n-1}X)+\rho(\kappa,0)_\tau(B_{n-1})\bigl(f(X)\bigr),\\
&f(F_nX)=-q^{-2}f(E_nX)+(1-q^{-2})
\rho(\kappa,0)_\tau(B_n^\tau)\bigl(f(F_{n-1}X)\bigr)
+\rho(\kappa,0)_\tau(\lbrack B_n^\tau,B_{n-1}\rbrack)
\bigl(f(X)\bigr)
\end{split}
\end{equation}
and
\begin{equation}\label{b}
\begin{split}
f(XE_{n-1})&=-f(XF_{n-1}),\\
f(XE_n)&=-f(XF_n)+(1-q^{-2})\vartheta_0(q^\sigma)f(XF_{n-1})
\end{split}
\end{equation}
for any $X\in \ug$, which follows easily
from Lemma \ref{technical}{\bf (i)} and the fact that $f$ is a vector
valued spherical function. Denote $I^\prime=(i_2,i_3,\ldots,i_r)$, which
is $<I$. Then it follows from \eqref{a} that
\[f_I(q^\lambda)\equiv
-q^{-2\delta_{i_1,n}}f(E_{i_1}F_{I^\prime}K^{\lambda-\delta}).
\]
By Lemma \ref{technical}{\bf (ii)} we now conclude that
\begin{equation*}
f_I(q^\lambda)\equiv
\begin{cases}
-q^{2(r-1)}f(F_{I^\prime}E_{i_1}K^{\lambda-\delta}),\qquad &\hbox{ if }\,\,
i_1=n-1,\\
-q^{-2}q^{2(r_2-1)}
f(F_{I^\prime}E_{i_1}K^{\lambda-\delta}),\qquad &\hbox{ if }\,\,
i_1=n.
\end{cases}
\end{equation*}
Continuing inductively while taking \eqref{b} into account, we obtain
\begin{equation}\label{oneway}
\begin{split}
f_I(q^\lambda)&\equiv (-1)^rq^{-2r_2+r_2(r_2-1)+2\sum_{s: i_s=n-1}(r-s)}
f(E_IK^{\lambda-\delta})\\
&\equiv (-1)^rq^{r_1+r_2(r_2-1)+2\sum_{s: i_s=n-1}(r-s)}
q^{-r_2\lambda_n-r_1\lambda_{n+1}+r\lambda_{n+2}}f(K^{\lambda-\delta}E_I).
\end{split}
\end{equation}
We follow a similar procedure to return from
$f(K^{\lambda-\delta}E_I)$ to $f_I(q^\lambda)$.
The starting point is
\[f(K^{\lambda-\delta}E_I)\equiv -f(K^{\lambda-\delta}E_{I^\prime}F_{i_1}),
\]
which follows from \eqref{b}. Moving $F_{i_1}$ to the left using
Lemma \ref{technical}{\bf (ii)} gives
\begin{equation*}
f(K^{\lambda-\delta}E_I)\equiv
\begin{cases}
-q^{2(r_1-1)}f(K^{\lambda-\delta}F_{i_1}E_{I^\prime}),\qquad &\hbox{ if }\,\,
i_1=n-1,\\
-q^{2(r-1)}f(K^{\lambda-\delta}F_{i_1}E_{I^\prime}),\qquad &\hbox{ if }\,\,
i_1=n.
\end{cases}
\end{equation*}
Repeating this procedure inductively gives
\begin{equation}\label{wayback}
\begin{split}
f(K^{\lambda-\delta}E_I)&\equiv (-1)^r q^{r_1(r_1-1)+2\sum_{s: i_s=n}(r-s)}
f(K^{\lambda-\delta}F_I)\\
&\equiv (-1)^rq^{r_1^2+2r_2+2\sum_{s: i_s=n}(r-s)}q^{-r\lambda_{n-1}
+r_1\lambda_n+r_2\lambda_{n+1}}f_I(q^\lambda).
\end{split}
\end{equation}
Combining \eqref{oneway} and \eqref{wayback} yields the
induction step \eqref{indstep}.
\end{proof}
\begin{lem}\label{groundstate1b}
The radial part $f|_T\in \mathbb{C}[u^{\pm 1}]$ of a vector valued
spherical function $f\in F_{(0,0,\kappa)}^{\sigma,\tau}$ is divisible
by $(1-q^{2r}u_{n-1}^{-1}u_n^{-1})$ for $r=1,\ldots,\kappa$.
\end{lem}
\begin{proof}
It is easy to check that $\Delta_{(n^\kappa)}$ is divisible by
$1-q^{2r}u_{n-1}^{-1}u_n^{-1}$ in $\mathbb{C}[u^{\pm 1}]$ for
$r=1,\ldots,\kappa$.
By the previous lemma it thus suffices to show that $P_{(n^\kappa)}$
is nonzero and relative prime to $1-q^{2r}u_{n-1}^{-1}u_n^{-1}$
for $r=1,\ldots,\kappa$. Clearly it is enough to prove this
in the classical limit $q=1$. We thus consider $q=e^h$ as a
formal parameter and we repeat part of the computations of the
proof of Lemma \ref{polardiv} for $J=(n^r)$, now modding out
to $h\mathbb{C}[[h]][z^{\pm 1}]$.

We simplify notations by writing
\begin{equation*}
\begin{split}
f(F_n^sK^{\lambda-\delta})&=
f_s(q^\lambda)r_{(\kappa^{n-2},\kappa-s,\kappa+s)},\\
f(E_nF_n^{s-1}K^{\lambda-\delta})&=
g_s(q^\lambda)r_{(\kappa^{n-2},\kappa-s,\kappa+s)}
\end{split}
\end{equation*}
for all $\lambda\in P_{2n}$ and $s\in\{1,\ldots,\kappa\}$, with
$f_s,g_s\in\mathbb{C}[[h]][z^{\pm 1}]$.
Lemma \ref{explicit} gives
\[\rho(\kappa,0)_\tau(\lbrack B_n^{\tau},B_{n-1}\rbrack)
f(F_{n}^{r-1}K^{\lambda-\delta})=
q^{2-\tau-3r}\lbrack r-1-\kappa\rbrack_q
f_{r-1}(q^\lambda)r_{(\kappa^{n-2},\kappa-r,\kappa+r)},
\]
with
\[\lbrack \alpha \rbrack_q=\frac{q^\alpha-q^{-\alpha}}{q-q^{-1}}
\]
the (symmetric) $q$-number. Identifying
$\mathbb{C}r_{(\kappa^{n-2},\kappa-r,\kappa+r)}\simeq \mathbb{C}$ as in the
proof of Lemma \ref{polardiv}, we obtain from the second formula
of \eqref{a},
\[
f_r=-q^{-2}g_{r-1}+q^{2-\tau-3r}\lbrack r-1-\kappa\rbrack_q
f_{r-1}\qquad \hbox{mod}\, h.
\]
By Lemma \ref{technical}{\bf (ii)}, we can move
$E_n$ to the other side in the expression
\[g_{r-1}(q^\lambda)=f(E_nF_n^{r-1}K^{\lambda-\delta}),
\]
after which we can use the second formula in \eqref{b} to
replace $E_n$ by $-F_n$ modulo $h$. This yields
\[g_{r-1}=-q^{2(r+1)}u_{n-1}^{-1}u_n^{-1}f_r\qquad \hbox{mod}\, h.
\]
Consequently,
\[(1-q^{2r}u_{n-1}^{-1}u_n^{-1})f_r=
q^{2-\tau-3r}\lbrack r-1-\kappa\rbrack_qf_{r-1}\qquad \hbox{mod}\, h,
\]
or, more precisely,
\[(1-u_{n-1}^{-1}u_n^{-1})f_r=(r-1-\kappa)f_{r-1}\qquad \hbox{mod}\, h
\]
for $r=1,\ldots,\kappa$. Thus $P_{(n^\kappa)}$
is nonzero, and
\begin{equation}\label{modh}
P_{(n^\kappa)}=C(1-u_{n-1}^{-1}u_n^{-1})^{-\kappa}
\Delta_{(n^\kappa)}\qquad \hbox{mod}\, h,
\end{equation}
for some nonzero constant $C\in\mathbb{C}$. Using the explicit expression
of $\Delta_{(n^\kappa)}$, it is easy to check that the right hand
side of \eqref{modh} (mod $h$) is relative prime to
$1-u_{n-1}^{-1}u_n^{-1}$ in
$\mathbb{C}[z^{\pm 1}]$, which completes the proof of the lemma.
\end{proof}
We are now in a position to determine the explicit form
of the radial part of the ground state when $\kappa_1=\kappa_2=0$.

\begin{cor}\label{kgroundstate}
Let $f\in F_{(0,0,\kappa)}^{\sigma,\tau}(0)$ be a ground state.
With a suitable
normalization of $f$, we have
\begin{equation}\label{gs2}
f|_T=u^{\delta(\kappa,0)}\prod_{1\leq i<j\leq n}
\bigl(q^2u_i^{-1}u_j,q^2u_i^{-1}u_j^{-1};q^2\bigr)_{\kappa}
\end{equation}
in $\mathbb{C}[u^{\pm 1}]$.
\end{cor}
\begin{proof}
We write $\mathcal{I}\in\mathbb{C}[u^{\pm 1}]$ for the right hand
side of \eqref{gs2}. Lemma \ref{groundstate1} and Lemma \ref{groundstate1b}
show that $f|_T$ is divisible by $\mathcal{I}$ in $\mathbb{C}[u^{\pm 1}]$,
so $f|_T=\mathcal{I}p$ for some
Laurent polynomial $p\in \mathbb{C}[u^{\pm 1}]$. It remains to show
that $p$ is a nonzero constant. For this we expand both sides of
the equality $f|_T=\mathcal{I}p$ in monomials $u^\mu$ ($\mu\in\Lambda_n$).

We write
\[p=\sum_\nu e_{\mu} u^\nu,
\]
and we denote $\nu_+\in\Lambda_n$ (respectively $\nu_-\in\Lambda_n$)
for a maximal (respectively minimal) element
with respect to the dominance order
$\leq$ such that $e_{\nu_+}\not=0\not=e_{\nu_-}$.
Observe that
\[\mathcal{I}=\sum_{-\delta(\kappa,0)\leq\nu\leq\delta(\kappa,0)}c_\nu u^\nu
\]
with $c_\nu\in \mathbb{C}$ and
$c_{-\delta(\kappa,0)}\not=0\not=c_{\delta(\kappa,0)}$.
Consequently, in the expansion of $\mathcal{I}p$ in monomials $u^\mu$
($\mu\in\Lambda_n$), the coefficient of
$u^{\delta(\kappa,0)+\nu_+}$
and of $u^{-\delta(\kappa,0)-\nu_-}$ is nonzero.

On the other hand, Lemma \ref{groundstateform} implies that
\[f|_T=\sum_{-\delta(\kappa,0)\leq\nu\leq\delta(\kappa,0)}d_\nu u^\nu
\]
with $d_\nu\in \mathbb{C}$ and
$d_{-\delta(\kappa,0)}\not=0\not=d_{\delta(\kappa,0)}$.
Hence the equality $f|_T=\mathcal{I}p$ forces $\nu_+\leq 0$ and $\nu_-\geq 0$.
This in turn implies $\nu_+=\nu_-=0$, hence
$p$ is a nonzero constant.
\end{proof}


\subsection{The ground state for $\kappa=0$.}

We start with the following preliminary lemma, which
will be convenient for several computations in this subsection.
Recall the definition of $\vartheta_k(s)$, see \eqref{lambdak}.

\begin{lem}\label{dynamical}
Let $W_1$ and $W_2$ be left $\ug$-modules and let $k\in\mathbb{Z}$.

{\bf (i)} If $w_0\in W_1$ is a vector satisfying $B_n^\sigma w_0=
\vartheta_0(q^{\sigma+k})w_0$ and $C_nw_0=w_0$, then
\[B_n^\sigma(w_0\otimes w^\prime)=w_0\otimes
B_n^{\sigma+k}w^\prime,\qquad \forall\,w^\prime\in W_2
\]
in the $\ug$-module $W_1\otimes W_2$.

{\bf (ii)} If $w_{\pm}\in W_1$ is a
vector satisfying $B_n^\sigma w_{\pm}=\vartheta_k(\pm q^{\pm
\sigma})w_{\pm}$ and $C_n w_{\pm}=q^kw_{\pm}$, then
\[B_n^\sigma(w_{\pm}\otimes w^\prime)=w_{\pm}\otimes
q^{-k}B_n^{\sigma\mp k}w^\prime,\qquad \forall\,w^\prime\in W_2
\]
in the $\ug$-module $W_1\otimes W_2$.
\end{lem}
\begin{proof}
Recall that the generator $B_n^\sigma\in\AAA_\sigma$ is explicitly given by
\[B_n^\sigma=y_nK_{n+1}^{-1}K_n^{-1}+K_n^{-1}x_nK_n^{-1}+
\vartheta_0(q^\sigma)K_n^{-2}.
\]
The proof of {\bf (i)} follows now immediately from the explicit
formula for $\Delta(B_n^\sigma)$ (see Proposition \ref{coalgebra})
and from the fact that $C_n=C_{n+1}$.
For the proof of {\bf (ii)}, let $w_{\pm}\in W_1$ be
as indicated in the lemma, and let $w^\prime\in W_2$.
Then Proposition \ref{coalgebra} implies
\begin{equation*}
\begin{split}
B_n^\sigma\bigl(w_{\pm}\otimes w^\prime\bigr)&=
w_{\pm}\otimes q^{-k}\bigl(y_nK_{n+1}^{-1}K_n^{-1}+K_n^{-1}x_nK_n^{-1}+
q^k\vartheta_k(\pm q^{\pm \sigma})K_n^{-2}\bigr)w^\prime\\
&=w_{\pm}\otimes q^{-k}B_n^{\sigma\mp k}w^\prime
\end{split}
\end{equation*}
since $q^{k}\vartheta_k(\pm q^{\pm \sigma})=\vartheta_0(q^{\sigma\mp k})$.
\end{proof}

In this subsection the realization of the fundamental,
simple $\ug$-modules within the
$q$-exterior algebra of the vector
representation of $\ug$ is used to explicitly construct the one-dimensional
$\AAA_\sigma$-modules $V(\kappa_1)_\sigma$ and $V(\kappa_2)_\sigma$
within simple $\ug$-modules.
As a consequence we obtain an
explicit expression for the radial part of the ground state for $\kappa=0$.

The vector representation of $\ug$ is the $2n$-dimensional
vector space $V$ with linear basis
$v_i$ ($i=1,\ldots,2n$) and with $\ug$-action defined by
\begin{equation}\label{vectoractie}
K_r^{\pm 1}v_j=q^{\pm \delta_{r,j}}v_j,\qquad
x_iv_j=\delta_{i+1,j}\,v_{j-1},\qquad
y_iv_j=\delta_{i,j}\,v_{j+1}
\end{equation}
for $r,j=1,\ldots,2n$ and $i=1,\ldots,2n-1$, with the convention that
$v_0=v_{2n+1}=0$. This is a realization of
the finite dimensional irreducible highest weight module
$L_{\epsilon_1}$ of $\ug$, with highest weight vector $v_1$.
Similarly, the dual vector representation is the
$2n$-dimensional vector space $V^*$ with basis $v_i^*$ ($i=1,\ldots,2n$)
and with $\ug$-action given explicitly by
\[
K_r^{\pm 1}v_j^*=q^{\mp \delta_{r,j}}v_j^*,\qquad
x_iv_j^*=-q^{-1}\delta_{i,j}\,v_{j+1}^*,\qquad
y_iv_j^*=-q\delta_{i+1,j}\,v_{j-1}^*
\]
for $r,j=1,\ldots,2n$ and $i=1,\ldots,2n-1$, with the convention that
$v_0^*=v_{2n+1}^*=0$ (it corresponds to the linear dual $V^*$ of $V$
viewed as $\ug$-module by $(X\phi)(v)=\phi(S(X)v)$ for
$\phi\in V^*$, $v\in V$ and $X\in \ug$).
This is a realization of
the finite dimensional irreducible highest weight module
$L_{-\epsilon_{2n}}$ of $\ug$, with highest weight vector $v_{2n}^*$.
For $s\in \mathbb{R}\setminus\{0\}$ we define $n$-dimensional subspaces
$W(s)=\hbox{span}\{w_i(s)\}_i\subseteq V$ and
$W^*(s)=\hbox{span}\{w_i^*(s)\}_i\subseteq V^*$ by
\[
w_i(s)=v_i+s^{-1}v_{2n+1-i},\qquad\qquad
w_i^*(s)=v_i^*+s^{-1}q^{2i-2n-1}\,v_{2n+1-i}^*
\]
for $i=1,\ldots,n$.

\begin{prop}\label{branchingvector}
The vector spaces $W(\pm q^{\pm \sigma})$ \textup{(}respectively
$W^*(\pm q^{\pm \sigma})$\textup{)} are inequivalent
simple $\AAA_\sigma$-submodules
of $V^\sigma$ \textup{(}respectively $V^*{}^\sigma$\textup{)}.
The irreducible decompositions of the semisimple
$\AAA_\sigma$-modules $V^\sigma$ and $V^*{}^\sigma$ are
\[V^\sigma=W(q^{\sigma})\oplus W(-q^{-\sigma}),\qquad\quad
V^*{}^\sigma=W^*(q^{\sigma})\oplus W^*(-q^{-\sigma}).
\]
\end{prop}
\begin{proof}
The action of the algebraic generators of $\mathcal{E}\subset \AAA_\sigma$
(see Lemma \ref{diagonal}) on
$w_i(s)$ and $w_i^*(s)$ is given by the formulas
\begin{equation}\label{Eactie}
\begin{split}
C_jw_i(s)&=q^{\delta_{i,j}}w_i(s),\qquad C_jw_i^*(s)=
q^{-\delta_{i,j}}w_i^*(s)\qquad\,\,\,\,\,\,\, (i,j=1,\ldots,n),\\
B_iw_i(s)&=w_{i+1}(s),\qquad\,\, B_iw_{i+1}^*(s)=-q^2w_i^*(s)
\qquad\,\,\,\,\,\,\,
(i=1,\ldots,n-1),\\
B_{2n+1-i}\,w_i(s)&=w_{i-1}(s),\qquad\,\,
B_{2n+1-i}\,w_{i-1}^*(s)=-w_i^*(s)\qquad
(i=2,\ldots,n),\\
B_jw_i(s)&=0,\qquad\qquad\quad B_jw_i^*(s)=0\quad
\qquad\qquad\qquad\quad \hbox{otherwise}.
\end{split}
\end{equation}
Hence $W(s)$ (respectively $W^*(s)$) is a $\mathcal{E}$-submodule of
$V$  (respectively $V^*$) which is isomorphic to the vector representation
(respectively dual vector representation) of
$U_q(\mathfrak{g}\mathfrak{l}(n))\simeq \mathcal{E}$,
cf. Lemma \ref{diagonal}. In particular, $W(s)$ and $W^*(s)$
are simple $\mathcal{E}$-modules.

Recall the notation \eqref{lambdak}. By a
direct computation one verifies that
\begin{equation}\label{Baction}
B_n^\sigma w_i(s)=\vartheta_{0}(q^\sigma)\,w_i(s),
\qquad\qquad B_n^\sigma
w_i^*(s)=\vartheta_{0}(q^\sigma)\,w_i^*(s)
\end{equation}
for all $i=1,\ldots,n-1$, and
\begin{equation}\label{extranumber}
B_n^\sigma w_n(\pm q^{\pm\sigma})=
\vartheta_{1}(\pm q^{\pm \sigma})\,w_n(\pm q^{\pm\sigma}),\qquad\quad
B_n^\sigma w_n^*(\pm q^{\pm\sigma})=
\vartheta_{-1}(\pm q^{\pm\sigma})\,w_n^*(\pm q^{\pm\sigma}).
\end{equation}
It is now clear that $W(\pm q^{\pm \sigma})\subset
V^\sigma$ are simple $\AAA_\sigma$-submodules, and that $V^\sigma=
W(q^\sigma)\oplus W(-q^{-\sigma})$ (and similarly for $V^*{}^\sigma$).
To prove that $W(q^\sigma)\not\simeq W(-q^{-\sigma})$ (resp.
$W^*(q^{\sigma})\not\simeq W^*(-q^{-\sigma})$) as $\AAA_\sigma$-modules,
it suffices to observe that the spectrum of $B_n^\sigma$
is different. This follows from the fact that
$\vartheta_k(q^\sigma)=\vartheta_k(-q^{-\sigma})$ for
$\sigma\in\mathbb{R}$ implies $k=0$.
\end{proof}

We now consider the $q$-exterior algebras of $V$ and $V^*$, cf. \cite{NYM}.
Let $T(V)$ and $T(V^*)$ be the tensor algebras of $V$ and $V^*$, respectively.
Let $J\subset T(V)$ be the graded, two-sided ideal
generated by the tensors $v_i\otimes v_i$ and
$v_r\otimes v_s+q^{-1}v_s\otimes v_r$ for $i=1,\ldots,2n$ and
$1\leq r<s\leq 2n$. Similarly, we write $J_*\subset T(V^*)$
for the graded, two-sided ideal generated by the tensors
$v_i^*\otimes v_i^*$ and $v_s^*\otimes v_r^*+q^{-1}v_r^*\otimes
v_s^*$ for $i=1,\ldots,2n$ and $1\leq r<s\leq 2n$.
The exterior algebras
$\Lambda(V)$ and $\Lambda(V^*)$ are defined to be the graded algebras $T(V)/J$
and $T(V^*)/J_*$, respectively.
The action of $\ug$ on $V$ and $V^*$ naturally
extends to a grading preserving, left $\ug$-module algebra
structure on the
exterior algebras $\Lambda(V)$ and $\Lambda(V^*)$.
We denote $\wedge$ for the products in $\Lambda(V)$ and $\Lambda(V^*)$.
Let $\Lambda^m(V)$ and $\Lambda^m(V^*)$ be the
$m$th graded pieces of $\Lambda(V)$ and $\Lambda(V^*)$ respectively,
then
\[\Lambda(V)=\bigoplus_{m=0}^{2n}\Lambda^m(V),\qquad
\Lambda(V^*)=\bigoplus_{m=0}^{2n}\Lambda^m(V^*)
\]
and
\[\Lambda^m(V)\simeq L_{(1^m,0^{2n-m})}, \qquad\Lambda^m(V^*)\simeq
L_{(0^{2n-m},-1{}^m)}
\]
as $\ug$-modules. For $I=\{i_1,\ldots,i_m\}\subseteq \{1,\ldots,2n\}$
with $i_1<i_2<\cdots<i_m$ we define
\[v_I=v_{i_1}\wedge v_{i_2}\wedge \cdots\wedge v_{i_m},\qquad
v_I^*=v_{i_1}^*\wedge v_{i_2}^*\wedge\cdots\wedge v_{i_m}^*.
\]
Then $\{v_I \, | \, \#I=m\}$ and $\{v_I^* \, | \, \#I=m\}$ are linear
bases of $\Lambda^m(V)$ and $\Lambda^m(V^*)$,
respectively. The action of $\ug$ on these basis elements is
explicitly given by
\begin{equation*}
\begin{split}
K_i^{\pm 1}v_I&=q^{\pm \delta(i,I)}v_I,\qquad
x_iv_I=v_{(I\setminus\{i+1\})\cup\{i\}},\qquad
y_iv_I=v_{(I\setminus\{i\})\cup\{i+1\}},\\
K_i^{\pm 1}v_I^*&=q^{\mp \delta(i,I)}v_I^*,\qquad
x_iv_I^*=-q^{-1}v_{(I\setminus\{i\})\cup\{i+1\}}^*,\qquad
y_iv_I^*=-qv_{(I\setminus\{i+1\})\cup\{i\}}^*,
\end{split}
\end{equation*}
where $\delta(i,I)=1$ if $i\in I$ and $=0$ otherwise,
and where $v_{(I\setminus\{k\})\cup\{l\}}=0$ and
$v_{(I\setminus\{k\})\cup\{l\}}^*=0$ if $k\not\in I$ or if $l\in I$.
Note in particular that $v_{\{1,\ldots,m\}}\in\Lambda^m(V)$
and $v_{\{2n-m+1,\ldots,2n\}}^*\in\Lambda^m(V^*)$ are highest weight
vectors. The formulas
\[\langle v_I,v_J\rangle=\delta_{I,J},\qquad
\langle v_I^*,v_J^*\rangle_*=q^{2\sum_{i\in I}i}\,\delta_{I,J}
\]
define $*$-unitary scalar products on
$\Lambda(V)$ and $\Lambda(V^*)$, respectively.

Using the notation \eqref{weightA}, we observe that
any one-dimensional $\AAA_\sigma$-module in $\Lambda^n(V)^\sigma$ (resp.
$\Lambda^n(V^*)^\sigma$)
necessarily lies in the subspace
$\Lambda^n(V)_{(1^n)}$ (resp. $\Lambda^n(V^*)_{(-1^n)}$).
This follows easily using Lemma \ref{diagonal}, since an
one-dimensional $U_q(\mathfrak{g}\mathfrak{l}(n))$-module is of highest
weight $(m^n)\in P_n^+$ for some integer $m$.
For $\underline{s}=(s_1,\ldots,s_n)$ an $n$-tuple of nonzero reals, we
define $\xi(\underline{s})\in \Lambda^n(V)_{(1^n)}$ and
$\xi^*(\underline{s})\in\Lambda^n(V^*)_{(-1^n)}$ by
\begin{equation*}
\begin{split}
\xi(\underline{s})&=w_1(s_1)\wedge w_2(s_2)\wedge\cdots\wedge
w_n(s_n),\\
\xi^*(\underline{s})&=w_1^*(s_1)\wedge w_2^*(s_2)\wedge \cdots
\wedge w_n^*(s_n).
\end{split}
\end{equation*}
It is also
convenient to use the opposite elements
$\xi_{op}(\underline{s})\in\Lambda^n(V)_{(1^n)}$ and
$\xi_{op}^*(\underline{s})\in\Lambda^n(V^*)_{(-1^n)}$, defined by
\begin{equation*}
\begin{split}
\xi_{op}(\underline{s})&=w_n(s_n)\wedge\cdots\wedge w_2(s_2)\wedge
w_1(s_1),\\
\xi^*_{op}(\underline{s})&=w_n^*(s_n)\wedge \cdots
\wedge w_{2}^*(s_{2})\wedge w_1^*(s_1).
\end{split}
\end{equation*}
For $1\leq i<j\leq n$ we have
\begin{equation*}
\begin{split}
 w_i(s_i)\wedge w_j(s_j)&=-q^{-1}w_j(s_j)\wedge w_i(q^{-2}s_i),\\
w_i^*(s_i)\wedge w_j^*(s_j)&=-qw_j^*(s_j)\wedge w_i^*(q^{2}s_i)
\end{split}
\end{equation*}
in $\Lambda(V)$ (resp. $\Lambda(V^*)$),
hence $\xi$ and $\xi_{op}$ (resp. $\xi^*$ and $\xi^*_{op}$)
are related by the formulas
\begin{equation}\label{relationopposite}
\begin{split}
\xi(q^{2(n-1)}s_1,\ldots,q^{2}s_{n-1},s_n)&=
\bigl(-q\bigr)^{n(1-n)/2}\xi_{op}(\underline{s}),\\
\xi^*(q^{-2(n-1)}s_1,\ldots,q^{-2}s_{n-1},s_n)&=
\bigl(-q\bigr)^{n(n-1)/2}\xi_{op}^*(\underline{s}).
\end{split}
\end{equation}
The norms of such vectors can be evaluated explicitly
by an easy induction argument. The result is
\begin{equation}\label{norms}
\begin{split}
\langle \xi(\underline{s}),\xi(\underline{t})\rangle&=
\prod_{i=1}^n\Bigl(1+s_i^{-1}t_i^{-1}q^{2(n-i)}\Bigr),\\
\langle \xi^*(\underline{s}),\xi^*(\underline{t})\rangle_*&=
q^{n(n+1)}\prod_{i=1}^n\Bigl(1+s_i^{-1}t_i^{-1}q^{2(i-n)}\Bigr).
\end{split}
\end{equation}
We use the abbreviation $\xi(t)=\xi((t^n))$ for
$t\in \mathbb{R}\setminus\{0\}$, where $(t^n)$ is the $n$-tuple with
$t$ in each entry. Similarly we write $\xi_{op}(t)$, $\xi^*(t)$ and
$\xi^*_{op}(t)$.
\begin{lem}\label{59}
{\bf (i)} Let $\underline{s}\in\bigl(\mathbb{R}\setminus \{0\}\bigr)^n$.
Then
$\hbox{span}_{\mathbb{C}}\{\xi_{op}(\underline{s})\}\subset
\Lambda^n(V)$ and $\hbox{span}_{\mathbb{C}}\{\xi^*_{op}(
\underline{s})\}\subset \Lambda^n(V^*)$ are
one-dimensional $\mathcal{E}$-submodules
if and only if $\underline{s}=(t^n)$ for some
$t\in\mathbb{R}\setminus\{0\}$.

{\bf (ii)} The two subspaces
$\hbox{span}\{\xi_{op}(\pm q^{\pm \sigma})\}\subset
\Lambda^n(V)^\sigma$ are one-dimensional $\AAA_\sigma$-submodules of
$\Lambda^n(V)^\sigma$.
The generator $B_n^\sigma\in\AAA_\sigma$
acts by
\[B_n^\sigma\,\xi_{op}(\pm q^{\pm\sigma})=
\vartheta_1(\pm q^{\pm \sigma})\xi_{op}(\pm q^{\pm \sigma}).
\]

{\bf (iii)} The two subspaces $\hbox{span}\{\xi^*_{op}(\pm
q^{\pm \sigma})\}\subset \Lambda^n(V^*)^\sigma$
are one-dimensional $\AAA_\sigma$-submodules of
$\Lambda^n(V^*)^\sigma$.
The generator $B_n^\sigma\in\AAA_\sigma$ acts by
\[
B_n^\sigma\,\xi^*_{op}(\pm q^{\pm\sigma})=
\vartheta_{-1}(\pm q^{\pm\sigma})\xi^*_{op}(\pm q^{\pm \sigma}).
\]
\end{lem}
\begin{proof}
{\bf (i)} Fix $i\in\{1,\ldots,n-1\}$.
By \eqref{vectoractie} we have
\[x_{2n-i}\,w_j(s)=y_i\,w_j(s)=0\qquad \forall
j\in\{1,\ldots,n\}\setminus \{i\}.
\]
Combined with \eqref{Eactie} and
Proposition \ref{coalgebra} we obtain
\begin{equation*}
\begin{split}
B_i\,\xi(\underline{s})&=w_1(s_1)\wedge\cdots\wedge
w_{i-1}(s_{i-1})\wedge B_i\,w_i(s_i)\wedge K_i^{-1}K_{2n-i}^{-1}\bigl(
w_{i+1}(s_{i+1})\wedge\cdots\wedge w_n(s_n)\bigr)\\
&=w_1(s_1)\wedge\cdots\wedge w_{i-1}(s_{i-1})
\wedge w_{i+1}(s_i)\wedge w_{i+1}(qs_{i+1})\wedge
w_{i+2}(s_{i+2})\wedge\cdots\wedge w_n(s_n).
\end{split}
\end{equation*}
Since
\[w_{i+1}(s_i)\wedge w_{i+1}(qs_{i+1})=(s_i^{-1}-q^{-2}s_{i+1}^{-1})\,
v_{2n-i}\wedge v_{i+1}
\]
in $\Lambda(V)$, we see that $B_i\,\xi(\underline{s})=0$
for all $i=1,\ldots,n-1$ if and only if $\underline{s}=(tq^{2(n-1)},
\ldots,tq^{2},t)$ for some $t\in\mathbb{R}\setminus\{0\}$.
Furthermore, $C_j\xi(\underline{s})=q\xi(\underline{s})$
for $j=1,\ldots,n$ by \eqref{Eactie}.
Then Lemma \ref{diagonal} and \eqref{relationopposite}
imply that $\xi_{op}(\underline{s})$ is a lowest
$\mathcal{E}\simeq U_q(\mathfrak{g}\mathfrak{l}(n))$-weight
vector of weight $(1^n)$ in $\Lambda^n(V)$
if and only if $\underline{s}=(t^n)$ for some $t\in\mathbb{R}\setminus\{0\}$,
which implies {\bf (i)} for $\xi_{op}$.
The proof for $\xi^*_{op}$ is similar.

For the proof of {\bf (ii)} and {\bf (iii)}, we note that
Lemma \ref{dynamical}{\bf (ii)}, \eqref{Eactie} and \eqref{Baction}
imply
\begin{equation*}
\begin{split}
B_n^\sigma\,\xi(\underline{s})&=
w_1(s_1)\wedge w_2(s_2)\wedge\cdots w_{n-1}(s_{n-1})\wedge
B_n^\sigma\,w_n(s_n),\\
B_n^\sigma\,\xi^*(\underline{s})&=
w_1^*(s_1)\wedge w_2^*(s_2)\wedge\cdots w_{n-1}^*(s_{n-1})
\wedge B_n^\sigma\,w_n^*(s_n).
\end{split}
\end{equation*}
Formula \eqref{extranumber} now completes the proof.
\end{proof}

In the following proposition we consider
the $\ug$-module $\Lambda^n(V)\otimes \Lambda^n(V^*)$. By the
Pieri formula (cf. \cite{DS}) its decomposition in irreducibles is given by
\begin{equation}\label{Pieri}
\Lambda^n(V)\otimes \Lambda^n(V^*)\simeq \bigoplus_{m=0}^n
L_{\omega_m^\natural},
\end{equation}
where $\omega_m=(1^m,0^{n-m})\in\Lambda_n^+$.
We denote $\pi_n:
\Lambda^n(V)\otimes\Lambda^n(V^*)\rightarrow L_{\omega_n^\natural}$
for the projection onto
$L_{\omega_n^\natural}$ along the decomposition \eqref{Pieri}.

Since $\delta(0,\kappa_1)=\kappa_1\omega_n$, \eqref{Pieri}
and our earlier analysis of branching rules easily imply that
the $V(\pm 1)_\sigma$-isotypical component of
$(\Lambda^n(V)\otimes \Lambda^n(V^*))^\sigma$ is
one-dimensional (and necessarily lies in the component
$L_{\omega_n^\natural}$), while the $V(0)_\sigma$-isotypical component of
$(\Lambda^n(V)\otimes \Lambda^n(V^*))^\sigma$ is $(n+1)$-dimensional.
In \cite{DS}, the $V(0)_\sigma$-isotypical component
was constructed using the explicit solution $J^\sigma$
of the reflection equation. In the following proposition
we give an explicit construction of one particular vector in the
$V(0)_\sigma$-isotypical component and we explicitly
construct the $V(\pm 1)_\sigma$-isotypical components of
$(\Lambda^n(V)\otimes \Lambda^n(V^*))^\sigma$.

\begin{prop}
Define $\eta_k(s)\in\Lambda^n(V)\otimes
\Lambda^n(V^*)$ \textup{(}$s\in\mathbb{R}\setminus\{0\}$,
$k\in\{-1,0,1\}$\textup{)}
by the formulas
\begin{equation*}
\begin{split}
\eta_{-1}(s)&=\xi_{op}(s)\otimes\xi^*_{op}(-qs^{-1}),\\
\eta_0(s)&=\xi_{op}(s)\otimes\xi^*_{op}(q^{-1}s),\\
\eta_{1}(s)&=
\xi_{op}(-s^{-1})\otimes \xi^*_{op}(qs).
\end{split}
\end{equation*}
Then $\mathbb{C}\{\eta_k(q^\sigma)\}\subset
\bigl(\Lambda^n(V)\otimes\Lambda^n(V^*)\bigr)^\sigma$
is a one-dimensional $\AAA_\sigma$-submodule, isomorphic to
$V(k)_\sigma$. Furthermore,
$\mathbb{C}\{\eta_k(q^\sigma)\}$ \textup{(}$k=\pm 1$\textup{)} and
$\pi_n(\mathbb{C}\{\eta_0(q^\sigma)\})$
are realizations of the one-dimensional $\AAA_\sigma$-modules
$V(k)_\sigma$ \textup{(}$k=\pm 1$\textup{)}
and $V(0)_\sigma$ in $L_{\omega_n^\natural}^\sigma$,
respectively.
\end{prop}
\begin{proof}
By Proposition \ref{coalgebra} and the previous lemma,
$\mathbb{C}\{\xi_{op}(s)\otimes \xi^*_{op}(t)\}\subset\Lambda^n(V)\otimes
\Lambda^n(V^*)$ is a copy of the trivial
$\mathcal{E}\simeq U_q(\mathfrak{g}\mathfrak{l}(n))$ module
for all $s,t\in\mathbb{R}\setminus\{0\}$.
Furthermore, by Lemma \ref{dynamical} and the previous lemma,
\[B_n^\sigma (\xi_{op}(\pm q^{\pm \sigma})\otimes \xi^*_{op}(t))=
\xi_{op}(\pm q^{\pm\sigma})\otimes q^{-1}B_n^{\sigma\mp 1}(\xi^*_{op}(t))
\]
in $\Lambda^n(V)\otimes \Lambda^n(V^*)$ for all
$t\in\mathbb{R}\setminus\{0\}$. The first statement
now follows from Lemma \ref{59} and the equalities
\[q^{-1}\vartheta_{-1}(-q^{-\sigma+1})=\vartheta_0(q^{\sigma-2}),\qquad
q^{-1}\vartheta_{-1}(q^{\sigma-1})=\vartheta_0(q^\sigma),\qquad
q^{-1}\vartheta_{-1}(q^{\sigma+1})=\vartheta_0(q^{\sigma+2}),
\]
which are the scalars by which $B_n^\sigma$ acts on $V(-1)_\sigma$,
$V(0)_\sigma$ and $V(1)_\sigma$, respectively.

We have already observed
that $\mathbb{C}\{\eta_{\pm 1}(q^\sigma)\}$ automatically
lies in the component $L_{\omega_n^\natural}^\sigma$ of the decomposition
\eqref{Pieri}.
In the decomposition of $\eta_0(q^\sigma)$
as sum of basis elements $v_I\otimes v_J^*$ ($\#I=\#J=n$),
the coefficient of $v_{\{1,\ldots,n\}}\otimes v_{\{n+1,\ldots,2n\}}^*$
is nonvanishing, hence $\pi_n(\eta_0(q^\sigma))\not=0$. This completes
the proof of the second statement of the proposition.
\end{proof}
We define a $*$-unitary
scalar product on $\Lambda^n(V)\otimes \Lambda^n(V^*)$ by
\[\langle v\otimes v^*,w\otimes w^*\rangle_1=\langle v,w\rangle
\langle v^*,w^*\rangle_*.
\]
The previous proposition implies that
\[f_k^{\sigma,\tau}(\cdot)
=\langle \cdot\,\pi_n(\eta_k(q^\sigma)),\eta_1(q^\tau)\rangle_1=
\langle \cdot\, \eta_k(q^\sigma),\eta_1(q^\tau)\rangle_1
\in F_{(1,k,0)}^{\sigma,\tau}(0)
\]
for $k=-1,0,1$ are ground states.
Their radial parts $f_{k}^{\sigma,\tau}|_T\in\mathbb{C}[u^{\pm 1}]$
can now be explicity evaluated as follows.
\begin{lem}
After suitable normalization, we have
\begin{equation*}
\begin{split}
f_{-1}^{\sigma,\tau}|_T&=
\prod_{i=1}^n\bigl(1-q^{1-\sigma+\tau}u_i^{-1}\bigr)
\bigl(1-q^{-3+\sigma-\tau}u_i\bigr),\\
f_{0}^{\sigma,\tau}|_T&=
\prod_{i=1}^n\bigl(1-q^{1-\sigma+\tau}u_i^{-1}\bigr)
\bigl(1+q^{-1-\sigma-\tau}u_i\bigr),\\
f_{1}^{\sigma,\tau}|_T&=
\prod_{i=1}^n\bigl(1+q^{1+\sigma+\tau}u_i^{-1}\bigr)
\bigl(1+q^{-3-\sigma-\tau}u_i\bigr)
\end{split}
\end{equation*}
in $\mathbb{C}[u^{\pm 1}]$.
\end{lem}
\begin{proof}
For $\lambda\in P_{2n}$, we write
$|\lambda|_n=\lambda_1+\lambda_2+\cdots +\lambda_n$.
By \eqref{relationopposite} we have
\begin{equation*}
\begin{split}
K^{\lambda-\delta}\,\xi_{op}(s)&=
(-1)^{\frac{n(n-1)}{2}}q^{|\lambda|_n-n^2}
\xi(sq^{\lambda_1-\lambda_{2n}-1},sq^{\lambda_2-\lambda_{2n-1}-1},\ldots,
sq^{\lambda_n-\lambda_{n+1}-1}),\\
K^{\lambda-\delta}\,\xi^*_{op}(s)&=
(-1)^{\frac{n(n-1)}{2}}q^{-|\lambda|_n+n^2}
\xi^*(sq^{1-\lambda_1+\lambda_{2n}},sq^{1-\lambda_2+\lambda_{2n-1}},
\ldots,sq^{1-\lambda_n+\lambda_{n+1}})
\end{split}
\end{equation*}
for all $\lambda\in P_{2n}$. By \eqref{relationopposite} and
\eqref{norms}, we conclude that
\begin{equation*}
\begin{split}
\langle K^{\lambda-\delta}\,\xi_{op}(s),\xi_{op}(t)\rangle&=
q^{|\lambda|_n-\frac{n(n+1)}{2}}\prod_{i=1}^n
\bigl(1+s^{-1}t^{-1}q^{1+\lambda_{2n+1-i}-\lambda_i}\bigr),\\
\langle K^{\lambda-\delta}\,\xi^*_{op}(s),\xi^*_{op}(t)\rangle_*&=
q^{-|\lambda|_n+\frac{3n(n+1)}{2}}
\prod_{i=1}^n\bigl(1+s^{-1}t^{-1}q^{-1+\lambda_i-\lambda_{2n+1-i}}\bigr)
\end{split}
\end{equation*}
for all $\lambda\in P_{2n}$. The lemma now easily follows from the
definition of the $\eta_k$'s.
\end{proof}

\begin{lem}\label{groundstate2}
Let $f\in F_{(\kappa_1,\kappa_2,0)}^{\sigma,\tau}(0)$ be a ground state.
After a suitable renormalization of $f$, we have
\[f|_T=u^{\delta(0,\kappa_1)}
\prod_{i=1}^n\bigl(-q^{1+\sigma+\tau}u_i^{-1};q^2\bigr)_{\kappa_1+\kappa_2}
\bigl(q^{1-\sigma+\tau}u_i^{-1};q^2\bigr)_{\kappa_1-\kappa_2}
\]
in $\mathbb{C}[u^{\pm 1}]$.
\end{lem}
\begin{proof}
We consider the $\ug$-module
$W=(\Lambda^n(V)\otimes \Lambda^n(V^*))^{\otimes \kappa_1}$,
with $*$-unitary inner product
\[\langle v_1\otimes\cdots\otimes v_{\kappa_1},
w_1\otimes\cdots\otimes w_{\kappa_1}\rangle_{\kappa_1}=
\prod_{i=1}^{\kappa_1}\langle v_i,w_i\rangle_1
\]
for all $v_i,w_j\in \Lambda^n(V)\otimes \Lambda^n(V^*)$.
Lemma \ref{dynamical} implies that the nonzero vector
\[
\chi_\tau=\eta_1(q^\tau)\otimes\eta_1(q^{\tau+2})\otimes
\cdots\otimes \eta_1(q^{\tau+2(\kappa_1-1)})\in W
\]
spans a copy of the one-dimensional $\AAA_\tau$-module
$V(\kappa_1)_\tau$ within $W^\tau$. By \eqref{Pieri} we have
\begin{equation}\label{Wdecomposition}
W\simeq L_{\delta(0,\kappa_1)^\natural}\oplus
\bigoplus_{\lambda\in P_{2n}^+: \lambda\prec
\delta(0,\kappa_1)^\natural}L_\lambda^{\oplus d_\lambda}
\end{equation}
as $\ug$-modules for certain multiplicities $d_\lambda\in
\mathbb{Z}_{\geq 0}$.  By Proposition \ref{branching},
$\mathbb{C}\{\chi_\tau\}$ lies in the component
$L_{\delta(0,\kappa_1)^\natural}$ of the decomposition
\eqref{Wdecomposition}.

For the realization of $V(\kappa_2)_\sigma$ in $W^\sigma$, we
need to consider the cases $0\leq\kappa_2\leq\kappa_1$ and
$-\kappa_1\leq\kappa_2\leq 0$ seperately.
For $0\leq\kappa_2\leq\kappa_1$, Lemma \ref{dynamical}
implies that the nonzero vector
\[\chi_\sigma^+=\eta_0(q^\sigma)^{\otimes (\kappa_1-\kappa_2)}\otimes
\eta_1(q^\sigma)\otimes \eta_1(q^{\sigma+2})\otimes\cdots\otimes
\eta_1(q^{\sigma+2(\kappa_2-1)})\in W
\]
spans a copy of the one-dimensional $\AAA_\sigma$-module
$V(\kappa_2)_\sigma$ in $W^\sigma$. If $p: W\rightarrow
L_{\delta(0,\kappa_1)^\natural}$ denotes the projection on
$L_{\delta(0,\kappa_1)^\natural}$ along the
decomposition \eqref{Wdecomposition},
then $p(\mathbb{C}\{\chi_\sigma^+\})$
is a copy of the one-dimensional $\AAA_\sigma$-module
$V(\kappa_2)_\sigma$ in $L_{\delta(0,\kappa_1)^\natural}^\sigma$ by
highest weight considerations, cf. Lemma \ref{highcontri}. Thus

\[
f(\cdot)=\langle \cdot\,p(\chi_\sigma^+),\chi_\tau\rangle_{\kappa_1}=
\langle \cdot\,\chi_\sigma^+,\chi_\tau\rangle_{\kappa_1}
\]
defines a ground state
$f\in F_{(\kappa_,\kappa_2,0)}^{\sigma,\tau}(0)$. Since
\[f|_T=
\prod_{i=1}^{\kappa_1-\kappa_2}f_{0}^{\sigma,\tau+2(i-1)}|_T
\prod_{j=1}^{\kappa_2}f_{1}^{\sigma+2(j-1),\tau+
2(\kappa_1-\kappa_2+j-1)}|_T,
\]
the result now follows from the previous lemma
by a straightforward manipulation of $q$-shifted factorials.
For $-\kappa_1\leq\kappa_2\leq 0$ the vector $\chi_\sigma^+$
should be replaced by
\[\chi_\sigma^-=\eta_0(q^\sigma)^{\otimes (\kappa_1+\kappa_2)}
\otimes \eta_{-1}(q^\sigma)\otimes\eta_{-1}(q^{\sigma-2})\otimes
\cdots\otimes\eta_{-1}(q^{\sigma-2(-\kappa_2-1)})\in W.
\]
We leave the remaining straightforward computations to the reader.
\end{proof}

Theorem \ref{maintheorem}{\bf (ii)} follows now immediately by
combining Lemma \ref{splitting}, Corollary \ref{kgroundstate} and
Lemma \ref{groundstate2}.


\section{Generalized quantum Schur orthogonality relations}

We freely use the notations
of the previous sections.
The normalized
Haar functional $h: \mathbb{C}_q[G]\rightarrow \mathbb{C}$ is the unique
linear map which intertwines the (left) regular $\ug$-action on
$\mathbb{C}_q[G]$ with the
trivial $\ug$-action on $\mathbb{C}$ and which maps $1\in
\mathbb{C}_q[G]$ to $1\in \mathbb{C}$. In particular,
$h$ is identically zero on the summands
$W(\lambda)$ ($\lambda\in P_{2n}^+\setminus \{0\}$)
in the Peter-Weyl decomposition \eqref{PW}.

The Haar functional $h$ defines
a pre-Hilbert space structure on $\mathbb{C}_q[G]$  by
\[\langle a,b\rangle_h=h(b^*a),\qquad \forall\, a,b\in \mathbb{C}_q[G].
\]
The quantum Schur
orthogonality relations
imply that the Peter-Weyl decomposition \eqref{PW} is
an orthogonal decomposition with respect to $\langle
\cdot,\cdot\rangle_h$.

We fix a scalar product $\langle
\cdot,\cdot\rangle_{\vec{\kappa}}$ on the simple $U_q(\mathfrak{k})$-module
$V(\kappa,\kappa_1)$ such that the representation map
$\rho(\kappa,\kappa_1)$
is $*$-unitary (of course, $\langle \cdot,\cdot\rangle_{\vec{\kappa}}$
only depends on the parameters $\kappa,\kappa_1$ of $\vec{\kappa}$,
and is unique up to scalar multiples).
{}From Corollary \ref{twist} it follows that the simple $\AAA_\tau$-module
$(\rho(\kappa,\kappa_1)_\tau, V(\kappa,\kappa_1)_\tau)$
is $*$-unitary with respect to
$\langle \cdot,\cdot\rangle_{\vec{\kappa}}$.
\begin{defi}
Let
$\pi=\pi_{\vec{\kappa}}^{\sigma,\tau}:
F_{\vec{\kappa}}^{\sigma,\tau}\times
F_{\vec{\kappa}}^{\sigma,\tau}\rightarrow \hbox{Hom}(\ug,\mathbb{C})$
be the sesquilinear form
\[\pi(f,g)(X)=\sum\, \langle f(K^{-\delta}X_{1}),
g(K^{-\delta}\omega(X_{2}))\rangle_{\vec{\kappa}},\qquad f,g\in
F_{\vec{\kappa}}^{\sigma,\tau},\,\,\, X\in \ug.
\]
\end{defi}
Let $f,g\in F_{\vec{\kappa}}^{\sigma,\tau}$. We write $f\cdot K^{-\delta}$
for the right regular action of $K^{-\delta}\in\ug$ on the first tensor
component of $f\in \mathbb{C}_q[G]\otimes V(\kappa,\kappa_1)_\tau$,
and similarly for $g\cdot K^{-\delta}$. If
\[f\cdot K^{-\delta}=\sum\,f_i\otimes v_i,\qquad
g\cdot K^{-\delta}=\sum\,g_j\otimes w_j,
\]
with $f_i, g_j\in \mathbb{C}_q[G]$ and
$v_i,w_j\in V(\kappa,\kappa_1)_\tau$, then the
definition of the pairing $\pi$ and \eqref{omegaster} imply that
\begin{equation}\label{writtenout}
\pi(f,g)(X)=\sum\,f_i(X_1)g_j^*(X_2)\langle
v_i,w_j\rangle_{\vec{\kappa}}=
\sum\,\bigl(f_ig_j^*\bigr)(X)\langle v_i,w_j\rangle_{\vec{\kappa}}
\end{equation}
for all $X\in \ug$. In particular, the image of $\pi$ is contained
in $\mathbb{C}_q[G]$. This observation can be strengthened as follows.
\begin{prop}\label{impi}
The image of the pairing $\pi_{\vec{\kappa}}^{\sigma,\tau}$
is contained in $\mathcal{H}^{\sigma,\tau}$.
\end{prop}
\begin{proof}
We use the characterization of the space $\mathcal{H}^{\sigma,\tau}$
as given in Remark \ref{alt}{\bf (ii)}.

Fix $f,g\in F_{\vec{\kappa}}^{\sigma,\tau}$.
Let $a\in \AAA_\sigma$ and write
$\Delta(a)=\sum a_{1}\otimes a_{2}$ with the
$a_{1}$'s from $\AAA_\sigma$ (which can be done
in view of Proposition \ref{coalgebra}).
Using that the character $\chi(\kappa_2)_\sigma:
\AAA_\sigma\rightarrow \mathbb{C}$
of $V(\kappa_2)_\sigma$ is $*$-unitary, $\chi(\kappa_2)_\sigma(a^*)=
\overline{\chi(\kappa_2)_\sigma(a)}$ for $a\in \AAA_\sigma$,
we obtain for $X\in \ug$,
\begin{equation*}
\begin{split}
\pi(f,g)(Xa)&=\sum\,\chi(\kappa_2)_\sigma(a_{1})\,
\langle f(K^{-\delta}X_{1}), g(K^{-\delta}\omega(X_{2})\omega(a_{2}))
\rangle_{\vec{\kappa}}\\
&=\sum\,\langle f(K^{-\delta}X_{1}),
g(K^{-\delta}\omega(X_{2})\omega(a_{2})a_{1}^*)\rangle_{\vec{\kappa}}\\
&=\epsilon(a)\pi(f,g)(X)
\end{split}
\end{equation*}
since
\[\sum\,\omega(a_{2})a_{1}^*=\Bigl(\sum\,a_{1}S(a_{2})\Bigr)^*=
\overline{\epsilon(a)}1
\]
by the antipode axiom for the Hopf algebra $\ug$. For the left invariance
with respect to $\AAA_\tau$, let $b\in \AAA_\tau$, and write again
$\Delta(b)=\sum\,b_1\otimes b_2$ with all $b_1$'s from $\AAA_\tau$.
Now $\Delta$ is a $*$-algebra homomorphism and $S$ is an anti-coalgebra
homomorphism, so
\[\Delta(\omega_\delta(b))=\sum\,\omega_\delta(b_2)\otimes
\omega_\delta(b_1).
\]
Furthermore, note that
$(\omega\circ\omega_\delta)(X)=K^{\delta}XK^{-\delta}$
for $X\in \ug$.
Combined with \eqref{Sinner},
we thus obtain for $X\in \ug$,
\begin{equation*}
\begin{split}
\pi(f,g)(\omega_\delta(b)X)&=
\sum\,\langle f(K^{-\delta}\omega_\delta(b_2)X_1),
g(b_1K^{-\delta}\omega(X_2))\rangle_{\vec{\kappa}}\\
&=\sum\,\langle f(K^{-2\delta}\omega(b_2)K^{\delta}X_1),
\rho(\kappa,\kappa_1)_\tau(b_1)
g(K^{-\delta}\omega(X_2))\rangle_{\vec{\kappa}}\\
&=\sum\,\langle f(b_1^*S^2(\omega(b_2))K^{-\delta}X_1),
g(K^{-\delta}\omega(X_2))\rangle_{\vec{\kappa}}\\
&=\overline{\epsilon(b)}\pi(f,g)(X),
\end{split}
\end{equation*}
since
\[\sum\,b_1^*S^2(\omega(b_2))=
\sum\,b_1^*S(b_2^*)=\epsilon(b^*)1
=\overline{\epsilon(b)}1.
\]
\end{proof}
Choose a vector $v\in \widetilde{V}(\kappa,\kappa_1)_\tau$
such that $\langle v,v\rangle_{\vec{\kappa}}=1$. Without loss of
generality we may and will assume that the implicit identification
$\widetilde{V}(\kappa,\kappa_1)_\tau\simeq\mathbb{C}$ in the
definition of the restriction map $|_T$ is realized by the
explicit map $zv\mapsto z$ ($z\in\mathbb{C}$).
\begin{cor}\label{Wsymmetry}
If $f,g\in F_{\vec{\kappa}}^{\sigma,\tau}$,
then $f|_T (g|_T)^*\in \mathbb{C}[u^{\pm 1}]^W$.
\end{cor}
\begin{proof}
Let $f,g\in F_{\vec{\kappa}}^{\sigma,\tau}$.
Since the Cartan elements $K^\lambda$ ($\lambda\in P_{2n}$)
are group-like and $*$-selfadjoint elements in $\ug$,
it follows from Lemma \ref{restriction} and from
the definition of $\pi$ that
\begin{equation*}
\begin{split}
\hbox{Res}_T(\pi(f,g))(q^{\lambda})&=\pi(f,g)(K^{\lambda})\\
&=\langle f(K^{\lambda-\delta}),
g(K^{-\lambda-\delta})\rangle_{\vec{\kappa}}\\
&=f|_T(q^{\lambda})(g|_T)^*(q^{\lambda})
\end{split}
\end{equation*}
for $\lambda\in P_{2n}$, hence
\[
\hbox{Res}_T(\pi(f,g))=f|_T (g|_T)^*
\]
in $\mathbb{C}[u^{\pm 1}]$. On the other hand,
$\hbox{Res}_T(\pi(f,g))\in \mathbb{C}[u^{\pm 1}]^W$ by the previous
proposition and by Proposition \ref{validity}.
\end{proof}
We define a linear functional
$h^{\sigma,\tau}: \mathbb{C}[u^{\pm 1}]\rightarrow \mathbb{C}$
by integrating against
the orthogonality measure of the Macdonald-Koornwinder polynomials
$P_{\mu}^{\sigma,\tau}$ ($\mu\in\Lambda_n^+$), see \eqref{Pzonal}.
We renormalize $h^{\sigma,\tau}$ so that
$h^{\sigma,\tau}(1)=1$. In particular, for real parameters
$\sigma,\tau$ satisfying
\begin{equation}\label{parametercondities}
q^{\epsilon\sigma+\epsilon^\prime \tau}<q^{-1},\qquad \forall\,
\epsilon,\epsilon^\prime\in\{\pm 1\},
\end{equation}
the functional $h^{\sigma,\tau}$ is given by
\[h^{\sigma,\tau}(p)=\frac{1}{N}\int_Tp(u)\Delta_{\sigma,\tau}(u)
\frac{du}{u},\qquad
p\in \mathbb{C}[u^{\pm 1}],
\]
with $\Delta_{\sigma,\tau}(u)=\Delta(u;-q^{\sigma+\tau+1}, -q^{-\sigma-\tau+1},
q^{\sigma-\tau+1}, q^{-\sigma+\tau+1};q^2,q^2)$ the
Macdonald-Koorn\-win\-der weight function and with positive normalization
constant $N=\int_T\Delta_{\sigma,\tau}(u)du$. Then
$h^{\sigma,\tau}$ gives rise to a pre-Hilbert
structure on $\mathbb{C}[u^{\pm 1}]$,
\[\langle p,r\rangle_{\sigma,\tau}=h^{\sigma,\tau}(r^*p),
\qquad p,r\in\mathbb{C}[u^{\pm 1}]
\]
with $*$-structure on $\mathbb{C}[u^{\pm 1}]$ defined by \eqref{starpol}.
Theorem \ref{maintheorem} for $\vec{\kappa}=\vec{0}$ and the
orthogonality relations for the
Macdonald-Koornwinder polynomials $P_\mu^{\sigma,\tau}$ imply
$h(f)=h^{\sigma,\tau}(f|_T)$ for $f\in F_{\vec{0}}^{\sigma,\tau}$.
By Proposition \ref{Hlink} this implies
\begin{equation}\label{same}
h(f)=h^{\sigma,\tau}\bigl(\hbox{Res}_T(f)\bigr),\qquad
\forall\, f\in \mathcal{H}^{\sigma,\tau}.
\end{equation}
Since $\hbox{Res}_T: \mathcal{H}^{\sigma,\tau}\rightarrow
\mathbb{C}[u^{\pm 1}]^W$ is a $*$-algebra isomorphism, we conclude
from \eqref{same} that
\[
\hbox{Res}_T: (\mathcal{H}^{\sigma,\tau}, \langle
  \cdot,\cdot\rangle_h)\rightarrow (\mathbb{C}[u^{\pm 1}]^W,
\langle \cdot,\cdot\rangle_{\sigma,\tau})
\]
is an isomorphism of pre-Hilbert spaces.

\begin{prop}\textup{(}\textup{
Generalized quantum Schur orthogonality
relations}\textup{)}. Let $\mu,\nu\in \Lambda_n^+$ and choose
elementary vector valued spherical functions
$f\in F_{\vec{\kappa}}^{\sigma,\tau}(\mu)$ and
$g\in F_{\vec{\kappa}}^{\sigma,\tau}(\nu)$.
If $\mu\not=\nu$, then $\langle f|_T, g|_T\rangle_{\sigma,\tau}=0$.
\end{prop}
\begin{proof}
Let $f,g\in F_{\vec{\kappa}}^{\sigma,\tau}$.
By the proof of Corollary \ref{Wsymmetry}, by Proposition \ref{impi}
and by \eqref{same} we have
\begin{equation}\label{innrel}
h\bigl(\pi(f,g)\bigr)=\langle f|_T, g|_T
\rangle_{\sigma,\tau}.
\end{equation}
Let $f\in F_{\vec{\kappa}}^{\sigma,\tau}(\mu)$ and
$g\in F_{\vec{\kappa}}^{\sigma,\tau}(\nu)$
with $\mu,\nu\in\Lambda_n^+$ and $\mu\not=\nu$.
We show that the left hand side of \eqref{innrel} vanishes.
Since $f\cdot K^{-\delta}\in W(\mu^\natural)\otimes
V(\kappa,\kappa_1)_\tau$ and $g\cdot K^{-\delta}\in W(\nu^\natural)\otimes
V(\kappa,\kappa_1)_\tau$, we may write
\[f\cdot K^{-\delta}=\sum\,f_i\otimes v_i,\qquad
g\cdot K^{-\delta}=\sum\,g_j\otimes w_j,
\]
with $f_i\in W(\mu^\natural)$, $g_j\in W(\nu^\natural)$ and
$v_i,w_j\in V(\kappa,\kappa_1)_\tau$.
By \eqref{writtenout} and by the quantum Schur orthogonality
relations, we conclude that
\[h\bigl(\pi(f,g)\bigr)=\sum\, h(f_ig_j^*)\langle v_i,w_j
\rangle_{\vec{\kappa}}=0,
\]
as desired.
\end{proof}

The main application of the generalized
quantum Schur orthogonality relations is the identification of the
vector valued spherical functions with
Macdonald-Koorn\-win\-der polynomials, as stated in Theorem
\ref{maintheorem}{\bf (iii)}.

\begin{cor}
Let $\mu\in\Lambda_n^+$ and choose
elementary vector valued spherical functions
$f_\mu\in F_{\vec{\kappa}}^{\sigma,\tau}(\mu)$ and $f_0\in
F_{\vec{\kappa}}^{\sigma,\tau}(0)$. Then
\[\frac{f_\mu|_T}{f_0|_T}=D\,
P_\mu\bigl(u;-q^{\sigma+\tau+1+\kappa_1+\kappa_2},
-q^{-\sigma-\tau+1},
q^{\sigma-\tau+1},q^{-\sigma+\tau+1+\kappa_1-\kappa_2};q^2,
q^{2\kappa+2}\bigr)
\]
for some nonzero constant $D$.
\end{cor}
\begin{proof}
It suffices to prove the corollary for real parameters
$\sigma,\tau$ satisfying \eqref{parametercondities}.
By Proposition \ref{triangular} we can normalize the elementary
vector valued spherical functions $f_\nu\in
F_{\vec{\kappa}}^{\sigma,\tau}(\mu)$ ($\nu\in\Lambda_n^+$)
such that
\begin{equation}\label{eig1}
p_\nu:=\frac{f_\nu|_T}{f_0|_T}=
m_{\nu}+\sum_{\stackrel{\nu^\prime\in\Lambda_n^+}{\nu^\prime<\nu}}
c_{\nu^\prime} m_{\nu^\prime}\in\mathbb{C}[u^{\pm 1}]^W
\end{equation}
for some constants $c_{\nu^\prime}\in\mathbb{C}$.

We fix now
arbitrary $\mu\in\Lambda_n^+$ and we prove that $p_\mu$ is the monic
Macdonald-Koornwinder polynomial of degree $\mu$ with suitable specialization
of the parameters. Formula \eqref{eig1} applied to $p_\mu$
shows that
\begin{equation}\label{eig1a}
p_\mu=
m_{\mu}+\sum_{\stackrel{\mu^\prime\in\Lambda_n^+}{\mu^\prime<\mu}}
c_{\mu^\prime} m_{\mu^\prime}\in\mathbb{C}[u^{\pm 1}]^W
\end{equation}
for some constants $c_{\mu^\prime}\in\mathbb{C}$.
By the generalized quantum Schur orthogonality relations
we furthermore have
\begin{equation}\label{eig2}
\begin{split}
\int_T p_\mu(u)p_\nu^*(u)W(u)\frac{du}{u}=0,\qquad \forall\,
\nu\in\Lambda_n^+:\,\nu\not=\mu
\end{split}
\end{equation}
with weight function $W(u)=f_0|_T(u)(f_0|_T)^*(u)\Delta_{\sigma,\tau}(u)$.
By Theorem \ref{maintheorem}{\bf (ii)}, the weight function $W(u)$
($u\in T$) can
be expressed in terms of
the Macdonald-Koornwinder weight function $\Delta(u;a,b,c,d;q,t)$ as
\[W(u)=C\,
\Delta\bigl(u;-q^{\sigma+\tau+1+\kappa_1+\kappa_2},-q^{-\sigma-\tau+1},
q^{\sigma-\tau+1},q^{-\sigma+\tau+1+\kappa_1-\kappa_2};q^2,
q^{2\kappa+2}\bigr)
\]
for some nonzero constant $C$. Properties \eqref{eig1}
and \eqref{eig2} thus imply
\begin{equation}\label{eig3}
\int_Tp_\mu(u)m_{\nu}(u)W(u)\frac{du}{u}=0\qquad \forall\,
\nu\in\Lambda_n^+:\,\nu<\mu.
\end{equation}
The proof now follows from the fact that
the Macdonald-Koornwinder polynomial
\[P_\mu\bigl(u;-q^{\sigma+\tau+1+\kappa_1+\kappa_2},
-q^{-\sigma-\tau+1},
q^{\sigma-\tau+1},q^{-\sigma+\tau+1+\kappa_1-\kappa_2};q^2,
q^{2\kappa+2}\bigr)\in\mathbb{C}[u^{\pm 1}]^W
\]
has been defined as the unique $W$-invariant Laurent polynomial satisfying
the properties \eqref{eig1a} and \eqref{eig3}.
\end{proof}


\end{document}